\documentclass[a4paper]{amsart}
\usepackage{graphicx,amsmath, amssymb, overpic}

\usepackage{color}
\input xy

\xyoption{all}

\newcommand{\ad}{\ensuremath{\operatorname{ad}}}
\newcommand{\Aut}{\ensuremath{\operatorname{Aut}}}

\newcommand{\Hom}{\ensuremath{\operatorname{Hom}}}
\newcommand{\Stab}{\ensuremath{\operatorname{Stab}}}

\newcommand{\dist}{\ensuremath{\operatorname{dist}}}
\newcommand{\Ball}{\ensuremath{\operatorname{Ball}}}
\newcommand{\Inn}{\ensuremath{\operatorname{Inn}}}

\newcommand{\F}{\mathbb{F}}
\newcommand{\N}{\mathbb{N}}

\newcommand{\G}{\Gamma}
\newcommand{\bs}{\backslash}
\newcommand{\forget}[1]{}

\newcommand{\Fqt}{\F_q(\!(t^{-1})\!)}
\newcommand{\bX}{\partial X}
\newcommand{\ep}{\varepsilon}

\newcommand{\cA}{\mathcal{A}}

\newcommand{\cU}{\mathcal{U}}
\newcommand{\la}{\langle}
\newcommand{\ra}{\rangle}

\newtheorem{theorem}{Theorem}[section]
\newtheorem{prop}[theorem]{Proposition}
\newtheorem{lemma}[theorem]{Lemma}
\newtheorem{corollary}[theorem]{Corollary}
\newtheorem{conjecture}{Conjecture}
\newtheorem{subcase}{Subcase}

\theoremstyle{definition}

\newtheorem*{Remark}{Remark}
\newtheorem{defn}{Definition}

\newtheorem{example}{Example}

\title{Lattices in complete rank 2 Kac--Moody groups} 

\author{Inna (Korchagina) Capdeboscq} \address{Mathematics Institute, Zeeman Building, University of Warwick, Coventry
CV4 7AL, UK} \email{I.Korchagina@warwick.ac.uk}\thanks{The work of the first author was partially supported by an EPSRC Grant No.
EP/EO49508 and partially by NSF Grant No. DMS-0635607.  This work of the second author was supported in part by NSF Grant No. DMS-0805206 and in part by EPSRC Grant No. EP/D073626/2, and the second author is now supported by ARC Grant No. DP110100440.}

\author{Anne Thomas}\address{School of Mathematics and Statistics, Carslaw Building F07, University of Sydney NSW 2006, Australia}\email{anne.thomas@sydney.edu.au}

\date{\today}

\begin{document}

\begin{abstract}  Let $\Lambda$ be a minimal Kac--Moody group of rank $2$ defined over the  finite field $\F_q$, where $q = p^a$ with $p$ prime.  Let $G$ be the topological Kac--Moody group obtained by completing $\Lambda$.  An example is $G=SL_2(K)$, where $K$ is the field of formal Laurent series over $\F_q$.  The group $G$ acts on its Bruhat--Tits building $X$, a tree, with quotient a single edge.  We construct new examples of cocompact lattices in $G$, many of them edge-transitive.  We then show that if cocompact lattices in $G$ do not contain $p$--elements,  the lattices we construct are the only edge-transitive lattices in $G$, and that our constructions include the cocompact lattice of minimal covolume in $G$.  We also observe that, with an additional assumption on $p$--elements in $G$, the arguments of Lubotzky \cite{L} for the case $G = SL_2(K)$ may be generalised to show that there is a positive lower bound on the covolumes of all lattices in $G$, and that this minimum is realised by a non-cocompact lattice, a maximal parabolic subgroup of $\Lambda$. \end{abstract}

\maketitle

\section{Introduction}\label{s:intro}

A classical theorem of Siegel~\cite{Siegel45} states that the minimum covolume among lattices in $G=SL_2(\mathbb{R})$ is~$\frac{\pi}{21}$, and determines the lattice which realises this minimum.  In the nonarchimedean setting, Lubotzky~\cite{L} determined the lattice of minimal covolume in $G=SL_2(K)$, where $K$ is the field $\Fqt$ of formal Laurent series over $\mathbb{F}_q$.

The group $G=SL_2(\Fqt)$ has, in recent developments, been viewed as the first example of a complete Kac--Moody group of rank $2$ over the finite field $\F_q$.  By definition, a complete Kac--Moody group is the completion of a minimal Kac--Moody group $\Lambda$ over a finite field, with respect to some topology.  We use the completion in the ``building topology", as discussed in, for example,~\cite{CR}.  Complete Kac--Moody groups are locally compact, totally disconnected topological groups, which may be thought of as infinite-dimensional analogues of semisimple algebraic groups (see Section~\ref{s:kac-moody} below for details).  

\subsection{Constructions of cocompact lattices}

Our first main result, Theorem \ref{t:cocompact existence} below, constructs many new cocompact lattices in rank $2$ complete Kac--Moody groups $G$.  It is interesting that there exist any cocompact lattices in such groups $G$, since for $n \geq 3$, Kac--Moody groups of rank $n$ do not admit cocompact lattices (with the possible exception of those whose root systems contain a subsystem of type $\tilde{A}_n$ --- see~\cite[Remark 4.4]{CM}).  In rank $2$, the previous examples of cocompact lattices in $G$ non-affine that are known to us are the free Schottky groups constructed by Carbone--Garland~\cite[Section 11]{CG}, some of the lattices constructed by R\'emy--Ronan \cite[Section 4.B]{RR}, which in rank $2$ are free products of finite cyclic groups, and some of the lattices obtained as centralisers of certain involutions by Gramlich--Horn--M\"uhlherr \cite[Section 7.3]{GHM}.
  
As we recall in Section \ref{s:kac-moody} below, the Kac--Moody groups $G$ that we consider have Bruhat--Tits building a regular tree $X$, and the action of $G$ on $X$ induces an edge of groups

\begin{figure}[ht]
\begin{picture}(20,10)(-20,0)
\put(-65,5){$\mathbb{G} = $}
\put(-30,4){\line(1,0){80}}
\put(-33,10){$P_1$}
\put(47,10){$P_2$}
\put(8,-5){$B$}
\put(-30,4){\circle*{3.5}}
\put(50,4){\circle*{3.5}}
\end{picture}
\end{figure}

\noindent where $P_1$ and $P_2$ are the standard parahoric subgroups of $G$ and $B=P_1 \cap P_2$ is the standard Iwahori subgroup.  The kernel of the $G$--action on $X$ is the finite group $Z(G)$, the centre of $G$ (see~\cite{CR}).  Now let $\G$ be a cocompact lattice in $G$ which acts transitively on the edges of $X$.  Then as we recall in Section \ref{s:lattices}, $\G$ is  the fundamental group $A_1 *_{A_0} A_2$ of an edge of groups 

\begin{figure}[ht]
\begin{picture}(20,10)(-20,0)
\put(-65,5){$\mathbb{A} = $}
\put(-30,4){\line(1,0){80}}
\put(-33,10){$A_1$}
\put(47,10){$A_2$}
\put(8,-5){$A_0$}
\put(-30,4){\circle*{3.5}}
\put(50,4){\circle*{3.5}}
\end{picture}
\end{figure}

\noindent with $A_0$, $A_1$ and $A_2$ finite groups (see Section \ref{s:bass-serre} for background on graphs of groups).  
 
We now state Theorem \ref{t:cocompact existence}, which constructs new cocompact lattices in $G$, most of them edge-transitive.  We discuss our assumption that $G$ has symmetric generalised Cartan matrix $A$ after the statement of Theorem \ref{t:cocompact existence}.  There are some exceptional constructions for small values of $q$ which we then record in Theorem~\ref{t:exceptions}.  Apart from the affine case $G = SL_2(\Fqt)$, there are no known linear representations of the groups $G$ in Theorems \ref{t:cocompact existence} and \ref{t:exceptions}.  

Our notation is as follows.  We write $C_n$ for the cyclic group of order $n$ and $S_n$ for the symmetric group on $n$ letters. Since for a finite field $\F_q$ and the root system $A_1$ there exist at most two corresponding finite groups of Lie type, namely $SL_2(\F_q)$ and $PSL_2(\F_q)$,  to avoid complications we  use Lie-theoretic notation and write $A_1(q)$ in both cases. (We will discuss this ambiguity whenever  necessary.)  We denote by $T$ a fixed maximal split torus of $G$ with $T\leq P_1\cap P_2$. Then $Z(G) \leq T$, and $T$ is isomorphic to a quotient of $\F_q^*\times \F_q^*$ (the particular quotient depending upon $G$).  Finally, the parahoric subgroups $P_1$ and $P_2$ of $G$ admit Levi decompositions (see Section \ref{s:rank 2} and, in particular, Proposition~\ref{p:completion} below), and for $i = 1,2$ we denote by $L_i$ a Levi complement of $P_i$.  The group $L_i$ factors as $L_i=M_iT$, where $M_i\cong A_1(q)$ is normalised by $T$, and we denote by $H_i$ a non-split torus of $M_i$ such that $N_T(H_i)$ is as big as possible.   We say that two edge-transitive cocompact lattices $\G=A_1 *_{A_0} A_2$ and $\G'=A_1' *_{A_0'} A_2'$ in $G$ are \emph{isomorphic} if $A_i \cong A_i'$ for $i = 0,1,2$ and the obvious diagram commutes.   

\begin{theorem}\label{t:cocompact existence}  Let $G$ be a complete Kac--Moody group of rank $2$ with symmetric generalised Cartan matrix $\left(\begin{matrix}2 & -m\\-m & 2 \end{matrix}\right)$, $m \geq 2$, defined over the finite field $\F_q$ of order $q=p^a$ where $p$ is prime.  Then in the following cases, the group $G$ admits edge-transitive cocompact lattices of each of the following isomorphism types.

\begin{enumerate} 
 \item \label{i: p = 2} If $p=2$ then $\Gamma=A_1*_{A_0}A_2$ where for $i=1,2$,  $A_i=A_0 \times H_i$ and $H_i \cong C_{q+1}$, and $A_0$ is any subgroup of $Z(G)$. 
 \item \label{i:PSL} If $p$ is odd and $L_i/Z(L_i)\cong PSL_2(q)$, assume also that $q\equiv 3 \pmod 4$.  Then $\Gamma=A_1*_{A_0}A_2$ where for  $i=1,2$, $A_i=A_0N_{M_i}(H_i)$, and $A_0$ is a subgroup of $Z(G)$ with $A_0\cap N_{M_i}(H_i)=Z(M_i)$, $i=1,2$.
\item \label{i:PGL} If $p$ is odd and $L_i/Z(L_i)\cong PGL_2(q)$, let \[ C_i:=C_{L_i}(H_i), \quad Q_i'\in\mathcal{S}yl_2(C_i),  \quad Q_i := \left( Q_i'\cap Z(L_i) \right) \in\mathcal{S}yl_2(Z(L_i))\] and let $Q_i^2$ be the unique subgroup of $Q_i$ of index $2$.
\begin{enumerate}
\item  If $q\equiv 1 \pmod 4$ and $Q_i^2\leq Z(G)$ then:
\begin{enumerate}
\item\label{i:PGL q = 1}  $\Gamma=A_1*_{A_0}A_2$ where for  $i=1,2$, $A_i=H_iQ_i'\langle t_i\rangle Z_0$ where  
$ t_i \in N_T(H_i)-C_T(H_i)$ is of order $2$, $Z_0\leq Z(G)$ and $A_0=Q_i\langle t_i\rangle Z_0$. 
\item\label{i:PGL q = 1 other}  Moreover, if $Q_i\leq Z(G)$, then also 
$\Gamma=A_1*_{A_0}A_2$ where for  $i=1,2$, $A_i=H_iQ_i'Z_0$ where  $Z_0\leq Z(G)$ and $A_0=Q_iZ_0\leq Z(G)$.
\end{enumerate}
\item\label{i:q = 3 mod 4} If  $q\equiv 3 \pmod 4$:
\begin{enumerate}
\item\label{i:PGL q = 3}  $\Gamma=A_1*_{A_0}A_2$ where for $i = 1,2$, $A_i=H_iQ_i'T_0Z_0$ with  $T_0\in\mathcal{S}yl_2(T)$, $Q_i'\cap T_0=Q_i$,  $Z_0\leq Z(G)$ and $A_0=T_0Z_0$.
\item\label{i:PGL q = 3 other} If $Z(M_i)\leq Z(G)$, then  also $\Gamma=A_1*_{A_0}A_2$ where either for $i = 1,2$, $A_i=H_iQ_i'A_0$ with  $Q'_i\cap A_0=Z(M_i)$ and $A_0\leq Z(G)$, or the description of $\Gamma$ in \eqref{i:PSL}  holds.
\end{enumerate}
\end{enumerate} 
\end{enumerate}

In all other cases, $G$ admits a cocompact lattice $\G'$ which acts on the tree $X$ inducing a graph of groups of the form
\begin{center}
\begin{overpic}{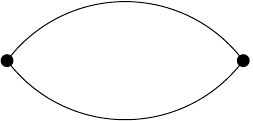}
\put(-70,21){$\mathbb{A}'=$}
\put(-10,21){$N_1$}
\put(100,21){$N_2$}
\put(47,3){$S$}
\put(47,49){$S$}
\end{overpic}
\end{center} where the finite groups $S$, $N_1$ and $N_2$ will be defined in Section~\ref{ss:define groups} below.  

\end{theorem}

\begin{Remark}\label{r:symmetric}  If $G$ has generalised Cartan matrix $A = \left(\begin{matrix}2 & -m\\-n & 2 \end{matrix}\right)$ which is not symmetric but satisfies $m, n \geq 2$, then the constructions in Theorem \ref{t:cocompact existence}\eqref{i: p = 2} and \eqref{i:PSL} above still hold.  In case \eqref{i:PGL} of Theorem \ref{t:cocompact existence}, there are analogous edge-transitive lattices with an even more involved description, which we omit. \end{Remark}

Our exceptional constructions for small values of $q$ are as follows.

\begin{theorem}\label{t:exceptions} Let $G$ be as in Theorem \ref{t:cocompact existence} above.  Then in the following cases, $G$ admits cocompact edge-transitive lattices of the following isomorphism types.  

When $p$ is odd and $q \equiv 1 \pmod 4$: 
\begin{enumerate} 
\item If $q=5$,  $\Gamma=A_1*_{A_0}A_2$ where for $i=1,2$, $A_i=A_0N_i$ where $N_i\cong A_1(3)$, $A_0\leq N_T(N_i)$ and $|N_i: N_i\cap A_0|=6$; and 
\item If $q=29$, $\Gamma=A_1*_{A_0}A_2$ where for $i=1,2$,   $A_i=A_0N_i$ where $N_i\cong A_1(5)$, $A_0\leq N_T(N_i)$ and $|N_i: N_i\cap A_0|=30$. 
\end{enumerate} 

When $p$ is odd and $q \equiv 3 \pmod 4$:
\begin{enumerate} 
\item If $q=7$ or $23$,   $\Gamma=A_1*_{A_0}A_2$ where for $i=1,2$, $A_i=A_0N_i$, $N_i\cong S_4\ \mbox{or}\ 2S_4$,  $A_0\leq N_T(N_i)$  and $|N_i:N_i\cap A_0|=q+1$ where $N_i\cap A_0$ is cyclic.
\item If $q=11$,   $\Gamma=A_1*_{A_0}A_2$  where for $i=1,2$, $A_i=A_0N_i$, and $A_0\leq N_T(N_i)$  with $|N_i:N_i\cap A_0|=12$, $N_i\cap A_0$ being cyclic, and one of the following holds:
\begin{enumerate}
\item $N_1\cong N_2\cong A_1(3)$; or
\item $N_1\cong N_2\cong  A_1(5)$.
\end{enumerate}
\item  If $q=19$ or $59$,   $\Gamma=A_1*_{A_0}A_2$ where  for $i=1,2$, $A_i=A_0N_i$, $N_i\cong  A_1(5)$, $A_0\leq N_T(N_i)$  and $|N_i:N_i\cap A_0|=q+1$ with $N_i\cap A_0$ being cyclic.
\end{enumerate}
\end{theorem}

We prove Theorems~\ref{t:cocompact existence} and~\ref{t:exceptions} in Section~\ref{s:cocompact}.  By the general theory of tree lattices (see Section~\ref{s:lattices}), each $\G$ or $\G'$ appearing in Theorems \ref{t:cocompact existence} and \ref{t:exceptions} yields a cocompact lattice in the full automorphism group $\Aut(X)$ of the tree $X$.  Since $\Aut(X)$ is much larger than the Kac--Moody group $G$, the key point in proving Theorems \ref{t:cocompact existence} and \ref{t:exceptions} is to show that each such $\G$ or $\G'$ embeds into $G$ as a cocompact lattice.  

For this, we develop an embedding criterion, Proposition \ref{p:embedding} in Section~\ref{s:embedding} below, which may be applied to construct lattices in any closed locally compact group acting transitively on the edges of a locally finite biregular tree.  Our embedding criterion generalises \cite[Lemma 3.1]{L}, which was used in \cite{L} to construct edge-transitive lattices in the affine case $G = SL_2(\Fqt)$.  We are able to provide somewhat simpler proofs even in that case by using Bridson and Haefliger's covering theory for complexes of groups \cite{BH} instead of Bass' covering theory for graphs of groups \cite{B}, since the theory for complexes of groups has a less complicated notion of morphism.   

Our construction of the lattice $\G'$ in Theorem~\ref{t:cocompact existence} generalises the construction in Example (6.2) of Lubotzky--Weigel~\cite{LW} of a cocompact lattice in $G = SL_2(\Fqt)$ when $q \equiv 1 \pmod 4$.  Example (6.2) of \cite{LW} uses the embedding criterion \cite[Theorem 5.4]{LW}, which relies upon Bass' covering theory \cite{B} and is specific to the affine case, while we apply our criterion Proposition~\ref{p:embedding} instead.  

\subsection{Conjectures}

Let $G$ be as in Theorem \ref{t:cocompact existence}.  Our further results, Theorems \ref{t:classification}, \ref{t:covolumes} and \ref{t:covolumes intro} below, depend upon conjectures about the behaviour of $p$--elements in $G$.  

In the affine case $G = SL_2(\Fqt)$, an element of $G$ has order $p$ if and only if it is unipotent.  Hence by Godement's Criterion, no $p$--element is contained in a cocompact lattice.  Moreover, any unipotent element of $G = SL_2(\Fqt)$ is contained in a conjugate of the upper unitriangular subgroup of $G$, which is an infinite group isomorphic to $\Fqt$.  On the other hand, Lubotzky \cite{L} showed that no cocompact lattice of $SL_2(\Fqt)$ contains non-trivial $p$--elements (analogous to the classical Godement's cocompactness criterion).  For general $G$, as we explain in Section \ref{s:motivation} below, there are many $p$--elements which cannot be contained in any cocompact lattice $\G < G$, since they belong to a conjugate of a canonical  subgroup $\cU$ of $G$ (see Section \ref{s:action ends} for the definition of $\cU$).  We thus make the following conjectures in the general case.

\begin{conjecture}\label{c:cocompact} Cocompact lattices in $G$ do not contain $p$--elements.
\end{conjecture}

\begin{conjecture}\label{c:unique end} Any $p$--element in $G$ is contained in a conjugate of the subgroup $\cU$ of $G$.
\end{conjecture}

\noindent In Section \ref{s:motivation}, we also explain why Conjecture \ref{c:unique end} implies Conjecture \ref{c:cocompact}.   

Theorems \ref{t:classification} and \ref{t:covolumes} below only require Conjecture \ref{c:cocompact}, while Conjecture \ref{c:unique end} is needed for Theorem \ref{t:covolumes intro}.   As in the remark following Theorem \ref{t:cocompact existence} above, similar results to Theorems \ref{t:classification}, \ref{t:covolumes} and \ref{t:covolumes intro} will hold with Cartan matrix $A =  \left(\begin{matrix}2 & -m\\-n & 2 \end{matrix}\right) $, $m, n \geq 2$, but we omit these more cumbersome statements.

\subsection{Classification of edge-transitive lattices}

Assuming Conjecture \ref{c:cocompact}, we are able to classify the edge-transitive lattices in $G$ up to isomorphism, as follows.

\begin{theorem}\label{t:classification} Let $G$ be as in Theorem \ref{t:cocompact existence}, and assume that Conjecture \ref{c:cocompact} holds.  Then every edge-transitive lattice in $G$ is isomorphic to one of the amalgamated free products of finite groups $A_1 *_{A_0} A_2$ in Theorems \ref{t:cocompact existence} and \ref{t:exceptions}.\end{theorem}

\noindent In other words, every edge-transitive lattice in $G$ which does not contain $p$--elements appears in Theorems \ref{t:cocompact existence} or \ref{t:exceptions}.  

The question of classifying amalgams of finite groups is, in general, difficult.  An \emph{$(m,n)$--amalgam} is an amalgamated free product $A_1 *_{A_0} A_2$ where $A_0$ has index $m$ in $A_1$ and index $n$ in $A_2$.  To avoid trivial examples, such an amalgam is said to be \emph{faithful} if $A_0$, $A_1$ and $A_2$ have no common normal subgroup.  A deep theorem of Goldschmidt~\cite{G} established that there are only 15 faithful $(3,3)$--amalgams of finite groups, and classified such amalgams.  Goldschmidt and Sims conjectured that when both $m$ and $n$ are prime, there are only finitely many faithful $(m,n)$--amalgams of finite groups (see~\cite{BK,F,G}).  This conjecture remains open, except for the case $(m,n)=(2,3)$, which was established by Djokovi\'c--Miller~\cite{DM}, and the work of  Fan~\cite{F}, who proved the conjecture when the edge group $A_0$ is a $p$--group, with $p$ a prime distinct from both $m$ and $n$.  

On the other hand, Bass--Kulkarni~\cite{BK} showed that if either $m$ or $n$ is composite, there are infinitely many faithful $(m,n)$--amalgams of finite groups.  The constructions in \cite{BK} may be viewed as giving infinitely many non-isomorphic edge-transitive lattices in the automorphism group of an $(m,n)$--biregular tree.  Thus if Conjecture \ref{c:cocompact} holds, there are in contrast only finitely many edge-transitive lattices in Kac--Moody groups $G$ as in Theorem \ref{t:cocompact existence}.  We note that, since the action of such $G$ on $X$ is not in general faithful, an amalgam $\G$ may embed as an edge-transitive lattice in $G$ even though it is not faithful (but its kernel will be at most the finite group $Z(G)$).

We prove Theorem \ref{t:classification} in Section \ref{s:classification}, making careful use of classical results on the subgroups of $SL_2(q)$, $PSL_2(q)$ and $PGL_2(q)$ of order coprime to $p$ and their actions on the projective line $\mathbb{P}^1(q)$.

\subsection{Covolumes of cocompact lattices}

If Conjecture \ref{c:cocompact} holds, there are also important consequences for the covolumes of lattices in $G$.   We recall in Section~\ref{s:lattices} below that the Haar measure $\mu$ on $G$ may be normalised so that the covolume $\mu(\G \bs G)$ of a cocompact lattice $\G$ with quotient graph of groups having two vertex groups $A_1$ and $A_2$ is equal to $|A_1|^{-1} + |A_2|^{-1}$.  Using this normalisation, we obtain the following.

\begin{theorem}\label{t:covolumes} Let $G$ be as in Theorem \ref{t:cocompact existence} above and assume that Conjecture \ref{c:cocompact} holds.  Then for $q \geq 514$ \[ \min \{ \mu(\G \bs G) \mid \G\mbox{ a cocompact lattice in }G\} = \frac{2}{(q+1)|Z(G)|\delta}\] where $\delta \in \{ 1, 2, 4\}$ (depending upon the particular group $G$).  Moreover, we construct a cocompact lattice which realises this minimum, and this lattice appears in Theorem~\ref{t:cocompact existence}. \end{theorem}

\noindent In other words, for $q$ large enough, among cocompact lattices in $G$ which do not contain $p$--elements the minimal covolume is $2/(q+1)|Z(G)|\delta$, and this minimum is achieved by one of the lattices in Theorem \ref{t:cocompact existence}.

Theorem \ref{t:covolumes} generalises Theorem 2 of Lubotzky \cite{L} and part of Proposition C  of Lubotzky--Weigel \cite{LW}, which together determine the cocompact lattices of minimal covolume in the affine case $G = SL_2(\Fqt)$.

We prove Theorem \ref{t:covolumes} in Section \ref{s:covolumes cocompact}.  The proof is very delicate and has many cases, involving consideration of how various subgroups of $SL_2(q)$, $PSL_2(q)$ and $PGL_2(q)$ of order coprime to $p$ can act on $X$.

\subsection{Covolumes of all lattices}

Finally, assuming the stronger Conjecture \ref{c:unique end}, we obtain a lower bound on the covolume of all lattices in $G$.  As we recall in Section \ref{s:kac-moody}, the minimal Kac--Moody group $\Lambda$ has twin buildings $X_+ \cong X_-$, and there are two completions $G_+ \cong G_-$ into which $\Lambda$ embeds, with $G_\pm$ acting on $X_\pm$.  As shown independently in Carbone--Garland \cite{CG} and R\'emy \cite{R}, a negative maximal parabolic subgroup $P_-$ of $\Lambda$ is a non-cocompact lattice in $G = G_+$.   We have:

\begin{theorem}\label{t:covolumes intro} Let $G$ be as in Theorem \ref{t:cocompact existence} above and assume that Conjecture \ref{c:unique end} above holds.   Then \[ \min \{ \, \mu(\G \bs G) \mid \G\mbox{ a lattice in }G \, \} = \frac{2}{(q+1)(q-1)|T|}.\] Moreover, this minimum is realised by the non-cocompact lattice $P_-$. \end{theorem}

Theorem \ref{t:covolumes intro}  generalises Theorem 1 of Lubotzky \cite{L}, which establishes that the lattice of minimal covolume in $G = SL_2(\Fqt)$ is the maximal  parabolic subgroup $SL_2(\F_q[t])$, a non-cocompact lattice.  More recently, Golsefidy \cite{G} has shown that for many Chevalley groups $\mathcal{G}$, the lattice of minimal covolume in $\mathcal{G}(\Fqt)$ is the non-cocompact lattice $\mathcal{G}(\F_q[t])$.  Theorem \ref{t:covolumes intro} and the results of \cite{L} and \cite{G} thus contrast with Siegel's original result~\cite{Siegel45} that the lattice of minimal covolume in $SL_2(\mathbb{R})$ is cocompact.  

The proof of Theorem \ref{t:covolumes intro} is similar to that of \cite[Theorem 1]{L}, after replacing the term ``unipotent element" with ``$p$--element".  We thus omit the proof.

\subsection*{Acknowledgements}

We are grateful to  Benson Farb, Richard Lyons, Chris Parker and Kevin Wortman for helpful conversations.
We would also like to thank Pierre-Emmanuel Caprace, Bertrand R\'emy, Ron Solomon
and an anonymous referee of an earlier version of this paper for suggestions and comments that improved this article.  
This work was completed  during the first author's  stay at the  Institute for Advanced Study and she would like to gratefully acknowledge the wonderful year she spent there.
The second author thanks the Mathematical Sciences Research Institute and the London Mathematical Society for travel support.
 
This material is based upon work supported by the National Science Foundation under  agreement No. DMS-0635607. Any opinions, findings and conclusions or recommendations expressed in this material are those of the authors and do not necessarily reflect the views of the National Science Foundation.

\section{Preliminaries}\label{s:preliminaries}

In Section \ref{s:trees} we recall basic definitions and establish notation for graphs and trees.   Section \ref{s:bass-serre} sketches  the theory of graphs of groups, including covering theory.  In Section \ref{s:lattices} we recall some theory of tree lattices.  Section \ref{s:kac-moody} then presents background on the Kac--Moody groups $G$ that we consider, and in Section \ref{s:finite groups} we recall some classical theorems concerning the finite subgroups of $SL_2(q)$, $PSL_2(q)$ and $PGL_2(q)$.

\subsection{Graphs and trees}\label{s:trees}

Let $A$ be a connected graph with sets $VA$ of vertices and $EA$ of oriented edges. The initial and terminal vertices of  $e \in EA$ are denoted by $\partial_0e$ and $\partial_1e$ respectively. The map $e \mapsto \overline{e}$ is orientation reversal, with $\overline{\overline{e}} = e$ and $\partial_{1-j}\overline{e} = \partial_j e$ for $j = 0,1$.  Given a vertex $a \in VA$, we denote by $E_A(a)$ the set of edges \[ E_A(a):=\{ e \in EA \mid \partial_0e = a \}\] with initial vertex $a$.

Let $A$ and $B$ be graphs.  A \emph{morphism of graphs} is a function $\theta:A \to B$ taking vertices to vertices and edges to edges, such that for every edge $e \in EA$, $\theta(\overline e) = \overline{\theta(e)}$ and $\theta(\partial_j(e)) = \partial_j(\theta(e))$ for $j = 0,1$.

Let $X$ be a simplicial tree with vertex set $VX$ and edge set $EX$.  A group $G$ is said to act \emph{without inversions} on $X$ if for all $g \in G$ and all edges $e \in EX$, if $g$ preserves the edge $e$ then $g$ fixes $e$ pointwise.

\begin{prop}[Serre, Proposition 19, Section I.4.3~\cite{S}]\label{p:finite_gp_tree} Let $A$ be a finite group acting without inversions on a simplicial tree $X$.  Then there is a vertex of $X$ which is fixed by $A$. \end{prop}

Now equip $X$ with a length metric $d$ by declaring each edge of $X$ to have length $1$.   Given an edge $e$ of $X$ and an integer $n \geq 0$, we define $\Ball(e,n)$ to be the union of the closed edges in $X$ all of whose points are distance at most $n$ from a point in the closed edge $e$.  By abuse of notation, we then define the distance $d(e,e')$ between edges $e$ and $e'$ of $X$ to be $0$ if $e = e'$, and to be $n \geq 1$ if $e'$ is in $\Ball(e,n)$ but not $\Ball(e,n - 1)$.

A \emph{geodesic line} in $X$ is an isometry $l:(-\infty,\infty) \to X$, and a \emph{geodesic ray} in $X$ is an isometry $r:[t,\infty) \to X$, $t \in \mathbb{R}$, such that $r(t)$ is a vertex of $X$; in this case we say that the ray \emph{begins at} the vertex $r(t)$.  We will often identify a geodesic line or ray with its image in $X$.  We say that two geodesic rays $r$ and $r'$ in $X$ are \emph{equivalent} if their intersection is infinite.  The \emph{boundary} $\bX$ of $X$, which is the same thing as the \emph{set of ends} of $X$, is the collection of equivalence classes of geodesic rays.  Given a geodesic ray $r$, we say that an end $\ep$ is \emph{determined by $r$} if $r$ belongs to the equivalence class $\ep$.

If $G$ is a group of isometries of $X$ which acts without inversions, then an element $g \in G$ is either \emph{elliptic}, meaning that $g$ fixes at least one vertex of $X$, or \emph{hyperbolic}, meaning that $g$ does not fix any vertex and acts as a translation along its \emph{axis}, a geodesic line in $X$ (see, for example, \cite[Chapter II.6]{BH}).  If $g$ is hyperbolic then $g$ generates an infinite cyclic subgroup of $G$.  

Let $g$ be a hyperbolic isometry of $X$, with axis $l$.  Then $g$ has exactly two fixed points in $\bX$, which we denote by $l(-\infty)$ and $l(\infty)$.  One of these fixed points is repelling and the other is attracting.

\subsection{Bass--Serre theory}\label{s:bass-serre}

A \emph{graph of groups} $\mathbb{A}=(A,\mathcal{A})$ over a connected graph $A$ consists of an assignment of vertex groups $\mathcal{A}_a$ for each $a \in VA$ and edge groups $\mathcal{A}_e = \mathcal{A}_{\overline{e}}$ for each $e \in EA$, together with monomorphisms $\alpha_e:\mathcal{A}_e \to \mathcal{A}_{\partial_0e}$ for each $e \in EA$.   

Any action of a group $\G$ on a tree $X$ without inversions induces a graph of groups over the quotient graph $A = \G \bs X$.  See for example~\cite{B} for the definitions of the \emph{fundamental group} $\pi_1(\mathbb{A},a_0)$ and the \emph{universal cover} $X = \widetilde{(A,a_0)}$ of a graph of groups $\mathbb{A}=(A,\mathcal{A})$, with respect to a basepoint $a_0 \in VA$.  The universal cover $X$ is a tree, on which $\pi_1(\mathbb{A},a_0)$ acts by isometries inducing a graph of groups isomorphic to $\mathbb{A}$.

In the special case that $\mathbb{A}$ is a graph of groups over an underlying graph $A$ which is a single edge $e$, we say that $\mathbb{A}$ is an \emph{edge of groups}.  Suppose $\partial_0e=a_1$ and $\partial_1e=a_2$.  Write $A_0$ for the edge group $\mathcal{A}_e$, and for $i = 1,2$ let $A_i$ be the vertex group $\mathcal{A}_{a_i}$.  The fundamental group $\pi_1(\mathbb{A},a_1)$ is then isomorphic to the free product with amalgamation $A_1 *_{A_0} A_2$, and the universal cover $X = \widetilde{(A,a_1)}$ is an $(m,n)$--biregular tree, where $m = [A_1:A_0]$ and $n = [A_2:A_0]$.  

We now adapt definitions from covering theory for complexes of groups (see \cite[Chapter III.$\mathcal{C}$]{BH}) to graphs of groups, and recall a necessary result from covering theory.  For the precise relationship between the category of graphs of groups and the category of complexes of groups over $1$--dimensional spaces, see~\cite[Proposition~2.1]{Th}.

\begin{defn}[Morphism of graphs of groups]\label{d:morphism} Let $\mathbb{A} = (A,\mathcal{A})$ and $\mathbb{B} = (B,\mathcal{B})$ be graphs of groups, with monomorphisms from edge groups to vertex groups respectively $\alpha_e:\mathcal{A}_e \to \mathcal{A}_{\partial_0e}$ for $e \in EA$ and $\beta_e:\mathcal{B}_f \to \mathcal{B}_{\partial_0f}$ for $f \in EB$.  Let $\theta:A \to B$ be a morphism of graphs.  A \emph{morphism of graphs of groups} $\Phi:\mathbb{A} \to \mathbb{B}$ over $\theta$ is given by: \begin{enumerate} \item a homomorphism $\phi_x:\mathcal{A}_x \to \mathcal{B}_{\theta(x)}$ of groups, for every $x \in VA \cup EA$; and \item an element $\phi(e) \in \mathcal{B}_{\partial_0(\theta(e))}$ for each $e \in EA$ such that the following diagram commutes, where $a = \partial_0 e$: \[ \xymatrix{ \cA_e    \ar[d]^-{\alpha_e}   \ar[r]_{\phi_e} & \mathcal{B}_{\theta(e)} \ar[d]^-{\ad(\phi(e)) \circ \beta_{\theta(e)} }\\ \cA_{a}   \ar[r]_{\phi_a} & \mathcal{B}_{\theta(a)} } \] \end{enumerate} \end{defn}

\begin{defn}[Covering of graphs of groups]\label{d:covering} With notation as in Definition~\ref{d:morphism} above, $\Phi:\mathbb{A} \to \mathbb{B}$ is a \emph{covering of graphs of groups} if in addition: \begin{enumerate} \item for each $x \in VA \cup EA$ the homomorphism $\phi_x:\mathcal{A}_x \to \mathcal{B}_{\theta(x)}$ is injective; and \item for each edge $f \in EB$ and each vertex $a \in VA$ with $\partial_0 f = b = \theta(a)$, the map \[ \Phi_{a/f}: \coprod_{e \in E_A(a) \cap \theta^{-1}(f)} \mathcal{A}_a/\alpha_e(\mathcal{A}_e) \to \mathcal{B}_b/\beta_f(\mathcal{B}_f) \] induced by $g \mapsto \phi_a(g)\phi(e)$ is bijective. \end{enumerate} \end{defn}

The result from covering theory that we will need is:

\begin{prop}[Bass, Proposition 2.7 of~\cite{B}]\label{p:coverings} Let $\mathbb{A}=(A,\mathcal{A})$ and $\mathbb{B}=(B,\mathcal{B})$ be graphs of groups.  Choose basepoints $a_0 \in A$ and $b_0 \in B$.  If there is a covering of graphs of groups $\Phi:\mathbb{A} \to \mathbb{B}$ over $\theta:A \to B$ with $\theta(a_0) = b_0$, then $\pi_1(\mathbb{A},a_0)$ injects into $\pi_1(\mathbb{B},b_0)$.\end{prop}

\subsection{Lattices in groups acting on trees}\label{s:lattices}

Let $G$ be a locally compact topological group with left-invariant Haar measure $\mu$.  Recall that a discrete subgroup $\Gamma \leq G$ is a \emph{lattice} if $\Gamma\backslash G$ carries a finite $G$--invariant measure, and is \emph{cocompact} if $\G \bs G$ is compact.

Now let $X$ be a locally finite tree and let $G$ be a closed, cocompact group of automorphisms of $X$, which acts without inversions and with compact open vertex stabilisers.  Then a subgroup $\G < G$ is discrete if and only if it acts on $X$ with finite vertex stabilisers, and the Haar measure $\mu$ on $G$ may be normalised so that the covolume of a discrete $\Gamma < G$ is given by $  \mu(\G \bs G) = \sum |\G_y|^{-1} $, where the sum is over a set of representatives $y$ of the vertices of the quotient graph $Y = \G \bs X$ (Bass--Lubotzky \cite{BL}).  Hence a discrete subgroup $\G < G$ is a lattice if and only if the series $\sum |\G_y|^{-1} $ converges. Moreover, a discrete subgroup $\G < G$ is a cocompact lattice in $G$ if and only if $Y$ is finite.

\subsection{Kac--Moody groups}\label{s:kac-moody}

We first in Section~\ref{s:basic_definitions} explain how one may associate, to a generalised Cartan matrix $A$ and an arbitrary field, a Kac--Moody group $\Lambda$, the so-called minimal or incomplete Kac--Moody group.  In Section~\ref{s:rank 2} we specialise to rank $2$ Kac--Moody groups over finite fields.  Section~\ref{s:completions} describes the completion $G$ of $\Lambda$ that appears in the statement of Theorem~\ref{t:cocompact existence} above and Section~\ref{s:action ends} recalls some results from \cite{CG} concerning the action of $G$ on the set of ends of its tree $X$.   Our treatment of Kac--Moody groups is brief and combinatorial, and partly follows Dymara--Januszkiewicz~\cite[Appendix TKM]{DJ}.  For a more sophisticated and general approach, using the notion of a ``twin root datum", we refer the reader to, for example, Caprace--R\'emy~\cite{CR}.

\subsubsection{Incomplete Kac--Moody groups}\label{s:basic_definitions}

Let $I$ be a finite set.  A \emph{generalised Cartan matrix} $A=(A_{ij})_{i,j \in I}$ is a matrix with integer entries, such that $A_{ii} = 2$, $A_{ij} \leq 0$ if $i \neq j$ and $A_{ij} = 0$ if and only if $A_{ji} = 0$.  (If $A$ is positive definite, then $A$ is the Cartan matrix of some
finite-dimensional semisimple Lie algebra.)  A \emph{Kac--Moody datum} is a $5$--tuple $(I,\mathfrak{h}, \{\alpha_i\}_{i \in I}, \{h_i\}_{i \in I},A)$ where $\mathfrak{h}$ is a finitely generated free abelian group, $\alpha_i \in \mathfrak{h}$, $h_i \in \Hom(\mathfrak{h},\mathbb{Z})$ and $A_{ij} = h_j(\alpha_i)$.  The set $\Pi = \{ \alpha_i \}_{i \in I}$ is called the set of \emph{simple roots}.

Given a generalised Cartan matrix $A$ as above, we define a \emph{Coxeter matrix} $M=(m_{ij})_{i,j\in I}$ as follows: $m_{ii} = 1$, and if $i \neq j$ then $m_{ij} = 2,3,4,6$ or $\infty$ as $A_{ij}A_{ji} = 0,1,2,3$ or is $\geq 4$. The associated \emph{Weyl group} $W$ is then the Coxeter group with presentation determined by $M$: \[ W = \langle \{w_i\}_{i \in I} \mid (w_iw_j)^{m_{ij}} \mbox{for $m_{ij} \neq \infty$} \rangle. \]  The Weyl group acts on $\mathfrak{h}$ via $w_i: \beta \mapsto \beta - h_i(\beta)\alpha_i$ for each $\beta \in \mathfrak{h}$ and each $i \in I$.  In particular, $w_i(\alpha_i) = -\alpha_i$ for each simple root $\alpha_i$.  The set  $\Phi$ of \emph{real roots} is defined by $\Phi = W \cdot \Pi$. In general, the set of real roots is infinite.

We will, not by coincidence, use the same terminology and notation for simple roots and real roots which are defined in the following combinatorial fashion.  Let $\ell$ be the word length on the Weyl group $W$, that is, $\ell(w)$ is the minimal length of a word in the letters $\{w_i\}_{i\in I}$ representing $w$.  The \emph{simple roots} $\Pi = \{ \alpha_i \}_{i \in I}$ are then defined by \[ \alpha_i = \{ w \in W \mid \ell(w_i w) > \ell(w) \}.\] The set $\Phi$ of \emph{real roots} is $\Phi = W\cdot \Pi = \{ w\alpha_i \mid w \in W, \alpha_i \in \Pi \}$, and $W$ acts naturally on $\Phi$.  The set $\Phi_+$ of \emph{positive roots} is $\Phi_+ = \{ \alpha \in \Phi \mid 1_W \in \alpha \}$, and the set of \emph{negative roots} $\Phi_-$ is $\Phi \bs \Phi_+$.  The complement of a root $\alpha$ in $W$, denoted $-\alpha$, is also a root.  As before, $w_i(\alpha_i) = -\alpha_i$ for each simple root $\alpha_i$.

We now define the split Kac--Moody group $\Lambda$ associated to a Kac--Moody datum as above, over an arbitrary field $k$.  The group $\Lambda$ may be given by a presentation, which is essentially due to Tits (see~\cite{Tits87}), and which appears in Carter~\cite{Carter}.  For simplicity, we state this presentation only for the simply-connected group $\Lambda_u$ and then discuss the general case.  Let $(I,\mathfrak{h}, \{\alpha_i\}_{i \in I}, \{h_i\}_{i \in I},A)$ be a Kac--Moody datum and $k$ a field.  The associated \emph{simply-connected Kac--Moody group $\Lambda_u$} over $k$ is generated by \emph{root subgroups} $U_\alpha = U_\alpha(k) = \langle x_\alpha(t) \mid t \in k \rangle$, one for each real root $\alpha \in \Phi$.  We write $x_i(u)=x_{\alpha_i}(u)$ and $x_{-i}(u)=x_{-\alpha_i}(u)$ for each $u \in k$ and $i \in I$, and put $\tilde{w}_i(u)=x_i(u)x_{-i}(u^{-1})x_i(u)$, $\tilde{w}_i=\tilde{w}_i(1)$ and $h_i(u)=\tilde{w}_i(u)\tilde{w}_i^{-1}$ for each $u \in k^*$ and $i \in I$.  A set of defining relations for the simply-connected Kac-Moody group $\Lambda_u$ is then: \begin{enumerate} \item\label{i:root_subgroups} $x_{\alpha}(t)x_{\alpha}(u)=x_{\alpha}(t+u)$, for all roots $\alpha \in \Phi$ and all $t,u \in k$.  \item\label{i:steinberg} If $\alpha, \beta\in\Phi$ is a prenilpotent pair of roots, that is, there exist $w,w'\in W$ such that $w(\alpha)\in\Phi_+$, $w(\beta)\in\Phi_+$, $w'(\alpha)\in\Phi_-$ and $w'(\beta)\in\Phi_-$, then for all $t,u \in k$: $$[x_{\alpha}(t),x_{\beta}(u)]=\prod_{\begin{array}{c}  i,j\in\mathbb{N}\\ i\alpha+j\beta\in\Phi\end{array}}  x_{i\alpha+j\beta}(C_{ij\alpha\beta}t^iu^j)$$ where the integers $C_{ij\alpha\beta}$ are uniquely determined by $i, j,  \alpha, \beta$, $\Phi$, and the ordering of the terms on the right-hand side. \item $h_i(t)h_i(u)=h_i(tu)$ for all $t,u\in k^*$ and all $i \in I$. \item $[h_i(t),h_j(u)]=1$ for all $t,u\in k^*$ and $i, j \in I$.  \item $h_j(u)x_i(t)h_j(u)^{-1}=x_i(u^{A_{ij}}t)$ for all $t \in k$, $u \in k^*$ and $i,j \in I$.\item $\tilde{w}_ih_j(u)\tilde{w}_i^{-1}=h_j(u)h_i(u^{-A_{ij}})$  for all $u \in k^*$ and $i,j \in I$.\item $\tilde{w}_ix_{\alpha}(u)\tilde{w}_i^{-1}=x_{w_i(\alpha)}(\epsilon u)$ where $\epsilon\in\{\pm 1\}$, for all $u \in k$. \end{enumerate}

By a result of P.-E.~Caprace (cf.  \cite[3.5(2)]{PEC}), any two split Kac--Moody groups of the same type defined over the same field are isogenic.  That is, if $\Lambda$ is any split Kac--Moody group associated to the same generalised Cartan matrix $A$ as $\Lambda_u$, and defined over the same field $k$, then there exists a surjective  homomorphism $i:\Lambda_u\rightarrow\Lambda$ with $\ker(i)\leq Z(\Lambda_u)$. The Kac--Moody group $\Lambda$ so constructed is sometimes called the  \emph{incomplete} or \emph{minimal} Kac--Moody group (for completions of $\Lambda$, see Section~\ref{s:completions} below).

A first example of an incomplete Kac--Moody group $\Lambda$ over a finite field is $\Lambda = SL_n(\F_q[t,t^{-1}])$, which is over the field $\F_q$, and is not simply-connected.

We now discuss several important subgroups of the Kac--Moody group $\Lambda$.  For any version (simply-connected or not), the \emph{unipotent} subgroup of $\Lambda$ is \[U = U_+ = \langle U_\alpha \mid \alpha \in \Phi_+ \rangle. \]   For $\Lambda_u$ simply-connected, the \emph{torus} \[T = \langle h_i(u) \mid i \in I, u \in k^* \rangle\] is isomorphic to the direct product of $|I|$ copies of $k^*$.  In general, the torus $T$ of $\Lambda$ is a homomorphic image of the direct product of $|I|$ copies of $k^*$.  For all $\Lambda$, we define $N$ to be the subgroup of $\Lambda$ generated by the torus $T$ and by the elements $ \{\tilde{w}_i\}_{i \in I} $ (where, in general as in the simply-connected case, $\tilde{w}_i = x_{\alpha_i}(1)x_{-\alpha_i}(1)x_{\alpha_i}(1)$ for all $i \in I$).  The \emph{standard Borel subgroup} $B = B_+$ of $\Lambda$ is defined by \[B = \langle T, U_+ \rangle = \langle T, U\rangle.  \] The group $B$ has decomposition $B = T \ltimes U_+ = T \ltimes U$ (see~\cite{Remy02}).

The subgroups $B$ and $N$ of $\Lambda$ form a $BN$--pair (also known as a Tits system) with Weyl group $W$, and hence $\Lambda$ has a Bruhat--Tits building $X$.   (In fact, the group $\Lambda$ has isomorphic twin buildings $X_+ \cong X_-$, associated to twin $BN$--pairs $(B_+,N)$ and $(B_-,N)$ respectively.) The chambers of $X$ correspond to the cosets of $B$ in $\Lambda$, hence $\Lambda$ acts naturally on $X$ with quotient a single chamber.  For each apartment $\Sigma$ of $X$, the chambers in $\Sigma$ are in bijection with the elements of the Weyl group $W$.  Each root $\alpha \subset W$ corresponds to a ``half-apartment".  The construction of the building $X$ for $\Lambda$ of rank $2$ is explained further in Section~\ref{s:rank 2} below.

\subsubsection{Rank $2$}\label{s:rank 2}

We now specialise to the cases considered in Theorem~\ref{t:cocompact existence} above.  Let $A$ be a generalised Cartan matrix of the form $A = \left(\begin{matrix}2 & -m\\-m & 2 \end{matrix}\right)$, with $m\geq 2$.  If $m = 2$, then $A$ is \emph{affine}, meaning that $A$ is positive semidefinite but not positive definite.  For all such $A$ (affine and non-affine) the associated Weyl group $W$ is \[ W = \langle w_1, w_2 \mid w_1^2, w_2^2 \rangle. \] That is, $W$ is the infinite dihedral group.  Let $\ell$ be the word length on $W$. The simple roots $\Pi = \{ \alpha_1,\alpha_2\}$ are then given by, for $i = 1,2$, \[\alpha_i = \{ w \in W \mid \ell(w_iw) > \ell(w)\} = \{1, w_{3-i}, w_{3-i}w_i, w_{3-i}w_iw_{3-i}, \ldots\}.\] The set $\Phi$ of real roots is $\Phi = \{ w\alpha_i \mid w \in W, i = 1,2 \}$.   

Now let $\Lambda$ be an incomplete Kac--Moody group with generalised Cartan matrix $A$, defined over a finite field $\F_q$, where $q = p^a$ with $p$ prime.  As $\Lambda$ is a group with $BN$--pair, as described above, for $i = 1,2$, the \emph{parabolic subgroup} $P_i$ of $\Lambda$ is defined by \[P_i = B \sqcup B\tilde{w}_iB.\]

Since $J_i=\{\alpha_i \}$ is a root system of type $A_1$, and thus is of finite type, now \cite[6.2]{Remy02} applies. Hence, the group $P_i$ has a \emph{Levi decomposition} $P_i = L_i \ltimes U_i$.  Here $U_i=U\cap U^{w_i}$ is called a \emph{unipotent radical} of $P_i$, and the group $L_i$ is called a \emph{Levi complement} of $P_i$.  The Levi complement factors as $L_i = TM_i$, where $T$ is the torus of $\Lambda$, and $M_i=\langle U_{\alpha_i}, U_{-\alpha_i}\rangle$, that is, $A_1(q)\cong M_i\triangleleft L_i$.

To describe the building $X$ for $\Lambda$, we first describe its apartments.  Let $\Sigma$ be the \emph{Coxeter complex} for the Weyl group $W$ (the infinite dihedral group).  Then $\Sigma$ is the one-dimensional simplicial complex homeomorphic to the line.  The set of real roots $\Phi$, described above, is in bijection with the set of half-lines in $\Sigma$.  Thus we may regard each real root as a geodesic ray in $X$.  So each real root determines an end of $X$.  

In particular, the set $\Phi^+$ of positive real roots is the disjoint union of the sets \[ \Phi^1_+ := \{ \alpha_1, w_1\alpha_2, w_1w_2\alpha_1, w_1w_2w_1\alpha_2, \ldots, (w_1w_2)^n\alpha_1, (w_1w_2)^n w_1 \alpha_2,\ldots \} \] and  \[ \Phi^2_+ := \{ \alpha_2, w_2\alpha_1, w_2w_1\alpha_2, w_2w_1w_2\alpha_1, \ldots, (w_2w_1)^n\alpha_2, (w_2w_1)^n w_2 \alpha_1,\ldots \} \] and the roots in $\Phi^+$ are in bijection with the set of half-lines in $\Sigma$ which contain the base chamber $B$.  See Figure \ref{f:roots} below.

\begin{figure}[ht]
\scalebox{0.9}{\includegraphics{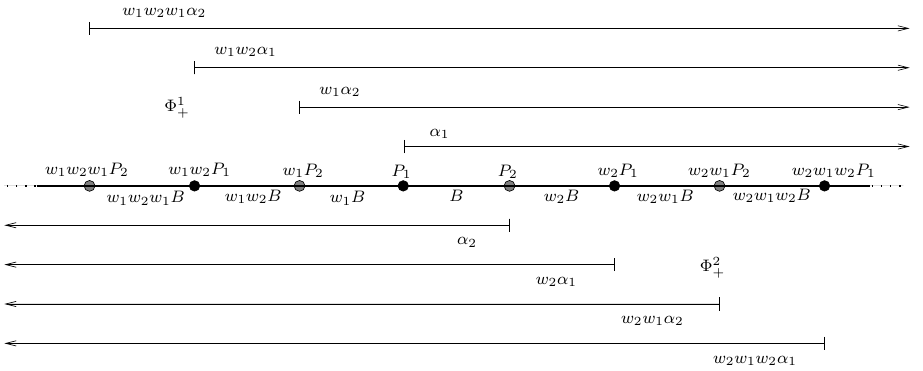}}
\caption{The sets of positive real roots $\Phi^1_+$ and $\Phi^2_+$, with each such root identified with a half-apartment of the standard apartment
$\Sigma$ containing $B$.}
\label{f:roots}
\end{figure}
 
The apartments of the building $X$ are copies of the Coxeter complex $\Sigma$ for $W$ and $X$ is a $(q+1)$--regular simplicial tree, with chambers the edges of this tree.  

\subsubsection{Completions of $\Lambda$}\label{s:completions}

We are finally ready to describe the main object of our study: the locally compact topological Kac--Moody groups.  In order to do this we will have to define a topological completion of the incomplete Kac--Moody group $\Lambda$.  There are several completions appearing in the literature.  For example, Carbone--Garland~\cite{CG} defined a representation-theoretic completion of $\Lambda$ using the ``weight topology".  A different approach by R\'emy and Ronan, appearing for instance in~\cite{RR}, is to use the action of $\Lambda$ on the building $X$, as follows.  The kernel $K$ of the $\Lambda$--action on $X$ is the centre $Z(\Lambda)$, which is a finite group when $\Lambda$ is over a finite field (R\'emy~\cite{Remy02}).  The closure of $\Lambda/K$ in the automorphism group of $X$ is then a completion of $\Lambda$.  For example, when $\Lambda = SL_n(\F_q[t,t^{-1}])$, the centre $Z(\Lambda)$ is the finite group $\mu_n(\F_q)$ of $n$th roots of unity in $\F_q$, and the completion in this topology is $SL_n(\F_q((t^{-1})))/\mu_n(\F_q) \cong PSL_n(\F_q((t^{-1})))$.  To avoid dealing with representation-theoretic constructions or with quotients, we are going to follow instead the completion in the building topology, defined by Caprace and R\'emy in~\cite{CR}.

So, let $\Lambda$ be an incomplete Kac--Moody group over a finite field, as defined in Section~\ref{s:basic_definitions} above.  We now describe the completion $G$ of $\Lambda$ which appears in Theorem~\ref{t:cocompact existence} (for $\Lambda$ with generalised Cartan matrix $A$ as in Section~\ref{s:rank 2} above).  Let $c_+=B_+$ be the chamber of the Bruhat--Tits building $X$ for $\Lambda$ which is fixed by $B=B_+$.  For each $n \in \mathbb{N}$, let \[U_{+,n} = \{ g \in U_+ \mid g \cdot c = c \mbox{ for each chamber $c$ such that $d(c,c_+) \leq n$ }\}.\] That is, $U_{+,n}$ is the kernel of the action of $U_+ = U$ on $\Ball(c_+,n)$.  Now define a function $\dist_+:\Lambda \times \Lambda \to \mathbb{R}_+$ by $\dist_+(g,h) = 2$ if $h^{-1}g \not \in U_+$, and $\dist_+(g,h) = 2^{-n}$ if $g^{-1}h \in U_+$ and $n = \max\{ k \in \mathbb{N} \mid g^{-1}h \in U_{+,k} \}$.  Then $\dist_+$ is a left-invariant metric on $\Lambda$ (see~\cite{CR}).  Let $G$ be the completion of $\Lambda$ with respect to this metric.  The group $G$ is called the \emph{completion of $\Lambda$ in the building topology}, and we will refer to $G$ as a \emph{topological Kac--Moody group}.  For example, when $\Lambda = SL_n(\F_q[t,t^{-1}])$, the topological Kac--Moody group $G$ is $G = SL_n(\F_q((t^{-1})))$.

Some properties of topological Kac--Moody groups that we will need are gathered in Proposition~\ref{p:completion} below.  We state these results only for $G$ as in Theorem~\ref{t:cocompact existence} above, although they hold more generally.

\begin{prop}\label{p:completion} Let $G$ be a topological Kac--Moody group as in Theorem~\ref{t:cocompact existence} above, with $G$ being the completion in the building topology of an incomplete Kac--Moody group $\Lambda$.  \begin{enumerate} \item\label{i:complete_G} $G$ is a locally compact, totally disconnected topological group. \item\label{i:Levi} Let $\hat{B}$, $\hat{U}$, $\hat{P_i}$ and $\hat{U_i}$ be the closures in $G$ of the subgroups $B=B_+$, $U=U_+$, $P_i$ and $U_i$ respectively of $\Lambda$.  Then $\hat{B} \cong T \ltimes \hat{U}$ and $\hat{P_i} \cong L_i \ltimes \hat{U_i}$. \item\label{i:building} $(\hat{B},N)$ is a $BN$--pair of $G$.  The corresponding building is canonically isomorphic to $X$, and so by abuse of notation we will denote it by $X$ as well.  The kernel of the action of $G$ on $X$ is the centre $Z(G)$, and $Z(G) = Z(\Lambda)$. \end{enumerate} \end{prop}

\noindent Items~\eqref{i:complete_G} and~\eqref{i:building} are established by Caprace--R\'emy in \cite{CR}, and  item~\eqref{i:Levi}  in \cite{CR} and \cite{PEC}.

Following the terminology used by R\'emy and Ronan \cite[p. 196]{RR} for topological Kac--Moody groups, we will refer to $\hat{B}$ as the \emph{(standard) Iwahori  subgroup} of $G$, and to $\hat{P_1}$ and $\hat{P_2}$ as the \emph{(maximal or standard) parahoric subgroups} of $G$.   To simplify notation, when the context is clear we will write $B$, $P_1$ and $P_2$ for the Iwahori and maximal parahoric subgroups of the topological Kac--Moody group $G$, rather than respectively $\hat{B}$, $\hat{P_1}$ and $\hat{P_2}$.

\subsubsection{Action of $G$ on ends}\label{s:action ends}

In this section we recall definitions and results from Carbone--Garland~\cite{CG} concerning the action of $G$ on the set of ends of $X$.  These will be used in Section \ref{s:motivation} below where we discuss our conjectures about $p$--elements in $G$.

For $i = 1,2$, let $\ep_i$ be the unique end of $X$ determined by the simple root $\alpha_i$ (see Figure \ref{f:roots} above).  We partition $\bX$ into two sets $\mathcal{E}_1$ and $\mathcal{E}_2$, so that $\ep_i \in \mathcal{E}_i$ for $i = 1,2$, as follows.  For $i = 1,2$, $\mathcal{E}_i$ consists of those ends represented by a geodesic ray which starts at the vertex $P_{3-i}$ and does not contain the edge $B$.

\begin{lemma}\label{l:U+_transitive_ends}  For any end $\varepsilon \in \mathcal{E}_1$, there is a $u_1 \in \hat{U}_+$ such that $u_1 \cdot \varepsilon = \varepsilon_1$.  For any end $\varepsilon' \in \mathcal{E}_2$, there is a $u_2 \in \hat{U}_+$ such that $u_2 \cdot \varepsilon' = \varepsilon_2$.  Moreover, $u_2$ may be chosen to fix $\varepsilon_1$. \end{lemma}

\begin{proof} This follows from~\cite[Lemma 14.1]{CG}.  In \cite{CG}, a different completion of minimal Kac--Moody groups was used, but the same proof goes through in the completion that we are using. \end{proof}

\noindent The following corollary of Lemma \ref{l:U+_transitive_ends} is the same statement as \cite[Corollary 14.1]{CG}.

\begin{corollary}\label{l:G_2-transitive_ends}  The group $G$ acts $2$--transitively on the set of ends of $X$. \end{corollary}

Now let $V_1$ be the elementary abelian $p$--subgroup of $U_+$ defined by \[ V_1 := \langle U_\alpha \mid \alpha \in \Phi^1_+ \rangle\]
and let $-V_2$ be the elementary abelian $p$--subgroup of $U_-$ defined by \[-V_2 := \langle U_\alpha \mid -\alpha \in \Phi^2_+ \rangle.\] Then the groups $V_1$ and $-V_2$ both fix the end $\varepsilon_1$. Notice that $[V_1,-V_2]=1$ and $V_1\cap -V_2=1$, and so $\langle V_1, -V_2\rangle=V_1\times -V_2$.  Let $\cU$ be the closed abelian  group \[\cU := \widehat{V_1\times  -V_2} \] and let \[ \mathcal{M} := \bigcap_{w \in W} w \hat{B} w^{-1}.\]  Note that $\mathcal{M}$ is precisely the subgroup of $G$ which fixes the standard apartment $\Sigma$ pointwise.  In particular, the torus $T$ is a (finite) subgroup of $\mathcal{M}$, of exponent $(q - 1)$.

Take $g_\tau\in G$ such that $g_\tau$ induces the element $\tau:=w_1w_2\in W$.  The element $g_\tau$ is hyperbolic, with axis the standard apartment $\Sigma$, translation length two edges, and with attracting fixed point $\varepsilon_2$ and repelling fixed point $\varepsilon_1$.   (In Figure \ref{f:roots} above, $g_\tau$ moves the standard apartment two edges to the left.)  Let $R$ be the infinite cyclic subgroup of $G$ generated by $g_\tau$.

The following decomposition of the stabiliser $G_{\ep_1}$ is from \cite{CG}, where the group $\mathcal{M}$ is denoted by  $\mathcal{B}_\mathcal{I}$ and $R$ is denoted by $T$.  As with Lemma \ref{l:U+_transitive_ends} above, the same proof applies even though we are using a different completion to \cite{CG}.  

\begin{theorem}[Theorem 14.1, \cite{CG}]\label{t:end stabiliser} The group $G_{\ep_1}$ can be expressed as $G_{\ep_1} = \cU\mathcal{M} R $. \end{theorem}

\subsection{Finite groups}\label{s:finite groups}

In our constructions of cocompact lattices in Kac--Moody groups and the proofs of Theorems \ref{t:classification} and \ref{t:covolumes}, we will need to look carefully at the finite subgroups of $G$. The following celebrated result of L.E.~Dickson and its corollary will be especially useful for us.

\begin{theorem}[Dickson, 6.5.1 of \cite{GLS3}] \label{t:Dickson} Let $K=PSL_2(q)$, where $q=p^a\geq 5$ and $p$ is a prime. Set $d=(2,q-1)$. Then $K$ has subgroups of the following isomorphism types in the indicated cases, and every subgroup of $K$ is isomorphic to a subgroup of one of the following groups: \begin{enumerate} \item Borel subgroups of $K$, which are Frobenius groups of order $q(q-1)/d$; \item Dihedral groups of orders $2(q-1)/d$ and $2(q+1)/d$; \item The groups $PGL_2(p^b)$ (if $2b\mid a$) and $PSL_2(p^b)$ (if $b$ is a proper divisor of $a$); \item The alternating group $A_5$, if $5$ divides $|K|$; \item The symmetric group $S_4$, if $8$ divides $|K|$; and \item The alternating group $A_4$. \end{enumerate} \end{theorem}

\begin{corollary}\label{c:Dickson} Let $K=SL_2(q)$, where $q=p^a$ with $p$  prime, and suppose $A$ is a proper  subgroup of $K$.

 If $p=2$ and $q+1$ divides $|A|$, then either $A\cong C_{q+1}$, a cyclic group of order $q+1$, or $A\cong D_{2(q+1)}$, a dihedral group of order $2(q+1)$.

If $p$ is odd and the image of $A$ in $K/Z(K)\cong PSL_2(q)$ has order divisible by $q+1$, then  $Z(K)=\langle -I\rangle\leq A$. Moreover, either $A$ is a subgroup of $K$ of order $2(q+1)$ such that $A/Z(K) \cong D_{q+1}$, a dihedral group of order $q+1$, or one of the following conditions hold:

\begin{enumerate}
\item $q=5$,  $A\cong SL_2(3)$;
\item $q=7$, $A\cong 2S_4$;
\item $q=9$, $A\cong SL_2(5) $;
\item $q=11$, $A\cong SL_2(3)$ or $A\cong SL_2(5)$;
\item $q=19$, $A\cong SL_2(5) $;
\item $q=23$, $A\cong 2S_4$;
\item $q=29$, $A\cong SL_2(5)$; or
%\item $q=47$, $A\cong 2S_4$
\item $q=59$, $A\cong SL_2(5)$.
\end{enumerate}

\end{corollary}
 
\begin{proof} Suppose that $p=2$. Then $d=1$ and $SL_2(q)=PSL_2(q)$. Assume  first that $q\geq 5$. Then if $q+1$ divides $|A|$, Theorem \ref{t:Dickson} above asserts that both $C_{q+1}$ and $D_{2(q+1)}$ are the obvious candidates for the role of $A$.  If not, $A$ would be one of the following groups: $A_4$, $S_4$ or $A_5$. Then $q+1$ would divide $12$, $24$ or $60$. Since $q$ is a power of $2$ and $q\geq 5$, this is not possible, proving the result.  Otherwise $q\in\{ 2,  4\}$, and the result follows immediately from the structures of $K=SL_2(2)\cong S_3$ and $K=SL_2(4)\cong A_5$.

Suppose now that $p$ is odd. This time $d=2$ and the image of $A$ in $PSL_2(q)$ is a group of order divisible by $q+1$.  Since $|A|$ is even while $K$ contains a unique involution $-I$, $\langle -I\rangle=Z(K)\leq A$.  If $q\geq 5$,  using the same argument as above, we obtain the desired conclusion. Otherwise $q=3$ and $K=SL_2(3)\cong Q_8C_3$, and the result follows immediately. \end{proof}
 
 We will also need another result  of L.E. Dickson about the subgroups of $PGL_2(q)$, stated in the following form in \cite{COTR}.
 
 \begin{theorem}[Theorem 2 of \cite{COTR}] \label{t:Dickson2}
 Let $q\equiv\epsilon\mod 4$ where $\epsilon=\pm1$. The subgroups of $PGL_2(q)$ are as follows.
 \begin{enumerate}
 \item Two conjugacy classes of cyclic subgroups $C_2$.
 \item One conjugacy class of cyclic subgroups of $C_d$ where $d\mid(q\pm\epsilon)$ and $d>2$.
 \item Two conjugacy classes of dihedral subgroups $D_4$ (only one class contained in the subgroup $PSL_2(q)$).
 \item\label{i:Dickson dihedral} Two conjugacy classes of dihedral subgroups $D_{2d}$, where $d\mid\frac{q\pm\epsilon}{2}$ and $d>2$
 (only one class contained in the subgroup $PSL_2(q)$).
 \item\label{i:Dickson dihedral 2} One conjugacy class of dihedral subgroups $D_{2d}$, where $(q\pm\epsilon)/d$ is an odd integer and $d>2$.
 \item Subgroups $A_4$, $S_4$ and $A_5$ when $q\equiv\pm1\mod 10$. There is only one conjugacy class of any of these types of subgroups and all lie in the subgroup $PSL_2(q)$, except for $S_4$ when $q\equiv\pm3\mod 8$.
 \item Subgroups $PSL(2,p^m)$, $PGL(2,p^m)$, the elementary Abelian group of order $p^m$ and a semidirect product of the elementary Abelian group of order $p^m$ and the cyclic group of order $d$, where $m\leq\log_pq$, $d\mid (q-1)$ and $d\mid (p^m-1)$.
 \end{enumerate}
 \end{theorem}

\section{Construction of cocompact lattices}\label{s:cocompact}

In this section we prove Theorems~\ref{t:cocompact existence} and~\ref{t:exceptions}, stated in the introduction.  We first establish our embedding criterion for cocompact lattices, Proposition~\ref{p:embedding}, in Section~\ref{s:embedding} below.  We then apply this criterion in Section~\ref{s:construction edge-transitive} to construct the edge-transitive lattices $\G$ in Theorems~\ref{t:cocompact existence} and \ref{t:exceptions}, and in Section~\ref{s:construction} to construct the lattices $\Gamma'$ in Theorem~\ref{t:cocompact existence}.

\subsection{Embedding criterion}\label{s:embedding}

Our embedding criterion applies not only to $G$ as in Theorem \ref{t:cocompact existence}, but to more general locally compact groups acting on trees, as follows.

Let $q_1$ and $q_2$ be positive integers and let $X$ be the $(q_1 + 1,q_2 + 1)$--biregular tree.  Let $G$ be any closed locally compact group of automorphisms of $X$ which acts on $X$ without inversions, with compact open vertex stabilisers $G_x$ for $x \in VX$ and with fundamental domain an edge $[x_1,x_2]$, where for $i = 1,2$ the vertex $x_i$ has valence $q_i + 1$. Denote by $P_i$ the stabiliser $G_{x_i}$ for $i = 1,2$, and let $B = P_1 \cap P_2$.  For notational convenience, we denote by $C$ the subgraph of $X$ with vertex set $\{x_1, x_2\}$ and edge set $\{ f,\overline{f} \}$, such that $\partial_0(f) = \partial_1(\overline{f}) = x_1$ and  $\partial_1(f) = \partial_0(\overline{f}) = x_2$.  Then $G$ is the fundamental group of an edge of groups $\mathbb{G}$ over $C$, as sketched in the introduction.

For some integer $n \geq 1$ dividing both $q_1 + 1$ and $q_2 + 1$, let $A = A_n$ be the graph with two vertices $a_1$ and $a_2$ and edge set $\{e_1,\ldots,e_n, \overline{e_1},\ldots,\overline{e_n}\}$, so that $\partial_0(e_j) = \partial_1(\overline{e_j}) = a_1$ and $\partial_1(e_j) = \partial_0(\overline{e_j}) = a_2$ for $1 \leq j \leq n$.  The case $n = 2$ is sketched in the statement of Theorem \ref{t:cocompact existence}.  We now state and prove a sufficient criterion for the fundamental group of a graph of groups over $A$ to embed in $G$ as a cocompact lattice.  

\begin{prop}\label{p:embedding} 

Suppose that there are finite groups $A_1 \leq P_1$ and $A_2 \leq P_2$ such that:
\begin{enumerate}
\item\label{i:n orbits} for $i = 1,2$, the group $A_i$ has $n$ orbits of equal size on $E_X(x_i)$;
\item\label{i:representatives} there are: \begin{itemize}\item representatives $f_1 = f$, $f_2$, \ldots, $f_{n}$ of the orbits of $A_1$ on $E_X(x_1)$ and $\hat{f}_{1} = \overline{f}$, $\hat{f}_2$, \ldots, $\hat{f}_{n}$ of the orbits of $A_2$ on $E_X(x_2)$; and \item elements $g_{1} = 1$, $g_2$,  \ldots, $g_{n} \in P_1$ and $\hat{g}_{1} = 1$, $\hat{g}_2$, \ldots, $\hat{g}_{n} \in P_2$;\end{itemize} such that for $1 \leq j \leq n$:
\begin{enumerate}
\item\label{i:g on f} $g_{j} \cdot f_{1} = f_{j}$ and $ \hat{g}_{j} \cdot \hat{f}_{1} = \hat{f}_{j}$;
\item\label{i:edge groups} $A_1 \cap B^{g_{j}} = A_1 \cap A_2 = A_2 \cap B^{\hat{g}_{j}}$; and
\item\label{i:g on edges} $(A_1 \cap A_2)^{g_{j}} = (A_1 \cap A_2)^{\hat{g}_{j}}$.
\end{enumerate}
\end{enumerate}

 Let $\mathbb{A}$ be the graph of groups over $A$ with: \begin{itemize} \item vertex groups $\mathcal{A}_{a_i} = A_i$ for $i = 1,2$; \item  edge groups $\mathcal{A}_{e_j} = \mathcal{A}_{\overline{e_j}} = A_1 \cap A_2$ for $1 \leq j \leq n$; \item each monomorphism $\alpha_{e_j}$ from an edge group $A_1 \cap A_2$ into $A_1$ inclusion composed with $\ad(g_{j}\hat{g}_{j}^{-1})$; and \item the monomorphisms $\alpha_{\overline{e_j}}$ from edge groups $A_1 \cap A_2$ into $A_2$ inclusions.  \end{itemize}

Then the fundamental group of the graph of groups $\mathbb{A}$ is a cocompact lattice in $G$, with quotient $A$.

\end{prop}

\noindent  Note that in the special case $n = 1$, where for $i = 1,2$ the group $A_i$ acts transitively on $E_X(x_i)$, condition \eqref{i:representatives} reduces to the requirement that $\Stab_{A_i}(x_{3 - i}) = A_1 \cap A_2$ for $i = 1,2$. 

\begin{proof} We construct a covering of graphs of groups $\Phi:\mathbb{A} \to \mathbb{G}$.  Since $A$ is a finite graph and the vertex groups $A_1$ and $A_2$ are finite, it then follows from our discussion of lattices in Section \ref{s:lattices} and Proposition \ref{p:coverings} above that the fundamental group of $\mathbb{A}$ is a cocompact lattice in $G$ with quotient the graph $A$.

Let $\theta:A \to C$ be the graph morphism given by $\theta(a_i) = x_i$ for $i = 1,2$, and $\theta(e_j) = f$ and $\theta(\overline{e_j}) = \overline{f}$ for $1 \leq j \leq n$.  We construct a morphism of graphs of groups $\Phi:\mathbb{A} \to \mathbb{G}$ over $\theta$ as follows.  For $i = 1,2$ let $\phi_{a_i}:\mathcal{A}_{a_i} \to P_i$ be the natural inclusion $A_i \hookrightarrow P_i$.  For $1 \leq j \leq n$ let $\phi_{e_j}:\mathcal{A}_{e_j} \to B$ be the composition of the natural inclusion $A_1 \cap A_2 \hookrightarrow B^{\hat{g}_{j}}$ with the map $\ad(\hat{g}_{j}^{-1}):B^{\hat{g}_{j}} \to B$.  Define $\phi(e_j) = g_{j}$ and $\phi(\overline{e_j}) = \hat{g}_{j}$.  Then it may be checked that $\Phi$ is indeed a morphism of graphs of groups.

To show that $\Phi$ is a covering, we first show that the map \[ \Phi_{a_1/f}: \coprod_{j = 1}^n \cA_{a_1}/\alpha_{e_j}(\cA_{e_j}) \to P_1/B \] induced by $g \mapsto \phi_{a_1}(g)\phi(e_j) = g g_{j}$ for  $g$ representing a coset of $\alpha_{e_j}(\cA_{e_j}) = (A_1 \cap A_2)^{g_{j}\hat{g}_{j}^{-1}} = A_1 \cap A_2$ in $\cA_{a_1} = A_1$ is a bijection.   For this, we note that since the edges $f_j = g_j \cdot f_1 = g_j \cdot f$ represent pairwise distinct $A_1$--orbits on $E_X(x_1)$, for all $g, h \in A_1$ and all $1 \leq j \neq j' \leq n$ the cosets $gg_{j}B$ and $hg_{j'}B$ are pairwise distinct.  The conclusion that $\Phi_{a_1/f}$ is a bijection then follows from the assumption that $A_1 \cap B^{g_{j}} = A_1 \cap A_2$.

The proof that the map \[  \Phi_{a_2/\overline{f}}: \coprod_{j = 1}^n \cA_{a_2}/\alpha_{\overline{e_j}}(\cA_{\overline{e_j}}) \to P_2/B \] is a bijection is similar.  We conclude that $\Phi:\mathbb{A} \to \mathbb{G}$ is a covering of graphs of groups, as desired. \end{proof}

\subsection{Construction of edge-transitive lattices}\label{s:construction edge-transitive}

Let $G$ be as in Theorem \ref{t:cocompact existence}.  We now apply Proposition \ref{p:embedding} in the case $n = 1$ to construct the edge-transitive amalgams $\G = A_1 *_{A_0} A_2$ described in Theorems \ref{t:cocompact existence} and \ref{t:exceptions}.  Throughout this section, $P_1$ and $P_2$ are the maximal parahoric subgroups of $G$, and for $i = 1,2$ the group $P_i$ is the stabiliser in $G$ of the vertex $x_i$ of $X$.  Recall that each $P_i$ has Levi decomposition $P_i = L_i \ltimes U_i$, and $L_i = T M_i$ where $T$ is a fixed maximal split torus of $G$ and $A_1(q) \cong M_i \lhd L_i$.  

We now prove parts \eqref{i: p = 2}, \eqref{i:PSL} and \eqref{i:PGL} of Theorem \ref{t:cocompact existence} in Sections \ref{proof: p = 2}, \ref{proof:PSL} and \ref{proof:PGL} respectively, and then prove Theorem \ref{t:exceptions} in Section \ref{proof:exceptions}.  In all cases below, $H_i$ is a non-split torus of $M_i$ such that $|N_T(H_i)|$ is as large as possible.

\subsubsection{Part \eqref{i: p = 2} of Theorem \ref{t:cocompact existence}}\label{proof: p = 2}

We begin with the case $p = 2$.  Since $q$ is even, $SL_2(q) \cong PSL_2(q)$ and so $M_i \cong SL_2(q)$.  Thus $H_i \cong C_{q+1}$.

\begin{lemma}\label{l:H_i transitive}  For $i = 1,2$, the group $H_i$ acts simply transitively on $E_X(x_i)$. \end{lemma}

\begin{proof}  The action of $M_i$ on $E_X(x_i)$ can be identified with the natural action of $SL_2(q)$ on the projective line $\mathbb{P}^1(q)$, in which any one-point stabiliser has trivial intersection with $H_i$.  \end{proof}

Now let $A_0$ be any subgroup of $Z(G)$ and for $i = 1,2$ let $A_i = A_0 \times H_i$.  Then by Lemma~\ref{l:H_i transitive}, condition \eqref{i:n orbits} of Proposition~\ref{p:embedding} is satisfied with $n = 1$, and $\Stab_{A_1}(x_2) = A_1 \cap A_2 = A_0 = \Stab_{A_2}(x_1)$, so condition \eqref{i:representatives} of Proposition \ref{p:embedding} also holds.  Proposition \ref{p:embedding} then implies that $\G = A_1*_{A_0} A_2$ is an edge-transitive lattice in $G$.  This completes the proof of Theorem~\ref{t:cocompact existence}\eqref{i: p = 2}.

\subsubsection{Part \eqref{i:PSL} of Theorem \ref{t:cocompact existence}}\label{proof:PSL}

If $p$ is odd, the two possibilities for $L_i / Z(L_i)$ are $L_i / Z(L_i) \cong PSL_2(q)$ or $L_i / Z(L_i) \cong PGL_2(q)$.  We first consider the case $L_i / Z(L_i) \cong PSL_2(q)$.  Assume also that $q \equiv 3 \pmod 4$.  Then since either $M_i\cong SL_2(q)$ or  $M_i \cong PSL_2(q)$, $H_i$ is cyclic of order $|Z(M_i)|(q+1)/2$  and $N_{M_i}(H_i)$ is a group of order $|Z(M_i)|(q+1)$. Moreover, $N_{M_i}(H_i)\cap T=Z(M_i)$.

\begin{lemma}\label{l:N_i transitive}  For $i = 1,2$, the group $N_{M_i}(H_i)$ acts  transitively on $E_X(x_i)$. \end{lemma}

\begin{proof}  Again, the action of $M_i$ on $E_X(x_i)$ can be identified with the natural action of $SL_2(q)$ on the projective line. In particular, a stabiliser of an edge $\Stab_{M_i}(e)$ in $M_i$ has order $\frac{q(q-1)}{2}|Z(M_i)|$. Hence $N_{M_i}(H_i)\cap \Stab_{M_i}(e)\leq N_{M_i}(H_i)\cap T=Z(M_i)$. Thus the length of the orbit of $N_{M_i}(H_i)$ on $E_X(x_i)$ is $|Z(M_i)|(q+1)/|Z(M_i)|=q+1$ which proves the result.
\end{proof}

For $i=1,2$, consider $Z(M_i)$. Then $Z(M_i)\leq T_0$ where $T_0\in\mathcal{S}yl_2(T)$. Since $q\equiv 3\pmod 4$ while $L_i/Z(L_i)\cong PSL_2(q)$, $|T_0|$ divides $4$ and $\exp(T)\leq 2$. In fact, $[T_0,M_i]=1$ for $i=1,2$. Thus $T_0\leq C_G(\langle M_1, M_2\rangle)=C_G(\Lambda)\leq Z(G)$, and so $Z(M_i)\leq Z(G)$ for $i=1,2$.  Now take $A_0$ to be any subgroup of $Z(G)$ with $Z(M_i)\leq A_0$ for $i = 1,2$.  Let $A_i: = A_0N_{M_i}(H_i)$. Then $A_0\cap N_{M_i}(H_i)=Z(M_i)$.  Lemma~\ref{l:N_i transitive}  implies that $A_i$ is transitive on $E_X(x_i)$. Moreover, $\Stab_{A_i}(x_{3-i})=A_0$. Now an easy application of Proposition~\ref{p:embedding} finishes the proof of Theorem~\ref{t:cocompact existence}\eqref{i:PSL}.

\subsubsection{Part \eqref{i:PGL} of Theorem \ref{t:cocompact existence}}\label{proof:PGL}

To complete the proof of Theorem~\ref{t:cocompact existence}, we consider the case $L_i / Z(L_i) \cong PGL_2(q)$ (with $p$ odd).  Define $$N_i:=N_{L_i}(H_i) \quad \mbox{and} \quad C_i:=C_{L_i}(H_i).$$  Then $H_i\leq C_i$, $C_i$ is a cyclic group of order dividing $(q^2-1)$ (more precisely, $|C_i|=(q+1)|Z(L_i)|$), $Z(L_i)\leq C_i$ and $|N_i:C_i|=2$ with $N_i=C_iN_T(H_i)$.   For any group $C$ such that $Z(L_i) \leq C \leq L_i$, denote by $\overline{C}$ the image of $C$ in $\overline{L}_i:=L_i/Z(L_i)$.  Then $\overline{C}_i\cong C_{q+1}$.

Assume first that $q\equiv 3\pmod 4$, so that $\frac{q-1}{2}$ is odd.  Then $C_i=C'_i\times Z'_i$ where $(|C_i'|, |Z_i'|)=1$, $Z'_i\leq Z(L_i)$ is the ``odd part" of $Z(L_i)$  (i.e.,  $|Z'_i|$ divides $\frac{q-1}{2}$ and $Z(L_i)=Z(M_i)Z_i'$),  and $\overline{C'}_i=\overline{C}_i$. Moreover,  $|C_i':H_i| = 2$.   In fact, if $Q_i'\in\mathcal{S}yl_2(C_i)$,  $C_i'=H_iQ_i'=H_i'\times Q_i'$ where $H_i'=O(H_i)\cong C_{\frac{(q+1)|Z(M_i)|}{|Q_i'|}}$ and $Q_i'\cap Z(L_i)=Z(M_i)$.

\begin{lemma}\label{l:C_i' transitive} 
For $i = 1,2$ the group $C_i'=H_iQ_i'$ acts transitively on $E_X(x_i)$. 
\end{lemma}

\begin{proof}
This time the action of $L_i$ on the set $E_X(x_i)$ is the natural action of $PGL_2(q)$ on the projective line.  As we noticed earlier the image $\overline{C_i'}$ of $C_i'$ in $\overline{L_i}$ is isomorphic to $C_{q+1}$ while the image $\overline{H}_i$ of $H_i$ is isomorphic to $C_{\frac{q+1}{2}}$.  Moreover, there exists $\overline{c_i'}\in \overline{C_i'}$ such that $\overline{C_i'}=\langle \overline{c_i'}\rangle$ and $1\neq \overline{c_i'}^2\in \overline{H_i}$. As $H_i\cap Z(L_i)=Z(M_i)$, $\Stab_{\overline{H_i}}(e)=1$ for any $e\in E_X(x_i)$, and so $\Stab_{\overline{C_i'}}(e)=1$ for $e\in E_{X}(x_i)$. Hence the length of the orbit of $\overline{C_i'}$ on $E_X(x_i)$ is $q+1$. The result follows immediately.
 \end{proof}

Now let $T_0$ be the Sylow $2$--subgroup of $T$. Then as $q\equiv 3\pmod 4$ while $|T|$ divides $(q-1)^2$, we observe that $|T_0|=2|Z(M_i)|$, $\exp(T_0)=2$, $T_0\cap Q_i'=Z(M_i)$ and $T_0\leq N_{L_i}(H_i)$.  Take $Z_0$ to be any subgroup of $Z(G)$ and $A_0 = T_0 Z_0$.   Let $A_1 = C_1' A_0$ and $A_2 = C_2' A_0$.  Lemma~\ref{l:C_i' transitive} then implies that condition~\eqref{i:n orbits} of Proposition~\ref{p:embedding} holds, with $n = 1$.  We also have $A_1 \cap A_2 = A_0$, and so condition~\eqref{i:representatives} of Proposition~\ref{p:embedding} is satisfied.  We have proved Theorem~\ref{t:cocompact existence}\eqref{i:PGL q = 3}.

To prove Theorem \ref{t:cocompact existence}\eqref{i:PGL q = 3 other}, assume that  $Z(M_i)\leq Z(G)$ for $i=1,2$.  Take $A_0$ to be a subgroup of $Z(G)$ with $Z(M_i)\leq A_0$ for $i=1,2$.  Then if $A_1 = C_1'A_0$ and $A_2 = C_2' A_0$, we similarly obtain an edge-transitive lattice $\G = A_1 *_{A_0} A_2$.   Notice that $C_i'  \cap A_0 = Z(M_i)$.  We remark also that $Z(M_i)\leq Z(G)$ is crucial in this case, for it allows for the condition  $A_1\cap A_2 = \Stab_{A_i}(x_{3-i})$ to be satisfied.  This completes the proof of part \eqref{i:q = 3 mod 4} of Theorem \ref{t:cocompact existence}. 

Suppose now that $q\equiv 1\pmod 4$. In particular,  $\frac{q+1}{2}$ is odd.  Then $C_i=H'_i\times Q_i'\times  Z'_i$ where the following holds:
\begin{enumerate}
\item  $H'_i\cong C_{\frac{q+1}{2}}$ and $H_i=H_i'\times Z(M_i)$; 
\item $Q'_i=O_2(C_i)$, in particular,  as before $Q_i'\in\mathcal{S}yl_2(C_i)$, $Q_i'$ is cyclic and $Q_i:= Q_i'\cap Z(L_i)$ is in $\mathcal{S}yl_2(Z(L_i))$; and 
\item $Z'_i\leq Z(L_i)$ with $|Z'_i|$ dividing $\frac{q-1}{2}$,  i.e., $Z_i'$ is a cyclic subgroup of odd order and whose order is the part of $|C_i|$ coprime to $(q+1)$. 
\end{enumerate}

Notice that $H_i\cap Q_i'=Z(M_i)$ and so $H'_i\times Q_i'=H_iQ_i'$. 

\begin{lemma}\label{l:C_i'' transitive} For $i = 1,2$ the group $H_iQ_i'$ acts transitively on $E_X(x_i)$. \end{lemma}

\begin{proof}
As in the previous case, the action of $L_i$ on the set $E_X(x_i)$ is the natural action of $PGL_2(q)$ on the projective line.  Now the image $\overline{H_iQ_i'}$ of $H_iQ_i'$ in $L_i/Z(L_i)$ is isomorphic to $C_{q+1}$ while the image $\overline{H}_i$ of $H_i$ is isomorphic to $C_{\frac{q+1}{2}}$.  Moreover, $\overline{Q_i'}\cap \overline{T}=1$, hence $\Stab_{\overline{H_iQ_i'}}(e)=1$ for any $e\in E_X(x_i)$.  Another application of the Orbit-Stabiliser Theorem immediately gives the result.
 \end{proof}

Suppose now that $Q_i\leq Z(G)$. Take $A_i:=H_iQ_i'Z_0$ with $Z_0\leq Z(G)$. Then $A_i$ acts transitively on $E_X(x_i)$ by Lemma~\ref{l:C_i'' transitive}, and $A_1\cap A_2=Q_iZ_0\leq Z(G)$ with $Q_iZ_0=\Stab_{A_i}(x_{3-i})$. Hence, all conditions of Proposition \ref{p:embedding} are satisfied, which proves part~\eqref{i:PGL q = 1 other} of~Theorem \ref{t:cocompact existence}.

Now let us impose a slightly weaker condition: instead of  assuming $Q_i\leq Z(G)$, let us suppose that $Q_i^2\leq Z(G)$ where $Q_i^2$ is the unique subgroup of $Q_i$ of index $2$ (i.e., $Q_i^2=\langle x_i^2\mid \langle x_i\rangle=Q_i \rangle$).  Now  take $A_i:=O_{\left(\frac{q-1}{|Q_i'|}\right)'}(N_i)=H_iQ_i'Z_0\langle t_i\rangle$ where $t_i\in N_T(H_i)-C_T(H_i)$ is of order $2$ and $Z_0\leq Z(G)$.  Then by Lemma~\ref {l:C_i'' transitive}, $A_i$ acts transitively on $E_X(x_i)$, as it contains $H_iQ_i'$. 

Finally, let $T_0\in\mathcal{S}yl_2(T)$.  If $Q_i\leq Z(G)$, $A_1\cap A_2=Q_i\langle t_i\rangle Z_0=(Z(G)\cap A_i)\la t_i\ra=\Stab_{A_i}(x_{3-i})$ as $Q_i\langle t_i\rangle$ contains all elements of order $2$ in $T_0$, and so $t_1\in Q_2\la t_2\ra$ while  $t_2\in Q_1\la t_1\ra$.  Assume though that $Q_i\not\leq Z(G)$.  Since $Q_i^2\leq Z(G)$,  $Q_i\la t_i\ra=\{ x\in T_0\mid x^2\in Z(G)\}$ for $i=1,2$, and so $Q_1\la t_1\ra= Q_2\la t_2\ra$. Therefore $A_1\cap A_2=Q_i\langle t_i\rangle Z_0=\Stab_{A_i}(x_{3-i})$ which finishes the proof of Theorem \ref{t:cocompact existence}.

\subsubsection{Proof of Theorem \ref{t:exceptions}}\label{proof:exceptions}

We finally construct the exceptional edge-transitive lattices in Theorem \ref{t:exceptions}.  We provide only a brief discussion, since our proof consists of carrying out for Kac--Moody groups $G$ as in Theorem \ref{t:cocompact existence} the constructions used by Lubotzky to prove \cite[Theorem 3.3]{L} for $SL_2(\Fqt)$.

The reason for the existence of these ``exceptional" lattices is that when $q$ is small enough, $L_i$ occasionally contains a subgroup $N_i$ distinct from the ones listed in the conclusions of Theorem \ref{t:cocompact existence}, yet still acting transitively on the set $E_X(x_i)$. In those cases we then take $A_i=N_iA_0$ with $A_0\leq N_T(N_i)$, $i=1,2$, and check that both conditions of  Proposition \ref{p:embedding} are satisfied with $n=1$. 

To determine $N_i$, we use the fact that $L_i$ acts on $E_X(x_i)$ as $\overline{L_i}=L_i/Z(L_i)$ (which is either $PSL_2(q)$ or $PGL_2(q)$) naturally acts on $\mathbb{P}^1(q)$. Thus we may as well look for subgroups $\overline{N_i}$ of $\overline{L_i}$ that act transitively on the points of $\mathbb{P}^1(q)$ and are different from $N_{\overline{L_i}}(\overline{H_i})$.  Since we are interested in constructing lattices without $p$--elements, we will also assume that $(q,|N_i|)=1$.  The results of Dickson together with the fact that the normaliser of a split torus of $\overline{L_i}$ is not transitive on the points of $\mathbb{P}^1(q)$, tell us immediately that $q$ is odd and $\overline{N}_i$ can only be isomorphic to one of the following groups: $A_4$, $S_4$ or $A_5$.
On the other hand, as $\overline{N}_i$ acts transitively on $E_X(x_i)$, the Orbit-Stabiliser Theorem implies that $|E_X(x_i)|=q+1$ divides $|\overline{N}_i|$.  Combining these two easy arguments gives us $q+1\leq 60$. In fact, when $(5, |L_i|)=1$, we have $q+1\leq 24$.

Thus we first list all the odd prime powers $q$ that are less  or equal to $59$. If $(5,q^2-1)=1$, then $q\leq 23$.  This leaves us with $5,7, 9, 11, 13, 17, 19, 23, 25, 29, 31, 41, 49$ and $59$.  Furthermore,  $q+1$ divides $|\overline{N_i}|$ and as $|\overline{N_i}|$ divides $2^3\cdot 3\cdot 5$,
we may exclude $13$, $17$, $25$, $31$, $41$ and $49$. Finally, if $q=9$, $q+1=10$, and so $\overline{N_i}\cong A_5$, which contradicts the fact that $(q,|N_i|)=1$.  We are now down to the same list of exceptional $q$ as in \cite{L}.  Hence, all that remains to do is for a given prime $q$ take a candidate for $\overline{N_i}$ for $i=1,2$ as in Theorem \ref{t:exceptions}, and check that the conditions of Proposition \ref{p:embedding} are satisfied with $n=1$. This follows immediately  as in \cite{L} by simply looking at the action of the appropriately chosen $N_i$ on the points of $\mathbb{P}^1(q)$. 

\subsection{Construction of the lattice $\mathbf{\Gamma'}$}\label{s:construction}

Let $G$ be as in Theorem~\ref{t:cocompact existence} and assume that we are \emph{not} in any of the cases in Theorem~\ref{t:cocompact existence} where $G$ admits an edge-transitive lattice.  We now construct a cocompact lattice $\G' < G$ which is the fundamental group of the graph of groups $\mathbb{A}'$ sketched in the statement of Theorem~\ref{t:cocompact existence}.  For this, we first in Section~\ref{ss:define groups} define certain finite subgroups $S$, $N_1$ and $N_2$ of $G$ and discuss their structure, then in Section~\ref{ss:applying embedding} construct $\G'$ by verifying that our embedding criterion, Proposition \ref{p:embedding} above, may be applied with $n = 2$, $A_1 = N_1$ and $A_2 = N_2$.

\subsubsection{The groups $S$, $N_1$ and $N_2$}\label{ss:define groups}

Let $P_1$ and $P_2$ be the standard maximal parahoric subgroups of $G$.  Then for $i = 1,2$, $P_i$ is the stabiliser in $G$ of a vertex $x_i$ of $X$ with $[x_1,x_2]$ an edge of $X$.  Recall that since $G$ is rank $2$ and has symmetric Cartan matrix, $P_1\cong P_2$. Moreover,   if $L_i$ is a Levi complement  of $P_i$, then $L_i=M_iT$ where $T\leq B\leq P_1\cap P_2$ is a torus of $G$ and $A_1(q)\cong M_i\triangleleft L_i$, where $A_1(q)$ is isomorphic to either $SL_2(q)$ or $PSL_2(q)$, depending upon $G$.  By assumption, $q \equiv 1 \pmod 4$, and either $L_i/Z(L_i)\cong PSL_2(q)$, or  if $L_i/Z(L_i)\cong PGL_2(q)$, $Q_i^2\not\leq Z(G)$ where $Q_i^2$ is the unique subgroup of index $2$ of the Sylow $2$--subgroup of $Z(L_i)$.

For $i = 1,2$ let $H_i$ be a fixed non-split torus of $M_i$ such that $N_T(H_i)$ is as big as possible. Then either $H_i\cong C_{\frac{q+1}{2}}$ or $C_{q+1}$, depending on whether $M_i\cong PSL_2(q)$ or $SL_2(q)$ respectively.  Also,  $N_T(H_i)/C_T(H_i)\cong C_2$ and $H_i\cap N_T(H_i)=Z(M_i)$.  Define \[S:=N_{T}(H_1)\cap N_T(H_2).\] Let us try to describe $S$ in more definite terms.  Let $Q$ be the Sylow $2$--subgroup of $T$ (it is unique since $T$ is abelian).  First, notice that if  $z\in N_T(H_i)$, $i=1,2$, is of odd order, then $[z,H_i]=1$ and $[z,M_i]=1$. Hence if $z\in S$ and $z$ is of odd order, $z\in C_G(\langle M_1, M_2\rangle)$ thus $z\in Z(G)$. It follows immediately that $Z(G)\leq S\leq Z(G)Q$.  Let us now investigate what happens when $z\in N_T(H_i)\cap Q$, $i=1,2$.

Take $x\in Q$  such that $x$ normalises but not centralises $H_i$ for some $i\in\{1, 2\}$. Then $x$ acts on $H_i$ as an element of order $2$, and so $x^2$ centralises $H_i$. It follows that $x^2$ centralises $M_i$.  Now consider $R_i:=\{x\in Q \mid x^2\in C_T(M_i)\}$. Then $ R_i\leq Q$ and $R_i\leq N_T(H_i)$.  Define $$Q_0:=R_1\cap R_2=\{x\in Q \mid x^2\in C_T(M_1)\cap C_T(M_2)\}=\{x\in Q \mid x^2\in Z(G)\}.$$  Clearly, $Q_0\leq S$.  On the other hand, take $s\in S\cap Q$. If $[s,H_i]=1$ for both $i=1,2$, then $[s, M_i]=1$ for $i=1,2$ implying $s\in Z(G)\cap Q\leq Q_0$.  Let $s\in S\cap Q$ be such that $[s,H_i]\neq 1$ for some $i\in\{ 1, 2\}$. As noticed above,  $s^2\in C_T(H_i)\leq C_T(M_i)$.  Hence, $s^2\in C_T(M_j)$, $\{ i, j\}=\{1, 2\}$. Therefore, $s^2\in Z(G)$. Thus $S\cap Q\leq Q_0$.  It follows that:

\begin{lemma}\label{l:S} $S=Z(G)Q_0$. \end{lemma}

Notice that $|N_S(H_i):C_S(H_i)|=2$. We also define \[ N_1: = SH_1 \quad \mbox{and} \quad N_2: = SH_2.\] 

\subsubsection{Application of embedding criterion}\label{ss:applying embedding}

By construction, for $i = 1,2$, $N_i$ is a finite subgroup of $P_i$, and $S = N_1 \cap N_2$.  We now verify that our embedding criterion, Proposition \ref{p:embedding} above, may be applied with $n = 2$, $A_1 = N_1$ and $A_2 = N_2$.  

Notice first that for $i = 1,2$, the intersection of $N_i$ with an edge stabiliser in $L_i$ is of index $\frac{q+1}{2}$.  The Orbit-Stabiliser Theorem yields immediately  that $N_i$ has $2$ orbits of equal size $\frac{1}{2}(q + 1)$ on $E_X(x_i)$.    That is, with $n = 2$, condition \eqref{i:n orbits} in the statement of Proposition \ref{p:embedding} above holds. 

Denote by $f_1$ the edge $[x_1,x_2]$ of $X$ and by $\hat{f}_1$ the edge $[x_2,x_1]$.  Choose an edge $f_{2} \in E_X(x_1)$ so that the edges $f_{1}$ and $f_{2}$ represent the two $N_1$--orbits on $E_X(x_1)$, and choose an edge $\hat{f}_{2} \in E_X(x_2)$ so that the edges $\hat{f}_{1}$ and $\hat{f}_{2}$ represent the two $N_2$--orbits on $E_X(x_2)$.  

The edges $f_{1}$ and $\hat{f}_{1}$ are fixed by $S$, since $S \leq T \leq B = P_1 \cap P_2$.  We claim that the edges $f_{2}$ and $\hat{f}_{2}$ may be chosen so that $S$ fixes both $f_{2}$ and $\hat{f}_{2}$.  To see this, consider first the action of $N_1$ on the edges $E_X(x_1)$.  Now $N_1\leq L_1$, and $L_1$ acts on the set $E_X(x_1)$ as on the points of projective line, that is, we observe this action via a homomorphism $\phi:L_1\rightarrow PGL_2(q)$. The kernel of this action is $\ker(\phi)=Z(L_1)=C_T(M_1)$.  We know that $N_1$ has $2$ orbits, say $\theta_1$ and $\theta_2$, in this action, each of length $\frac{q+1}{2}$ which is odd.  Assume that the fixed points of $S$ all lie inside the same orbit of $N_1$, say $\theta_1$.  Then $S$ would act fixed-point free on $\theta_2$. Now, $S\ker(\phi)/\ker(\phi)\cong S/S\cap\ker(\phi)$ and as $|S\ker(\phi)/\ker(\phi)|=2$, $S$  would have a fixed point on $\theta_2$, a contradiction. Hence we may choose the edge $f_{2} \in E_X(x_1)$ so that $f_{2}$ is fixed by $S$.  Similarly, we may choose $\hat{f}_{2} \in E_X(x_2)$ to be fixed by $S$. 

Now let $g_{1} = \hat{g}_{1} = 1_G$.  Consider the fixed points of $S$ on $E_X(x_i)$, $i=1,2$.  Since $|S\cap M_i:Z(M_i)|=2$, they are the two  points fixed by the whole of $T$.  Choose $g_2\in N_{P_1}(T)$  that  represents $w_1\in W$.  Then $g_2\cdot f_1=f_2$.  Similarly, we may choose $\hat{g_2}\in N_{P_2}(T)$ that represents $w_2$ and such that $\hat{g}_2\cdot\hat{f}_1=\hat{f}_2$.  Then (2a) in Proposition \ref{p:embedding} above holds.  Let $\tau :=g_{2}\hat{g}_{2}^{-1}$.  We observe that $S^\tau = S$, since by Lemma~\ref{l:S} above $S=Z(G)Q_0$, a characteristic subgroup of $T$ which is therefore $N_G(T)$--invariant.  Hence $(N_1 \cap N_2)^{g_2} = (N_1 \cap N_2)^{\hat{g}_2}$, and so \eqref{i:g on edges} in Proposition \ref{p:embedding} above is satisfied.
 
To show that \eqref{i:edge groups} in Proposition \ref{p:embedding} above holds, we must show that $N_1 \cap B =  N_1 \cap B^{g_2} = N_2 \cap B = N_2 \cap B^{\hat{g}_2} = S$.  Since $N_1 \leq L_1 \leq P_1$, we have that $N_1 \cap B = N_1 \cap (B \cap L_1)$.  Now, $B \cap L_1$ is isomorphic to a Borel subgroup $TU_0$ of $L_1$, where $U_0\cong E_{p^a}$, the elementary abelian group of exponent $p$ and order $q=p^a$, is normalised by $T$. On the other hand $N_1$ is a finite subgroup of $L_1$. The order of $N_1$ is $|S|\frac{q+1}{2}$ and it divides $|T|\frac{q+1}{2}$. Moreover, $(|S|,\frac{q+1}{2})=1$. Therefore, numerical reasons imply that $N_1 \cap B$ is a finite group whose order divides $|T|$ and is actually at most $|S|$. But $S\leq N_1$ and $S\leq T\leq B$.  Hence $N_1 \cap B = S$ as required.  The argument that $N_2 \cap B = S$ is similar.

Since $S$ fixes the edge $f_2 = g_2 \cdot f_1$, we have $S \leq B^{g_2}$.  The argument that $N_1 \cap B^{g_2} = S$ is then similar to the previous paragraph.  Finally, $S$ also fixes the edge $\hat{f}_2 = \hat{g}_2 \cdot \hat{f}_1$, and again by similar arguments we conclude that $N_2 \cap B^{\hat{g}_2} = S$.  Therefore all hypotheses of Proposition \ref{p:embedding} above are satisfied with $n = 2$, $A_1 = N_1$ and $A_2 = N_2$, and so the fundamental group $\Gamma'$ of the graph of groups $\mathbb{A}$ as sketched in the statement of Theorem \ref{t:cocompact existence} is a cocompact lattice in $G$ with quotient the graph $A$.

\section{Discussion of conjectures}\label{s:motivation}

Let $G$ be as in Theorem \ref{t:cocompact existence}.  In this section we motivate Conjectures \ref{c:cocompact} and \ref{c:unique end}, stated in the introduction, and explain why Conjecture \ref{c:unique end} implies Conjecture \ref{c:cocompact}.   

We will use the following well-known property of cocompact lattices.

\begin{prop}[see p. 10 of Gelfand--Graev--Piatetsky-Shapiro~\cite{GGPS}] \label{t:gelfand} Let $G$ be any locally compact topological group,  and $\Gamma$ a cocompact lattice in $G$.  If $u\in\Gamma$, then $$u^G := \{ gug^{-1} \mid g \in G\}$$ is a closed subset of $G$. \end{prop}

\noindent We will also need the next result, which follows from the reference Ronan \cite[Theorem 6.15]{Ronan}. 

\begin{lemma}\label{l:fixing balls} Let $\alpha$ be a positive real root whose corresponding geodesic ray in $\Sigma$ begins at a vertex distance at least $n+1$ from both $P_1$ and $P_2$.  Then $U_\alpha$ fixes the ball of radius $n$ about the edge $B$.\end{lemma}

Recall the definition of the closed abelian group $\cU$ from Section \ref{s:action ends} above.  Observe that every torsion element of $\cU$ has order $p$.  We now show that nontrivial elements of conjugates of $\mathcal{U}$ cannot be contained in cocompact lattices.  It is then immediate that Conjecture \ref{c:unique end} implies Conjecture \ref{c:cocompact}.

\begin{corollary}\label{c:contraction}  Let $u$ be a nontrivial element of any $G$--conjugate of $\mathcal{U}$.  Then there is a $g \in G$ such that  \[ \lim_{n \to \infty } g^nug^{-n} = 1_G.\]  Hence by Proposition \ref{t:gelfand} above, $u$ is not contained in any cocompact lattice $\G < G$. \end{corollary}

\begin{proof}  It suffices to prove the result for $u \in \mathcal{U}$.  Let $g_\tau \in G$ induce the translation $\tau = w_1w_2$, as defined in Section \ref{s:action ends} above, and for each positive integer $n$ let $f_n:G  \to G$ be given by $f_n(x) = g_\tau^n x g_\tau^{-n}$.  We will show that $f_n(u) \to 1_G$.

By definition, $u$ is the limit in the building topology of a sequence of elements $u_k$ in the elementary abelian $p$--group $V_1 \times -V_2$.  Since $u$ fixes the end $\ep_1$, there is a root $\alpha \in \Phi^1_+ \cup \Phi^2_-$ such that $u$ fixes pointwise the half-line $\alpha$.  Now there are at most finitely many root groups in $-V_2$ which fix such an $\alpha$, namely those root groups $U_\beta$ with $\alpha \subseteq \beta \subseteq -\alpha_2$.   We may thus assume without loss of generality that for all $u_k$, all of the root groups $U_\beta$ in $-V_2$ on which $u_k$ has a nontrivial projection satisfy $\alpha \subseteq \beta \subseteq -\alpha_2$.  Hence for all $n$ large enough, $f_n(u) \in \hat{V_1}$.  So we reduce to the case that $u \in \hat{V}_1$.

For each root $\alpha \in \Phi^1_+$, consider the group \[f_n(U_\alpha) = g_\tau^{n}U_\alpha g_\tau^{-n} = U_{\tau^{n}\alpha}.\]  Since $\alpha \in \Phi^1_+$ and $\tau$ acts by translation by two edges with repelling fixed point $\ep_1$, the distance from the initial vertex of the ray $\tau^{n}\alpha$ to the edge $B$ is at least $2n \geq n + 1$.  Lemma \ref{l:fixing balls} then implies that $f_n(U_\alpha)$ fixes pointwise the ball of radius $n$ about $B$. Therefore $f_n(V_1)$ fixes this ball pointwise.  By definition of the building topology, it follows that $f_n(\hat{V}_1) = g_\tau^{n}\hat{V}_1g_\tau^{-n}$ fixes pointwise the ball of radius $n$ about $B$.  Hence for all $u \in \hat{V_1}$, $f_n(u) \to 1_G$, as required. \end{proof}

\section{Classification of edge-transitive lattices}\label{s:classification}

Let $G$ be as in Theorem \ref{t:cocompact existence}.  In this section we assume Conjecture \ref{c:cocompact}, that is, that cocompact lattices in $G$ do not contain $p$--elements, and under this assumption prove Theorem \ref{t:classification}, which says that the edge-transitive lattices in $G$ are exactly those described in Theorems \ref{t:cocompact existence} and \ref{t:exceptions}. 

To fix notation throughout this section, for $i = 1,2$ let $x_i$ be the vertex of the tree $X$ which is fixed by the standard parahoric $P_i$.  We have that $P_i$ has Levi decomposition $P_i = L_i \ltimes U_i$ (see Proposition~\ref{p:completion}) where $L_i = T_i M_i$.

In order to classify the edge-transitive lattices in $G$, we will rely on the following converse to the case $n = 1$ of Proposition~\ref{p:embedding}, which follows from Bass--Serre theory and the discussion of lattices in Section \ref{s:lattices}.  Our notation is the same as in Proposition~\ref{p:embedding}.

\begin{lemma}\label{l:converse n=1} Suppose $\G$ is a cocompact lattice in $G$ with fundamental domain the edge $[x_1,x_2]$.  For $i = 1,2$, let $A_i = \Stab_\G(x_i)$.  Then $\G \cong A_1 *_{A_1 \cap A_2} A_2$ and $A_1$ and $A_2$ are finite subgroups of $G$ such that for $i = 1,2$, 
\begin{enumerate} 
\item\label{i:transitive} $A_i$ acts transitively on $E_X(x_i)$; and
\item\label{i:stab} $\Stab_{A_i}(x_{3 - i}) = A_1 \cap A_2$.
\end{enumerate}
\end{lemma}
  
Now suppose $\G$ is a cocompact edge-transitive lattice in $G$.   Our classification of edge-transitive lattices is up to isomorphism, so by Lemma~\ref{l:converse n=1} we may assume without loss of generality that $\G \cong A_1 *_{A_0} A_2$, where for $i = 1,2$ the group $A_i$ is a finite subgroup of the standard parahoric $P_i$, and $A_0 = A_1 \cap A_2$.  By Conjecture \ref{c:cocompact}, the groups $A_1$ and $A_2$ do not contain any $p$--elements.  Also, since the finite group $A_0$ has index $(q+1)$ in both $A_1$ and $A_2$, for $i = 1,2$ the group $A_i$ has order divisible by $(q + 1)$.  

We continue by listing the groups $A_1$ and $A_2$ which satisfy the restrictions stated so far.  Recall from Section~\ref{s:kac-moody} that if $L = L_i$ is the Levi factor of a standard maximal parahoric subgroup $P = P_i$ of $G$, then $L=TM$ where $T$ is the torus of $G$ and $A_1(q)\cong M\triangleleft L$.

\begin{prop}\label{p:char} Let $A = A_i$ be a finite subgroup of $P = P_i$ such that $|A|$ is coprime to $p$ and is divisible by $(q+1)$.  We may assume without loss of generality that $A$ is a subgroup of the Levi factor $L$.  Then one of the following conditions hold:
\begin{enumerate}
\item\label{i:p = 2 X} If $p=2$, then $H\leq A\leq C_T(M)H$ where $H\cong C_{q+1}$ is a non-split torus of $M\cong SL_2(q)$.
\item\label{i:p odd X} If $p$ is odd and $A\leq MZ(L)$, then $A\cong T_0H$ where $T_0\leq N_T(H)$ and $H\leq M$. More precisely, if $A_1(q)$ is universal  (that is, if $M\cong SL_2(q)$), then $H$ is isomorphic to a subgroup  listed in the conclusions to Corollary~\ref{c:Dickson} above.  Otherwise, $M\cong PSL_2(q)$ and  $H\cong H'/\langle -I\rangle$ where $H'$ is a conclusion to Corollary~\ref{c:Dickson} above. But if $p=3$ and $q=9$, $H/Z(H)\not\cong A_5$.
\item\label{i:p odd X last case} If $p$ is odd, $L/Z(L)\cong PGL_2(q)$ and $A\not\leq MZ(L)$, then  either $AZ(L)/Z(L)\leq H_0$ where  $H_0\cong D_{2(q+1)}$, or $q\in\{5, 11\}$ and $AZ(L)/Z(L)\cong S_4$.
\end{enumerate}
\end{prop}

\begin{proof}  Since $P = P_i$ has Levi decomposition $P = L \ltimes U$ where $U$ is pro--$p$ (see Proposition \ref{p:completion}), we may assume without loss of generality that $A \leq L$.  
Parts \eqref{i:p = 2 X} and \eqref{i:p odd X} then follow immediately from Corollary~\ref{c:Dickson} above. 
 Part \eqref{i:p odd X last case} follows from Theorem~\ref{t:Dickson2} above and the fact that subgroups isomorphic to $D_{q+1}$ from Part \eqref{i:Dickson dihedral} of Theorem ~\ref{t:Dickson2} 
 are either contained in subgroups isomorphic to $D_{2(q+1)}$ from Part \eqref{i:Dickson dihedral 2} of the same result, or possibly $q\equiv\pm 3 \pmod 8$ and $A\cong S_4$. Now \cite{Atlas} confirms the result. \end{proof}

The remaining proof of Theorem \ref{t:classification} is divided into the same cases as in the statement of Theorem \ref{t:cocompact existence}.

\subsection{Case \eqref{i: p = 2} of Theorem \ref{t:cocompact existence}}

We first prove Theorem~\ref{t:classification} in the case $p = 2$, by showing that every edge-transitive lattice in $G$ has the form described in Theorem \ref{t:cocompact existence}\eqref{i: p = 2}.  By Lemma~\ref{l:converse n=1} and the discussion before Proposition~\ref{p:char}, any edge-transitive lattice $\G < G$ satisfies $\G \cong A_1 *_{A_1 \cap A_2} A_2$ where for $i = 1,2$, $A_i$ is a finite subgroup of $P_i$ of order divisible by $(q+1)$ and coprime to $p$.  But then by  Proposition~\ref{p:char}\eqref{i:p = 2 X}, for $i = 1,2$ we have $H_i \leq A_i \leq C_{T}(M_i)H_i$ where $H_i \cong C_{q+1}$ is a non-split torus of $M_i$.  (Note that since $H_i$ acts transitively on $E_X(x_i)$, Lemma~\ref{l:converse n=1}\eqref{i:transitive} automatically holds.)  Now observe that $$A_0 := A_1\cap A_2\leq C_T(M_1)\cap C_T(M_2)\leq C_T(\langle M_1, M_2\rangle)\leq Z(G),$$ as required.  By Lemma~\ref{l:converse n=1}\eqref{i:stab}, it follows that $A_i = A_0H_i$, which completes the proof of Theorem \ref{t:classification} in this case.  

\subsection{Case \eqref{i:PSL} of Theorem \ref{t:cocompact existence}}  

We are now in the case that $p$ is odd and $L_i/Z(L_i) \cong PSL_2(q)$.  To prove Theorem \ref{t:classification} in this case, assume first that $q \equiv 3 \pmod 4$.  As in the proof of Theorem~ \ref{t:cocompact existence}, we may conclude that in this case $Z(M_i)\leq Z(G)$ for $i=1,2$.  Then for $i = 1,2$, we are considering subgroups  $A_i$ of $M_i$ whose order  in $M_i/Z(M_i)\cong PSL_2(q)$ is divisible by $(q+1)$.  By Proposition~\ref{p:char}\eqref{i:p odd X}, either $A_i=A_0N_{M_i}(H_i)$ where $H_i$ is a non-split torus of $M_i$ and $A_0\leq N_T(H_i)$, or $A_i=A_0N_i$ where $N_i\cong H$ above and $A_0\leq N_T(N_i)$.  Then in all cases $A_0\leq C_T(M_i)$ for $i=1,2$, and so $A_0\leq C_G(\langle M_1, M_2\rangle)\leq Z(G)$ and $A_0\cap N_{M_i}=Z(M_i)$ as required.  

We now assume that $q\equiv 1 \pmod 4$.  

\begin{lemma} The group $A_i=A_0N_{M_i}(H_i)$ does not act transitively on $E_X(x_i)$.  \end{lemma}

\begin{proof}
To see this, in its action on $E_X(x_i)$,  $A_i=A_0N_{M_i}(H_i)$ intersects a one-point stabiliser $B_i$ of $L_i$ in a subgroup of index $4/d$, that is, $$|A_i:A_i\cap B_i|=|A_0N_{M_i}(H_i):A_0N_{M_i}(H_i)\cap T\Stab_{M_i}(x_{3-i})|=|N_{M_i}(H_i):N_{M_i}(H_i)\cap T|=4/d,$$ where $d=2$ if $M_i\cong SL_2(q)$ and $d=1$ if $M_i\cong PSL_2(q)$.  Hence the length of the orbit of $A_i$ in its action on $E_X(x_i)$ is at most $\frac{2(q+1)/d}{4/d}=\frac{q+1}{2}$. And so $A_i$ is not transitive on $E_X(x_i)$. \end{proof}

Now Lemma~\ref{l:converse n=1}\eqref{i:stab} implies that $G$ does not contain edge-transitive cocompact lattices unless possibly one of the following holds: $q=5$ and $B_1\cong B_2\cong A_1(3)$, or  $q=29$ and $B_1\cong B_2\cong A_1(5)$.  These exceptional cases are constructed in Theorem \ref{t:exceptions}.

\subsection{Case \eqref{i:PGL} of Theorem \ref{t:cocompact existence}}

We are now in the case that $p$ is odd and $L_i/Z(L_i)\cong PGL_2(q)$.  We will use the notation introduced in the proof of Theorem~\ref{t:cocompact existence}. Recall that $C_i=C_{L_i}(H_i)$, $H_i\leq C_i$ and  $C_i$ is a cyclic group of order $(q+1)|Z(L_i)|$, $Z(L_i)\leq C_i$ and
$|N_i:C_i|=2$. Moreover, if $\overline{C}_i$ is the image of $C_i$ in $\overline{L}_i=L_i/Z(L_i)$,  $\overline{C}_i\cong C_{q+1}$ and $N_i=C_iN_T(H_i)$.  Using part \eqref{i:p odd X last case} of Proposition~\ref{p:char} we obtain that if $A_i$ is a subgroup we are looking for,  $\overline{A_i}=A_iZ(L_i)/Z(L_i)$ must be isomorphic to a subgroup of $D_{2(q+1)}$ whose order is divisible by $(q+1)$ or if $q\in\{5, 11\}$, $\overline{A_i}\cong S_4$.
Assume now that $\Gamma$ is not an exceptional lattice, i.e., we are not in the latter case.  Then $A_i$ is a subgroup of $N_i:=N_{L_i}(H_i)$. This is so because $N_i$ is the full pre-image of $D_{2(q+1)}$ in $L_i$.

Suppose first that $q\equiv 3\pmod 4$. In terms of  notation from the proof of Theorem \ref{t:cocompact existence}, $C_i=C'_i\times Z'_i$.  Suppose first that $A_i\leq N_i$ is such that $\overline{A}_i\cong D_{2(q+1)}$.  Then $A_i\geq C_i'$ and $A_i\geq T_0$. This is true as $T_0\leq N_{T}(H_i)$,
$\overline{C'_iT_0}=\overline{N}_i$ and if a subgroup of $N_i$ does not contain $C'_iT_0$, its image in $L_i$ would be properly contained in $D_{2(q+1)}$.  On the other hand, $A_i\leq N_i=C'_iT_0Z'_i$. Now, $A_0$ is a subgroup of $A_i$ of index $(q+1)$, and as $A_0=A_1\cap A_2\leq P_1\cap P_2=B$ while $(|A_0|,p)=1$, $A_0\leq T$.  It follows that $A_0\leq T_0Z'_i$.  Notice that $g\in A_0$ is of odd order if and only if $g\in Z'_i$.  Hence, $g\in O(A_0)$ implies $g\in Z(L_1)\cap Z(L_2)\leq C_G(\langle M_1, M_2\rangle)\leq Z(G)$.  We may now conclude that $A_0\leq T_0Z(G)$ which gives us the lattices from \eqref{i:PGL q = 3} of Theorem~\ref{t:cocompact existence}.

Let us now investigate which subgroups of $N_i$ will give us other edge-transitive lattices in $G$.  We are  looking for $A_i\leq N_i$ that acts edge-transitively on $E_X(x_i)$ and is such that $(q+1)$ dividing $|A_i|$.  Suppose first that $A_i$ is a subgroup of $C_i$ such that
$\overline{A}_i\cong C_{q+1}$.  Then $C'_i\leq A_i\leq C_i'Z'_i$. Recall that $A_0$ is a subgroup of $A_i$ of index $(q+1)$, and as before one can easily see that $A_0\leq T$.  Now, $C_i'Z_i'\cap T=Z(M_i)Z_i'$ and $Z(M_i)\leq C_i'$. Again, $g\in A_0$ is of odd order if and only if $g\in Z_i'$, implying $g\in Z(G)$.  Assume now that $Z(M_i)\neq 1$. It follows that $A_0=Z(M_1)(A_0\cap Z'_1)=Z(M_2)(A_0\cap Z'_2)$.  Because the Sylow $2$--subgroup of $A_0$ is normal in $A_0$, it follows that $Z(M_1)=Z(M_2)$, and so  $Z(M_i)\leq Z(G)$ for $i=1,2$.  If $Z(M_i)\not\leq Z(G)$, we will obtain that condition \eqref{i:stab} of Lemma~\ref{l:converse n=1} is violated.  Therefore, if $Z(M_i)\leq Z(G)$, we obtain some lattices from 
Theorem~\ref{t:cocompact existence}\eqref{i:PGL q = 3 other}.  Finally, the only other subgroups of $N_i$ that act transitively on $E_X(x_i)$ are subgroups that map onto $\overline{N_{M_i}(H_i)}\cong D_{q+1}$. These are the subgroups that map into $\overline{M_i}$ and we already described them in \eqref{i:PSL} of Theorem \ref{t:cocompact existence}: $A_i=N_{M_i}(H_i)Z_0$ where $Z_0\leq Z(G)$. However, $N_{M_1}(H_1)Z_0\cap N_{M_2}(H_2)Z_0= Z(M_i)Z_0$ for $i=1,2$, and it follows immediately that $Z(M_1)=Z(M_2)\leq Z(G)$, giving  the remaining lattices from \eqref{i:q = 3 mod 4} of Theorem \ref{t:cocompact existence}.

Suppose now that $q\equiv 1\pmod 4$.  This time $C_i=H_iQ_i'\times Z'_i$.  Suppose first that $A_i\leq N_i$ is such that $\overline{A}_i\cong D_{2(q+1)}$.  Then $A_i\geq H_iQ_i'\langle t_i\rangle$ where $t_i$ is an involution in $N_{T}(H_i)-C_T(H_i)$.  This holds because if a subgroup of $N_i$ does not contain $H_iQ_i'\langle t_i\rangle$, its image in $L_i$ would be properly contained in $D_{2(q+1)}$.  It follows that $A_i\leq H_iQ_i'\langle t_i\rangle O(Z'_i)$.  As in the previous case, $A_0$ is a subgroup of $A_i$ of index $(q+1)$ such that $(|A_0|,p)=1$ and $A_0\leq T$.  It follows that $A_0\leq (Q_i'\cap T)\langle t_i\rangle O(Z'_i)$.  As before $g\in A_0$ is of odd order  if and only if $g\in O(Z'_i)$, implying $g\in Z(G)$.  We may  now conclude that $A_0\leq (Q_i'\cap T)\langle t_i\rangle Z(G)$.  Now, $Q_i\leq A_i$ and in fact $Q_i\leq A_0$. Hence, the Sylow $2$--subgroup of $A_0$ is $Q_i\langle t_i\rangle$. It follows that $Q_1\langle t_1\rangle=Q_2\langle t_2\rangle$. 
%Let $2^b$ be the exponent of $Q_i$, $i=1,2$.  Denote by $Q^2_i:=\Omega_{b-1}(Q_i)$ 
As before let $Q_i^2$ be the subgroup of $Q_i$ generated by the squares of its elements.  It follows that $Q^2_1=Q^2_2$, and so $Q^2_i\leq Z(L_1)\cap Z(L_2)\leq Z(G)$.  If this condition fails, it is clearly impossible to have an edge-transitive lattice, while if it holds, the conditions of Proposition \ref{p:embedding} are satisfied. We therefore have obtained the lattices from \eqref{i:PGL q = 1} of Theorem \ref{t:cocompact existence}.

Let us now investigate which subgroups of $N_i$ will give us other edge-transitive lattices in $G$.  Again we are in quest of $A_i\leq N_i$ that acts edge-transitively on $E_X(x_i)$ (in particular,  with $(q+1)$ dividing $|A_i|$).  Suppose first that $A_i$ is a subgroup of $C_i$ such that $\overline{A}_i\cong C_{q+1}$.  Then $Q_i'H_i\leq A_i\leq C_iH_i$. Now $A_0$ is a subgroup of $A_i$ of index $(q+1)$, and as before it follows that $A_0\leq T$. Thus $|Q_i':Q_i'\cap A_0|=2$.  Again, $g\in A_0$ is of odd order if and only if $g\in Z_i'$, implying $g\in Z(G)$.  Hence, $A_0=(Q_i'\cap A_0)(A_0\cap Z'_i)=Q_i(A_0\cap Z'_i)$.  Because the Sylow $2$--subgroup of $A_0$ is normal in $A_0$, $Q_1=Q_2$.  In particular, the usual argument implies $Q_i\leq Z(G)$ for $i=1,2$.  If $Q_i\not\leq Z(G)$, we will obtain that condition \eqref{i:stab} of Lemma~\ref{l:converse n=1} is violated.  Therefore we obtain the remaining conclusions  of Theorem~\ref{t:cocompact existence}.  Since there are no other subgroups of $N_i$ that act transitively on
$E_X(x_i)$, we are done with the proof of Theorem~\ref{t:classification}.

\section{Minimality of covolumes among cocompact lattices}\label{s:covolumes cocompact}

Let $G$ be as in Theorem \ref{t:cocompact existence}.  In this section we assume Conjecture \ref{c:cocompact}, that is, that cocompact lattices in $G$ do not contain $p$--elements, and we prove Theorem~\ref{t:covolumes}, stated in the introduction, on the minimal covolume of cocompact lattices in $G$.  

We will subdivide our proof of Theorem \ref{t:covolumes} between the cases where $G$ admits an edge-transitive lattice, discussed in Sections \ref{s:minimality edge-transitive}--\ref{s:subcase 2}, and the remaining cases, discussed in Sections \ref{s:minimality non-edge-transitive}--\ref{s:non subcase 2}.  Many arguments become rather tedious for small values of odd $q$.  We thus assume that $q$ is large enough (in particular, our arguments work for $q\geq 300$ for the edge-transitive case, and $q\geq 514$ for the remaining cases).  Our discussion can be carried out in similar fashion when $p$ is odd and $q\leq 514$, but we decided to skip it in order to have ``cleaner" statements.  Hence our further assumption is that $q\geq 514$.

\subsection{Edge-transitive case}\label{s:minimality edge-transitive}

The main results are Lemmas~\ref{l:minimalcovcomcase} and~\ref{l:minimalcovnoncomcase} below, which show that, in the cases where $G$ admits an edge-transitive lattice, an edge-transitive lattice is a cocompact lattice of minimal covolume in $G$.  Thus by Theorem \ref{t:classification}, a cocompact lattice of minimal covolume is listed in the statement of Theorem \ref{t:cocompact existence}, and so the minimal covolume among cocompact lattices may be easily computed.

In the proofs below, we assume by contradiction that there is a cocompact lattice $\G$ of covolume strictly less than that of the edge-transitive lattice $\G_0$ of smallest covolume, and then show by careful case analysis that in fact $\mu(\G \bs G) \geq \mu(\G_0 \bs G)$.  Extending these arguments would, we believe, show that any cocompact lattice of minimal covolume in $G$ must be edge-transitive.  Since we did not need this further statement to determine the strict lower bound on covolumes given by Theorem \ref{t:covolumes}, and the proofs are already lengthy, we have not carried out this extension.

As before, for $i=1,2$, let  $P_i$ be a maximal parahoric subgroup of $G$, with $P_i$ the stabiliser in $G$ of a vertex $x_i$ of $X$.  Recall that $P_1\cong P_2$, and that if $L_i$ is a Levi factor of $P_i$, then $L_i=TM_i$ where $T\leq B\leq P_1\cap P_2$ is a torus of $G$, and $A_1(q)\cong M_i\triangleleft L_i$.  Now $M_i$ is normalised by $T$, and $T\cap M_i$ induces what are called inner-diagonal automorphisms on $M_i\cong A_1(q)$.  On the other hand there are various possibilities for the action of elements of $T-T\cap M_i$ on $M_i$. In particular, either ``none of them" or ``some of them" induce non-trivial outer-diagonal automorphisms on $M_i$.  This results in the following subcases:

\begin{subcase}[``none of them"]\label{sub1} For $i = 1,2$,  $L_i/Z(L_i)\cong PSL_2(q)$.
\end{subcase}
\begin{subcase}[``some of them"]\label{sub2} For $i = 1,2$,  $L_i/Z(L_i)\cong PGL_2(q)$.
\end{subcase}

%$$ \bf{Subcase\ 1}:\ \it{For} \ i=1,2,  L_i/Z(L_i)\cong PSL_2(q),\ \mbox{and}$$ $$ \bf{Subcase\ 2}:\ \it{For}  \ i=1,2,\  L_i/Z(L_i)\cong PGL_2(q).$$  

We further subdivide our discussion of the edge-transitive case based on this observation, with Subcase~\ref{sub1} being considered in Section \ref{s:subcase 1} and Subcase~\ref{sub2} in Section \ref{s:subcase 2}.

\subsection{Edge-transitive case, Subcase \ref{sub1}}\label{s:subcase 1}

In this case  $L_i=M_i\circ T_i$, that is, $L_i$ is a central (commuting) product of $M_i$ and  $T_i=C_T(M_i)$. It is possible but not necessary  that $T_i\cap M_i=1$.

\begin{example}  Let $G$ have generalised Cartan matrix $A=\left(\begin{matrix}2 & -2\\-2 & 2 \end{matrix}\right)$.
\begin{enumerate}
\item  Let  $p=2$ and $G=G_u$, the universal version of the group. Then $G$ is a central extension of $SL_2(\Fqt)$ by $\F_q^{\times}$, and so $L_i\cong C_{q-1}\times PSL_2(q)$ with $T_i\cap M_i=1$ and $|T_i|=q-1$.
\item Let  $p$ be an odd prime, and $G\cong SL_2(\Fqt)$. Then $L_i\cong SL_2(q)$ with $T_i=T\cap M_i=\langle -I\rangle\cong C_2$.
\end{enumerate}
\end{example}

An interesting and unusual consequence is the following observation that we already used in the proofs of Theorems \ref{t:cocompact existence} and \ref{t:classification}, but now will state more explicitly.

\begin{corollary}\label{c:centre}
Let $G$ be as in Theorem~\ref{t:cocompact existence}.  Suppose further that $q\equiv 3\mod 4$ and  $L_i/Z(L_i)\cong PSL_2(q)$ for $i=1,2$.  If $M_i\cong SL_2(q)$, then $Z(M_i)\leq Z(G)$, and in particular, $Z(G)\neq 1$.
\end {corollary}

\begin{proof}
We are now in the conclusions of Theorem~\ref{t:cocompact existence}\eqref{i: p = 2} and \eqref{i:PSL}.
If $M_i\cong SL_2(q)$, then because of the structure of $H_i$, we have $A_0\cap H_i=A_0\cap M_i= Z(M_i)$.  Hence $A_0\geq Z(M_i)\cong C_2$. But $A_0\leq Z(G)$, proving the result.
\end{proof}

Let us now discuss the question of covolumes.  Let $\Gamma = A_1 *_{A_0} A_2$ be an edge-transitive lattice in $G$.  Then by Theorem \ref{t:classification}, $\Gamma$ is one of the conclusions listed in Theorem~\ref{t:cocompact existence}\eqref{i: p = 2} or \eqref{i:PSL} above.  As recalled in Section \ref{s:lattices} above, the covolume of $\Gamma$ in $G$ may be calculated as follows: $$\mu(\G \bs G) = \frac{1}{|A_1|} + \frac{1}{|A_2|}=\frac{1}{(q+1)|A_0|} + \frac{1}{(q+1)|A_0|}=\frac{2}{(q+1)|A_0|}.$$  In all  conclusions to Theorem~\ref{t:cocompact existence}\eqref{i: p = 2} and \eqref{i:PSL}, the edge group $A_0$ satisfies $A_0\leq Z(G)$.  It follows that among all the \emph{edge-transitive} cocompact lattices in $G$, the lattice $\Gamma_0$ with edge group $A_0 = Z(G)$ has the \emph{smallest possible covolume}.

Now take $\Gamma$ to be a cocompact, not necessarily edge-transitive, lattice in $G$.  What happens then?

\begin{lemma}\label{l:minimalcovcomcase}  Let $G$ be as in Theorem~\ref{t:cocompact existence}\eqref{i: p = 2} or \eqref{i:PSL}.  In fact, if $p$ is odd, suppose that $q > 300$.  Assume that Conjecture~\ref{c:cocompact} holds.  Then an edge-transitive lattice is a cocompact lattice in $G$ of minimal covolume. \end{lemma}

\noindent Combined with the discussion above, Lemma~\ref{l:minimalcovcomcase} proves Theorem~\ref{t:covolumes} for this subcase.  Note that in the statement of Theorem~\ref{t:covolumes} in the introduction, $\delta = 1$ in this case.  

\begin{proof}[Proof of Lemma~\ref{l:minimalcovcomcase}] Since $Z(G)$ is finite, without loss of generality we may assume that $Z(G)=1$, hence all lattices act faithfully.  We have already constructed an edge-transitive lattice $\G_0$ of covolume $\mu(\G_0 \bs G)=\frac{2}{|Z(G)|(q+1)} = \frac{2}{q + 1}$.  We will show that this is a lattice of minimal covolume among the cocompact lattices of $G$.  In order to do so, assume that there exists a cocompact lattice $\G$ in $G$ whose covolume is strictly smaller than that of $\G_0$.  

The quotient graph $S := \G \bs X$ is a bipartite graph containing at least two adjacent vertices.  Let us call them $x_1$ and $x_2$.    Then $$\mu(\G\bs G)=\sum_{s\in S}\frac{1}{|\G_{s}|}\geq \frac{1}{|\G_{x_1}|}+\frac{1}{|\G_{x_2}|}.$$ 
 Moreover, as $\G_{x_i}$ is finite, by Proposition~\ref{p:finite_gp_tree} above, we may assume without loss of generality that $\G_{x_i}\leq P_i$.  Hence by abuse of notation, we will for $i = 1,2$ denote by $x_i$ the vertex of $X$ stabilised by $P_i$, so that the edge $[x_1,x_2]$ of $X$ is the edge stabilised by $B$.  By Conjecture \ref{c:cocompact}, we have $(|\G_{x_i}|,p)=1$.  Thus in fact we may suppose that $\G_{x_i}\leq L_i$.  Notice that as $T\leq P_1\cap P_2$, we have $\G\cap T\leq\G\cap P_i=\G_{x_i}$.  It follows that $\G_{x_1}\cap T=\G_{x_2}\cap T$ and since $(p,|\G_{x_i}|)=1$, $\G_{x_1}\cap\G_{x_2}\leq T$.

We now consider two further subcases depending on the value of $p$.  The case $p = 2$ is considered in Section \ref{s:minimality subcase 1 p = 2} and the case $p$ odd in Section \ref{s:minimality subcase 1 p odd}.

\subsubsection{Proof of Lemma \ref{l:minimalcovcomcase} when $p = 2$}\label{s:minimality subcase 1 p = 2}

In this case we have that $M_i\cong SL_2(q)\cong PSL_2(q)$.  Thus $L_i=M_i\times T_i$ with $T_i\leq T$ and $T_i$ isomorphic to a subgroup of $C_{q-1}$.  Since $q=2^a$ for some $a\in\N$, $(q-1)$ is odd. Moreover, as we already remarked, since we are assuming that Conjecture~\ref{c:cocompact} holds, $\G_{x_i}$ is a subgroup of $L_i$ of odd order.

We first observe that if $T_i \neq 1$ then $[T_i, M_j] \neq 1$, where $\{i,j\}=\{1,2\}$.  Otherwise, $T_i\leq C_G(\langle M_1, M_2\rangle)\leq Z(G)=1$, a contradiction.

Next we note that $\G_{x_i}\cap T_i \neq 1$ for at least one of $i=1,2$.  Suppose that $\G_{x_i}\cap T_i=1$.  Then $$\G_{x_i}\cong\G_{x_i}T_i/T_i\leq M_iT_i/T_i\cong M_i.$$  Since by Dickson's Theorem, the only subgroups of odd order of $M_i$ are isomorphic to subgroups of either $C_{q-1}$ or $C_{q+1}$, it follows that  $|\G_{x_i}|\leq (q+1)$. Therefore if $\G_{x_i}\cap T_i=1$ for both $i=1,2$, we would have that   $\mu(\G\bs G)\geq \frac{2}{(q+1)} = \mu(\G_0 \bs G)$, contradicting our assumption that the covolume of $\G$ is strictly smaller than the covolume of $\G_0$.

Thus without loss of generality we may assume that there exists $1\neq y_1\in\G_{x_1}\cap T_1$.  By the observations above, $o(y_1)\mid(q-1)$ and $[y_1,M_2]\neq 1$. In fact, $\langle y_1\rangle$ acts faithfully on $M_2$ inducing inner automorphisms.  Note also that as $y_1\in\G_{x_1}\cap T$, in fact, $y_1\in\G_{x_2}$.

If there exists $1\neq y\in\G_{x_2}$ with $o(y)\mid (q+1)$, then as $\langle y, y_1\rangle$ acts faithfully on $M_2\cong PSL_2(q)$, Dickson's Theorem implies that $\langle y, y_2\rangle$ contains an element of order $2$, a contradiction with the cocompactness of $\G$.  Hence, for every $g\in\G_{x_2}$, $o(g)\mid (q-1)$, and in fact $g\in T$.  Thus $\G_{x_2}\leq T$, and so $\G_{x_2}\leq\G_{x_1}$.  Now, if $\G_{x_2}\cap T_2\neq 1$, then going through the same discussion, but switching the roles of $x_1$ and $x_2$, we will obtain that $\G_{x_1}\leq \G_{x_2}$, implying that $\G_{x_1}=\G_{x_2}$. Thus if all the vertices $s$ in the quotient graph $S = \G \bs X$ have this property (i.e., $\G_s\cap T_s\neq 1$),  we will be getting the same group as a vertex group at every step, and so $\G$ will be finite, a contradiction.  Hence, we may assume that $\G_{x_2}\cap T_2=1$, and  it follows that $\G_{x_2}$ is a subgroup of $T$ isomorphic to a subgroup of $C_{q-1}$, with $[\G_{x_2},M_2]\neq 1$. In particular, $\G_{x_2}\leq\G_{x_1}$.

Let us look at the possibilities for $\G_{x_1}$. If $\G_{x_1}\leq T$, it follows that $\G_{x_1}\leq\G_{x_2}$ implying $\G_{x_1}=\G_{x_2}$.  In particular, $|\G_{x_i}|\leq (q-1)$ for $i=1,2$, which violates the minimality of covolume of $\G$.  Therefore,  $\G_{x_1}\not\leq T$. Now Dickson's Theorem implies that the only other choice is for $\G_{x_1}\leq H_1\times T_1$ where $H_1$ is a non-split torus of $M_1$.  If there are more than two vertices adjacent to $x_1$ in $S$, then $\mu(\G\bs G)\geq \frac{2}{q-1}>\mu(\G_0\bs G)$, a contradiction.  Hence on the quotient graph, $x_2$ is the only vertex neighbouring $x_1$. Moreover, all other vertices of $S$ are neighbours of $x_2$ and their stabilisers are of the same type as that  of $x_1$.  In particular, it follows that $\G_{x_2}$ is a normal subgroup of all other vertex groups, which contradicts the fact that $\G$ is a faithful lattice, since $Z(G)=1$.

\subsubsection{Proof of Lemma \ref{l:minimalcovcomcase} when $p$ is odd}\label{s:minimality subcase 1 p odd}

In this case, $q=p^a$ for some $a\in\N$, $q\equiv 3\pmod 4$ and we also assume $q > 300$.

Since we are working in the situation where $Z(G)=1$, Corollary~\ref{c:centre} implies that $Z(M_i)=1$.  Assume there exists an involution $t_i\in T_i$.  If $[t_i,M_j]\neq 1$, then  $t_i\in T$ would induce a non-trivial automorphism on $M_j$.  But $q\equiv 3\mod 4$,  and so all the non-trivial involutory inner automorphisms of $M_j$ come from the elements of a non-split torus, a contradiction.  Hence $t_i\in C_G(\langle M_i, M_j\rangle\leq Z(G)=1$, contradicting our assumption.  Therefore $L_i\cong M_i\times T_i$ with  $M_i\cong PSL_2(q)$ and $T_i$ isomorphic to a subgroup of $C_{\frac{q-1}{2}}$.  Moreover, as in the previous case, if  $T_i\neq 1$, it must act faithfully on $M_j$ inducing inner automorphisms, for otherwise $Z(G)\neq 1$.

If $\G_{x_i}\cap T_i=1$, then $$\G_{x_i}\cong\G_{x_i}T_i/T_i\leq M_iT_i/T_i\cong M_i.$$ In particular,  $|\G_{x_i}|\leq (q+1)$. Therefore if $\G_{x_i}\cap T_i=1$ for both $i=1,2$, it follows that $\mu(\G\bs G)\geq\mu(\G_0\bs G)$, a contradiction.  Hence, without loss of generality we may assume that there exists $1\neq y_1\in\G_{x_1}\cap T_1$ such that $\langle y_1\rangle=\G_{x_1}\cap T_1$.  Notice that $o(y_1)\mid\frac{q-1}{2}$ and $\langle y_1\rangle$ acts faithfully on $M_2$ inducing non-trivial inner automorphisms.  Recall that as  $y_1\in T\cap\G_{x_1}$, $y_1\in\G_{x_2}$. And so $\G_{x_2}$ acts on $M_2$ either as a subgroup of $N_{M_2}(M_2\cap T)$, or as a subgroup of $K_2$ where $K_2\in\{ S_4, A_5\}$. Notice that in the latter case $o(y_1)\in\{3,5\}$.

Assume first that $\G_{x_2}\cap T_2=1$. Then $\G_{x_2}$ is isomorphic to a subgroup of $M_2$.   If $\G_{x_2}$ is isomorphic to a subgroup of $K$, then using the previous paragraph we obtain that $\mu(\G\bs G)\geq\frac{1}{60}+\frac{1}{5\cdot 2\cdot(q+1)/2}>\frac{2}{q+1}$ for $q>47$. Since this obviously contradicts the minimality of covolume of $\G$, $\G_{x_2}$ must be isomorphic to a subgroup of $N_{M_2}(M_2\cap T)$.  Since we have full information about the action of $L_i$ on the set of edges $E_X(x_i)$ coming out of $x_i$, let us first consider the action of $\G_{x_1}$ on $E_X(x_1)$.  Since $T_1$ fixes $E_X(x_1)$ pointwise, we are interested in the action of the projection of $\G_{x_1}$ on $M_1$. Assume first that $\G_{x_1}$ acts transitively on $E_X(x_1)$.  Then in the quotient graph $x_1$ has a unique neighbour $x_2$.  Consider the set of  neighbouring vertices of $x_2$ in $S$.  If $x_1$ is its unique neighbour, $VS=\{x_1, x_2\}$ and since $\Gamma_{x_1}$ acts transitively on $E_X(x_1)$, the lattice $\Gamma$ is edge-transitive and we are done. Hence there are at least two neighbouring vertices of $x_2$ in $S$ (one already represented by $x_1$).  Let $z_1,\ldots, z_k$ be representatives of the other neighbouring vertices of $x_2$ in $S$. If each $\G_{z_i}$ acts transitively on the edges $E_X(z_i)$ coming out of $z_i$,  we are looking at the whole $VS=\{z_1, \ldots , z_k,  x_2, x_1\}$.  If  for some $i$, $C_{\Gamma_{z_i}}(M_{z_i})=1$, $\Gamma_{z_i}$ is isomorphic to a subgroup of $M_{z_i}$ and so $|\Gamma_{z_i}|\leq (q+1)$. Hence $\mu(\Gamma\bs G)\geq\mu(\Gamma_0\bs G)$, a contradiction.  Therefore for all $i$,  $C_{\Gamma_{z_i}}(M_{z_i})\neq 1$.

Now denote by $y_{z_i}$ an element of $\Gamma_{z_i}$ such that $\langle y_{z_i}\rangle=C_{\Gamma_{z_i}}(M_{z_i})$. Then  just like $y_1$, $y_{z_i}$ acts faithfully on $M_2$ and $y_{z_i}\in\Gamma_{x_2}$. But $\Gamma_{x_2}$ is a subgroup of the normaliser of the split torus of $M_2$, and so $[y_1, y_{z_i}]=1$. Again using the fact that we know how $L_2$ acts on $E_X(x_2)$,  we observe that $y_1$ then must fix the edge $(x_2, z_i)$.  It follows that $y_1$ fixes all vertices of $S$ and so is in the kernel of the action of $\Gamma$, which since $\G$ acts faithfully is a contradiction.  Therefore, either $\G_{x_1}$ or some $\G_{z_i}$ does not act transitively on the corresponding set of edges, and so by taking $x_1=z_i$ if necessary, we are in the following case: $\G_{x_1}$ has at least two orbits on $E_X(x_1)$.

Therefore, either $x_1$ has a unique neighbouring vertex $x_2$ in $S$ but the number of edges between them in $S$ is greater than one, or there are at least two neighbouring vertices of $x_1$ in the quotient graph $S$.  Let us discuss these two cases. In the former one, if $VS=\{x_1, x_2\}$, then $\langle y_1\rangle\triangleleft\G$, which is a contradiction as $\G$ is a faithful lattice since $Z(G) = 1$.  It follows that $|VS|>2$ and $x_2$ has at least two neighbouring vertices $z_1, \ldots, z_k$.  Not to contradict the minimality of $\Gamma$ we may now assume that $C_{\Gamma_{z_i}}(M_{z_i})\neq 1$.  Using the same argument as above, we obtain that $y_1\in\Gamma_{z_i}$.  Assume that $y_1$ acts on $M_{z_i}$ as an element of $K_{z_i}\in\{S_4, A_5\}$. Suppose first that $o(y_1)\leq 15$.  Let us evaluate the covolume of $\Gamma$. Let $a_2\in\Gamma_{x_2}$ be such that $\langle a_2\rangle\leq\Gamma_{x_2}=\langle a_2\rangle\rtimes\langle s_2\rangle\cong D_{2\cdot o(a_2)}$ and $\langle a_2\rangle=T\cap\Gamma_2$.  Then $a_2\in \Gamma_{x_1}$ and $y_1\in\langle a_2\rangle$. If $o(y_1)=o(a_2)$, then $|\Gamma_{x_2}|\leq 30$ which immediately contradicts the minimality of covolume of $\Gamma$ ($\frac{1}{30}\geq \frac{2}{q+1}$ for $q>59$).  Hence, $o(a_2)>o(y_1)$ and so $\Gamma_{x_1}$ acts on $M_1$ as either a subgroup of $K_1\in\{S_4, A_5\}$, or as a subgroup of $N_{M_1}(M_1\cap T)$.  

In the former case $|\Gamma_{x_1}|\leq 15\cdot 60$ and in particular, $|\G_{x_1}\cap T|\leq 15\cdot 5=75$.  As we noticed earlier,  $\G_{x_1}\cap T=\G_{x_2}\cap T$, and so $|\G_{x_2}\cap T|\leq 75$. Since $\G_{x_2}$ is isomorphic to a subgroup of a normaliser of split torus of $M_2$, $|\G_{x_2}|\leq 2\cdot 75=150$. If follows that $\mu(\G\bs G)>\frac{2}{q+1}$ for $q\geq 257$, again a clear contradiction. In the latter case $|\Gamma_{x_1}|\leq o(y_1)2\frac{o(a_2)}{o(y_1)}\leq (q-1)$. Since $|\Gamma_{x_2}|\leq (q-1)$, we again get a contradiction with the minimality of covolume of $\Gamma$. Hence $o(y_1)>15$ and so there exists $1\neq y_1'\in\langle y_1\rangle$ such that $[y_1',M_{z_i}]=1$ for all $i$.  We then use $y_1'$ in the place of $y_1$ in all the previous subgroups.  Therefore $y_1$ either centralises $M_{z_i}$ or acts on it as a subgroup of a normaliser of a split torus of $M_{z_i}$. In both situations, $\langle y_1\rangle$ is normal in $\Gamma_{z_i}$ for $i=1,\ldots , k$.  Continuing with this argument, we obtain that $y_1$ is in the kernel of the action of $\Gamma$, a contradiction.  

Therefore we are in the  case that $x_1$ has more than one neighbouring vertex in $S$.  One is $x_2$ and let $z$ be among the other neighbouring vertices of $x_1$.  If $\G_z\cap T_z=1$, then using $|\Gamma_{x_2}|$ and $|\Gamma_z|$, we obtain a contradiction with the minimality of covolume of $\Gamma$.  Therefore $\G_z\cap T_z\neq 1$ and we may take $x_2=z$.

Thus we may assume that  $\G_{x_2}\cap T_2\neq 1$. Hence there exists $y_2\in \G_{x_2}\cap T_2$ with $\langle y_2\rangle=\G_{x_2}\cap T_2$.  As before  we notice that $o(y_2)\mid\frac{q-1}{2}$,  $y_2\in \G_{x_1}$ and $\langle y_2\rangle$ acts faithfully on $M_1$ via inner automorphisms.  As before now Dickson's Theorem  allows us to conclude that either $\G_{x_1}$ acts on $M_1$ as a subgroup of $K_1$ with $K_1\in\{S_4, A_5\}$, or $\G_{x_1}$ acts on $M_1$ as  a subgroup of $N_{M_1}(M_1\cap T)$.  In the former case $o(y_2)\in\{3, 5\}$.  

Let us begin with the case when $\G_{x_1}$ acts on $M_1$ as a subgroup of $K_1$.  If $\G_{x_2}$ acts on $M_2$ as a subgroup of $K_2$, then $\mu(\G \bs G)\geq\frac{2}{60\cdot 5}>\frac{2}{q+1}=\mu(\G_0\bs G)$ for $q\geq 300$, a contradiction.   Hence $\G_{x_2}$ acts on $M_2$ as a subgroup $N_{M_2}(M_2\cap T)$.  Recall that $|VS|>2$.  Let  $v_1, \ldots , v_k$ be the neighbours of $x_1$ in $VS-\{x_2\}$. Since $y_1$ fixes every edge in $E_X(x_1)$, it follows that $y_1$ acts faithfully on $M_{v_i}$ and holding a discussion similar to the above one with $v_i$ in place of $x_2$, we may assume that $\G_{v_i}$ acts on $M_{v_i}$ as a subgroup of a normaliser of a split torus of $M_{v_i}$. It follows that $\langle y_1\rangle $ is normal in each $\G_{v_i}$.  Now let $z_1, \ldots , z_m$ be the neighbours of $x_2$ in $S-\{x_1\}$.  Let us consider $\G_{z_i}=P_{z_i}\cap\G$.  If $\G_{z_i}\cap C_G(M_{z_i})=1$,  then there is at most one such vertex, otherwise we would contradict the minimality of covolume of $\Gamma$.  Hence we may assume that if it happens, $i=1$, that is,  $\G_{z_1}\cap C_G(M_{z_1})=1$. Then we may further assume that $T\leq P_{z_1}$.  Thus $y_1, y_2\in\Gamma_{z_1}$. If $\Gamma_{z_1}$ acts on $M_{z_1}$ as a subgroup of $K_{z_1}\in\{ S_4, A_5\}$, then $|\Gamma_{z_1}|\leq 60$, which is a contradiction, as always ($\frac{1}{60}\geq\frac{2}{q+1}$ for $q>120$).  Hence $\Gamma_{z_1}$ acts on $M_{z_1}$ as a subgroup of a normaliser of a split torus of $M_{z_1}$ (it is not possible to have a normaliser of a non-split torus, as $o(y_j)$, $j=1,2$, is co-prime to $(q+1)$).  It follows that $\langle y_1\rangle$ is a normal subgroup of $\Gamma_{z_1}$.  Hence, for $i>1$, there exists $y_{z_i}\in  C_G(M_{z_i})$ of order dividing $\frac{q-1}{2}$. But this element sits in the kernel of action of $L_{z_i}$ on $E_X(z_i)$ and therefore, $y_{z_i}\in \Gamma_{x_2}$.  On the other hand by the usual argument, $y_2$ acts faithfully on $M_{z_i}$ and so $[y_2, y_{z_i}]=1$. It follows that $\langle y_{z_i}\rangle$  is normal in $\G_{x_2}$.

Finally, as $C_{\Gamma_{x_2}}(M_2)$ stabilises $(x_2,z_i)$, it follows that $C_{\G_{x_2}}(y_{z_i})\leq\G_{z_i}$.  It follows that $y_1\in\G_{z_i}$ and so $\langle y_1, y_2\rangle\leq \G_{z_i}$.  Assume that $y_1$ acts on $M_{z_i}$ as a subgroup of $K_{z_i}\in\{S_4, A_5\}$.  If $o(y_1)\leq 15$, then $|\G_{x_1}|\leq 15\cdot 60$ and in particular, $|\G_{x_1}\cap T|\leq 15\cdot 5=75$.  As we remarked earlier, $\G_{x_1}\cap T=\G_{x_2}\cap T$, and so $|\G_{x_2}\cap T|\leq 75$. Since $\G_{x_2}$ acts on $M_2$ as a subgroup of a normaliser of split torus, $|\G_{x_2}|\leq 2\cdot 75=150$, which leads to the usual contradiction as $\frac{1}{15\cdot 60}+\frac{1}{150}\geq\frac{2}{q+1}$ for $q\geq 257$.  Hence, $o(y_1)>15$ and so there exists $y_1'\in\langle y_1\rangle$ with $[y_1',M_{z_i}]=1$ for all $i$.  In this case we  will use $y_1'$ in the place of $y_1$ in all the previous subgroups.  By iterating this argument we may show that $\langle y_1\rangle\triangleleft\G$ which is a contradiction.

 We are now reduced to the last possible situation: $\G_{x_1}$ acts on $M_1$ as a subgroup of $N_{M_1}(M_1\cap T)$. Notice that because of the symmetry between $x_1$ and $x_2$, to finish the analysis it remains to consider the case when $\G_{x_2}$ acts on $M_2$ as a subgroup of $N_{M_2}(M_2\cap T)$.

But in this case $\langle \G_{x_1}, \G_{x_2}\rangle \leq N$. Hence we may move to the next vertex $y$ on our graph. Using the previous argument we obtain again that the only possibility will be $\G_y\leq N$, and so on and so forth. Therefore, in the end of this case, the only possible conclusion will be $\G\leq N$.  This is a contradiction, as $N$ is not a cocompact lattice of $G$, nor does it contain any cocompact lattice.
\end{proof}

\subsection{Edge-transitive case, Subcase \ref{sub2}}\label{s:subcase 2}

We are now in the case that $G$ admits edge-transitive lattices and $T$ induces non-trivial outer-diagonal automorphisms on $M_i$ for  $i=1,2$. As $M_i\cong A_1(q)$,  $p$ is $\bf{odd}$, for if $p=2$,  $A_1(q)=SL_2(q)=PSL_2(q)$ does not admit non-trivial outer-diagonal automorphisms.  It follows that $L_i$ is  isomorphic to a homomorphic image of $GL_2(q)$. In particular,  $L_i=T_iM_i\langle t_i\rangle$ where $T_i=C_T(M_i)$, $t_i\in T-T_i$ and $L_i/T_i\cong PGL_2(q)$.  Moreover, if $q\equiv 3\pmod 4$, $T_i/T_i\cap M_i$ is a cyclic group of odd order and $t_i\in T$ is an involution, while if $q\equiv 1\pmod 4$, $1\neq t^2_i\in M_iT_i$.

As in Subcase \ref{sub1} in Section \ref{s:subcase 1} above, we investigate minimality of covolumes.  Assume that $q\geq 300$ and let $\Gamma = A_1 *_{A_0} A_2$ be an edge-transitive lattice in $G$.  Reading carefully through the statement of Theorem~\ref{t:cocompact existence}\eqref{i:PGL}, we observe that $|A_0|=\delta|Z_0|$ where $Z_0\leq Z(G)$ and $\delta\in\{1, 2, 4\}$ depends on the structure of $G$ and $\Gamma$.  

%Notice that if $q\equiv 1\pmod 4$, then if $Q_i\not\leq Z(G)$, $\delta=4$, while if $Q_i\leq Z(G)$, $\delta=2$ in case (A) and $\delta_0=1$ in case (B).  On the other hand, if $q\equiv 3\pmod 4$, then if $Z(M_i)\not\leq Z(G)$, $\delta_0=4$, while if $Z(M_i)\leq Z(G)$, $\delta_0=2$ in case (i) and $\delta_0=1$ in case (iii).

Now, among all the edge-transitive cocompact lattices in $G$, choose $\Gamma_0=A_1 *_{A_0} A_2$ such that $|A_0|$ is as large as possible.  It follows that if $q\equiv 1\pmod 4$, $\Gamma_0$ is described in \eqref{i:PGL q = 1} of Theorem~\ref{t:cocompact existence}, while if $q\equiv 3\pmod 4$, $\Gamma_0$ is described in \eqref{i:PGL q = 3} of Theorem~\ref{t:cocompact existence}, and in both cases $Z_0=Z(G)$. In particular, $|A_0|=\delta|Z(G)|$ with $\delta\in\{2, 4\}$.  Therefore for any other edge-transitive lattice $\G = A_1 *_{A_0} A_2$ in $G$,  we have  $$\mu(\G \bs G)\geq \mu(\G_0 \bs G)=\frac{2}{(q+1)|Z(G)|\delta}$$ where $\delta\in\{ 2, 4\}$ as described above.  And so among all the \emph{edge-transitive} cocompact lattices in $G$, the lattice $\Gamma_0$ with edge group $A_0$ of order $|Z(G)|\delta$ has the \emph{smallest possible covolume}.

Now take $\Gamma$ to be a cocompact, not necessarily edge-transitive, lattice in $G$.  What happens then?

\begin{lemma}
\label{l:minimalcovnoncomcase}
Let $G$ be as in Theorem~\ref{t:cocompact existence}\eqref{i:PGL}, $q\geq 300$, and assume that $G$ admits an edge-transitive lattice. Assume further that Conjecture~\ref{c:cocompact} holds.  Then an edge-transitive lattice is a cocompact lattice in $G$ of minimal covolume.
\end{lemma}

\noindent Again, the discussion above together with Lemma~\ref{l:minimalcovnoncomcase} proves Theorem~\ref{t:covolumes} for this case.

\begin{proof}[Proof of Lemma~\ref{l:minimalcovnoncomcase}]  We consider the case $q \equiv 3 \pmod 4$ in Section \ref{s:subcase 2 q = 3} and then the case $q \equiv 1 \pmod 4$ in Section \ref{s:subcase 2 q = 1}. 

\subsubsection{Proof of Lemma \ref{l:minimalcovnoncomcase} when $q \equiv 3$ (mod 4)}\label{s:subcase 2 q = 3}

Consider $L_i=M_iT$.  Then $M_i\triangleleft L_i$, $T_i=C_T(M_i)$ is a homomorphic image of $C_{q-1}$, $T_i/Z(M_i)$ is a cyclic group of order dividing $\frac{q-1}{2}$ and there exists an involution $t_i\in T$ that induces an outer-diagonal ($PGL_2$-) automorphism on $M_i$.  Notice that even if $Z(G)=1$, this time it is possible to have $Z(M_i)\neq 1$.  However,  if $Z(G)=1$, $T_i$ acts faithfully on $M_j$ where $\{i,j\}=\{1,2\}$ for $C_{T_i}(M_j)\leq C_T(\langle M_1, M_2\rangle)\leq Z(G)$.  Again we assume without loss of generality that $Z(G)=1$, and we suppose that $\G$ is a cocompact lattice in $G$ whose covolume is strictly less than that of $\G_0$ (discussed above).  

A few more comments: if $Z(M_i)=1$,  $|T|_2=2$ and so there exists a unique involution $t\in T$ that induces a $PGL_2$--automorphism on both $M_1$ and $M_2$. As above, if $Z(M_i)=1$, $\delta=2$, while if $|Z(M_i)|=2$, $\delta=4$.  The covolume of $\G_0$ is then $\frac{1}{(q+1)\delta}$. Put $\delta_0:=\frac{\delta}{2}$.

If $\G_{x_i}\cap T_i\leq Z(M_i)$, then $$\G_{x_i}Z(M_i)/Z(M_i)\cong\G_{x_i}T_i/T_i\leq M_iT/T_i\cong PGL_2(q).$$  In particular,  $|\G_{x_i}|\leq (q+1)\delta$. Therefore if $\G_{x_i}\cap T_i\leq Z(M_i)$ for both $i=1,2$, it follows that $\mu(\G\bs G)\geq\mu(\G_0 \bs G)$, a contradiction.  Hence  without loss of generality we may assume that there exists $1\neq y_1\in\G_{x_1}\cap T_1-Z(M_1)$ such that $\langle y_1\rangle=O(\G_{x_1}\cap T_1)$.  (Recall that for a finite group $H$, $O(H)$ denotes the largest normal subgroup of $H$ of odd order.)  Notice that $o(y_1)\mid\frac{q-1}{2}$ and $\langle y_1\rangle$ acts faithfully  on $M_2$ via  inner automorphisms.  Recall that as $y_1\in T\cap\G_{x_1}$, $y_1\in\G_{x_2}$. And so $\G_{x_2}$ acts on $M_2$ either as a subgroup of a normaliser of a split torus of $M_2\langle t_2\rangle$, or as a subgroup of $K_2$ where $K_2/Z(K_2)\in\{ S_4, A_5\}$. Note that in the latter case $o(y_1)\in\{3,5\}$.

Assume first that $\G_{x_2}\cap T_2\leq Z(M_2)$. Then $\G_{x_2}$ is isomorphic to a subgroup of $M_2\langle t_2\rangle$.  If $\G_{x_2}$ is isomorphic to a subgroup of $K_2$, then using the previous paragraph we obtain that $\mu(\G\bs G)\geq\frac{1}{60\delta_0}+\frac{1}{5\cdot 2\cdot(q+1)\delta_0}>\frac{2}{2(q+1)\delta_0}=\frac{2}{(q+1)\delta}$ for $q>53$. Since this obviously contradicts the minimality of covolume of $\G$,  $\G_{x_2}$ must be isomorphic to a subgroup of a normaliser of a split torus in $M_2\langle t_2\rangle$.

Since we have full information about the action of $L_i$ on the set of edges $E_X(x_i)$ coming out of $x_i$, $i=1,2$, let us first consider the action of $\G_{x_1}$ on $E_X(x_1)$.  Since $T_1$ fixes $E_X(x_1)$ pointwise, we are interested in the action of the projection of $\G_{x_1}$ on $M_1\langle t_1\rangle$.  Assume first that $\G_{x_1}$ acts transitively on $E_X(x_1)$.  Then in the quotient graph $S$,  $x_1$ has a unique neighbour vertex $x_2$.  Consider the set of  neighbours of $x_2$ in $S$.  If $VS=\{x_1, x_2\}$, the quotient graph corresponds to the edge-transitive lattice, and we are done.  Therefore $x_2$ has at least two neighbours vertices in $S$.  Let $z_1, \ldots , z_k$ be representatives of the neighbouring vertices of $x_2$ in $S$ other than $x_1$.  If each $\G_{z_i}$ acts transitively on the set of edges $E_X(z_i)$ coming out of $z_i$, then we are looking at the whole $VS=\{z_1, \ldots , z_k, x_2, x_1\}$.  If for some $i$, $C_{\Gamma_{z_i}}(M_{z_i})\leq Z(M_i)$, $\Gamma_{z_i}$ is isomorphic to a subgroup of $PGL_2(q)$ and so $|\Gamma_{z_i}|\leq(q+1)\delta$. Hence $\mu(\Gamma\bs G)\geq\mu(\Gamma_0\bs G)$, a contradiction.  Therefore for all $i$, $C_{\Gamma_{z_i}}(M_{z_i})\not\leq Z(M_i)$.  Denote by $y_{z_i}$ an element of $\Gamma_{z_i}$ such that $\langle y_{z_i}\rangle=O(C_{\Gamma_{z_i}}(M_{z_i}))$.  Then just like $y_1$, $y_{z_i}$ acts faithfully on $M_2$ and $y_{z_i}\in\Gamma_{x_2}$.  But $\Gamma_{x_2}$ is a subgroup of the normaliser of the split torus of $M_2$, and so $[y_1, y_{z_i}]=1$.  Again using the fact that we know how $L_2$ acts on the set $E_X(x_2)$, we observe that $y_1$ must fix the edge $(x_2, z_i)$.  It follows that $y_1$ fixes all vertices of $S$, and so is in the kernel of the action of $\Gamma$, a contradiction. Therefore, either $\Gamma_{x_1}$, or some $\Gamma_{z_i}$ does not act transitively on the corresponding set of edges. Thus by taking $x_1=z_i$ if necessary, we reduce to the following case: $\Gamma_{x_1}$ has at least two orbits on $E_X(x_1)$.

Therefore, either $x_1$ has a unique neighbouring vertex $x_2$ in $S$ but the number of edges between $x_1$ and $x_2$ in $S$ is greater than one, or there are at least two neighbouring vertices of $x_1$ in the quotient graph.  Let us discuss these two cases.  In the former one, if $|VS|=2$, $\langle y_1\rangle\triangleleft\Gamma$ which is a contradiction since $\Gamma$ is a faithful lattice. It follows that in the former case $|VS|>2$  and so $x_2$ has more than one neighbouring vertex in $S$: $z_1, \ldots , z_k$.   Again not to contradict the minimality of covolume of $\Gamma$, we may assume that $C_{\Gamma_{z_i}}(M_{z_i})\not\leq Z(M_i)$. Using the same argument as above we obtain that $y_1\in\Gamma_{z_i}$.  Assume that $y_1$ acts on $ M_{z_i}$ as an element of $K_{z_i}$ where $K_{z_i}/Z(K_{z_i})\in\{S_4, A_5\}$.

Suppose first that $o(y_1)\leq 15$. Let us evaluate the covolume of $\Gamma$. Let $a_2\in\Gamma_{x_2}$ be such that $\langle a_2\rangle=O(\Gamma_{x_2})$. Then $y_1\in\langle a_2\rangle$. If $o(y_1)=o(a_2)$, $|\Gamma_{x_2}|\leq 30\delta$ which immediately contradicts the minimality of covolume of $\Gamma$ ($\frac{1}{30\delta}>\frac{2}{(q+1)\delta}$ as long as $q>59$).  Hence, $o(a_2)>o(y_1)$.  Now, consider an action of $\Gamma_{x_1}$ on $M_1$.  If it acts as a subgroup of $K_1$ with $K_1/Z(K_1)\in\{S_4, A_5\}$, then $|\Gamma_{x_1}|\leq 15\cdot 60\delta_0$ and in particular, $|\G_{x_1}\cap T|\leq 15\cdot 5\delta=75\delta$.  But $\G_{x_1}\cap T=\G_{x_2}\cap T$, and so $|\G_{x_2}\cap T|\leq 75\delta$. As $\G_{x_2}$ acts on $M_2$ as a subgroup of a normaliser of split torus, $|\G_{x_2}|\leq 2\cdot 75\delta=150\delta$.  Therefore, $\mu(\G\bs G)\geq\frac{1}{900\delta_0}+\frac{1}{150\delta}\geq\frac{2}{(q+1)\delta}$ for $q\geq 225$, a contradiction with the minimality of covolume of $\Gamma$.  If it acts as a subgroup of a normaliser of a non-split torus, then together with the action $a_2$, we will get that $\Gamma_{x_1}$ will act as a subgroup containing $p$--elements, a contradiction. Thus it can only act on $M_1$ as a subgroup of a normaliser of a split torus, and so $|\Gamma_{x_1}|\leq o(y_1)2\frac{o(a_2)}{o(y_1)}\delta\leq (q-1)\delta$. Since $|\Gamma_{x_2}|\leq (q-1)\delta$, we again get a contradiction with the minimality of the covolume of $\Gamma$.  Hence, $o(y_1)>15$ and so there exists $y_1'\in\langle y_1\rangle$ such that $[y_1',M_{z_i}]=1$ for all $i$.  In this case we will use $y_1'$ instead of $y_1$ in all the previous discussions.  Therefore for all $i$, either $y_1$ centralises $M_{z_i}$ or acts on it as a normal subgroup of a normaliser of a split torus.  In both situations, $\langle y_1\rangle$ is normal in $\Gamma_{z_i}$, for $i=1, \ldots , k$.  Continuing with this argument we obtain that $\langle y_1\rangle$ is a normal subgroup of $\Gamma$ which implies that $y_1$ is in the kernel of the action of $\Gamma$ on $X$, a contradiction.

Therefore we may suppose that $x_1$ has more than one neighbouring vertex in $S$. One is $x_2$, and let $z$ be among the other neighbouring vertices of $x_1$.  If $C_{\Gamma_{z}}(M_{z})\leq Z(M_z)$, then using $|\Gamma_{x_2}|$ and $|\Gamma_z|$, we obtain a contradiction with the minimality of $\mu(\Gamma\bs G)$. Therefore $C_{\Gamma_z}(M_z)\not\leq Z(M_z)$ and we may take $x_2=z$. Therefore, we may assume that $\Gamma_{x_2}\cap T_2\not\leq Z(M_2)$.  It follows that there exists $y_2\in\Gamma_2\cap T_2$ with $\langle y_2\rangle=O(\Gamma_2\cap T_2)$. As before notice that $y_2\in\Gamma_{x_1}$ and $\langle y_2\rangle $ acts faithfully on $M_1$ via inner automorphisms.  Also, the subgroup structure of $L_1$ implies that either $\Gamma_{x_1}$ acts on $M_1$ as a subgroup of $K_1$ where $K_1/Z(K_1)\in\{S_4, A_5\}$, or as a subgroup of a normaliser of a split torus of $M_1\langle t_1\rangle$.  Assume  the former, i.e.,  let $\Gamma_{x_1}$ act on $M_1$ as a subgroup of $K_1$.  If $\Gamma_{x_2}$ acts on $M_2$ as a subgroup of $K_2$ where $K_2/Z(K_2)\in\{S_4, A_5\}$, then $\mu(\Gamma\bs G)\geq\frac{2}{5\cdot 60\delta_0}>\frac{2}{(q+1)\delta}=\mu(\Gamma_0\bs G)$ for $q\geq 300$, a contradiction.  Hence, $\Gamma_{x_2}$ acts on $M_2$ as a subgroup of a normaliser of a split torus. Recall that $|VS|>2$ and $x_1$ has neighbouring vertices in $S$ other than $x_2$. Let us call them $v_1, \ldots , v_k$. Arguing for $v_i$ as we did for $x_2$, we obtain that $\langle y_1\rangle$ acts faithfully on $M_{v_i}$ and that $\Gamma_{v_i}$ acts on $M_{v_i}$ as a normaliser of a non-split torus.  It follows that $\langle y_1\rangle$ is normal in each $\Gamma_{v_i}$.

Now let $z_1, \ldots , z_m$ be the neighbouring vertices of $x_2$ in $S-\{x_1\}$. Let us consider $\Gamma_{z_i}=\Gamma\cap P_{z_i}$.  If $C_{\Gamma_{z_i}}(M_{z_i})\leq Z(M_{z_i})$ for some $i$, then there is at most one such vertex, for otherwise we would contradict the minimality of covolume of $\Gamma$.  Assume there exists such a vertex $z_i$, then without loss of generality assume that $i=1$, i.e.,  $C_{\Gamma_{z_1}}(M_{z_1})\leq Z(M_{z_1})$ and for $i>1$, $C_{\Gamma_{z_i}}(M_{z_i})\not\leq Z(M_{z_i})$. Further, we may assume that $T\leq P_{z_1}$.  Then $y_1, y_2\in\Gamma_{z_1}$.  If $\Gamma_{z_1}$ acts on $M_{z_1}$ as a subgroup of $K_{z_1}$ where $K_{z_1}/Z(K_{z_1})\in\{S_4, A_5 \}$, then $|\Gamma_{z_1}|\leq 60\delta_0$, which leads to the usual contradiction. Thus $\Gamma_{z_1}$ acts on $M_{z_1}$ as a subgroup of a normaliser of a split torus, and so $\langle y_1\rangle\triangleleft \Gamma_{z_1}$.  Now, let us look at $\Gamma_{z_i}$ for $i\geq 2$ (or it is possible that the previous case does not happen, then we are looking at $\Gamma_{z_i}$ for $i\geq 1$).  It follows  that for each $i$ there exists $y_{z_i}$ such that $\langle y_{z_i}\rangle=O(C_{\Gamma_{z_i}}(M_{z_i}))$.  Since $y_{z_i}$ sits in the kernel of $L_{z_i}$ in its action on $E_X(z_i)$, $y_{z_i}$ fixes every vertex neighbouring $z_i$, and in particular fixes $x_2$. It follows that $y_{z_i}\in\Gamma_{x_2}$. On the other hand, by the usual argument, $y_2$ acts faithfully on $M_{z_i}$ implying $[y_2, y_{z_i}]=1$. Therefore, $\langle y_{z_i}\rangle\triangleleft\Gamma_{x_2}$.  Finally, as  $C_{\Gamma_{x_2}}(M_{2})$ stabilises every $(x_2, z_i)$, it follows that  $C_{\Gamma_{x_2}}(y_{z_i})\leq\Gamma_{z_i}$. It follows that $y_1\in\Gamma_{z_i}$ and so $\langle y_1, y_2\rangle\leq\Gamma_{z_i}$.  Suppose  $y_1$ acts on $M_{z_i}$ as an element of a subgroup of $K_{z_i}$ with $K_{z_i}/Z(K_{z_i})\in\{S_4, A_5\}$.  If $o(y_1)\leq 15$, $|\G_{x_1}|\leq 15\cdot 60\delta_0$ and in particular, $|\G_{x_1}\cap T|\leq 15\cdot 5\delta=75\delta$.  Now, as we noticed earlier, $\G_{x_1}\cap T=\G_{x_2}\cap T$, and so $|\G_{x_2}\cap T|\leq 75\delta$. As $ \G_{x_2}$ acts on $ M_2$ as a subgroup of a normaliser of split torus, $|\G_{x_2}|\leq 2\cdot 75\delta=150\delta$, leading to the usual contradiction for $q\geq 300$.   Hence, $o(y_1)>15$ and so there exists $y_1'\in\langle y_1\rangle$ such that $[y_1',M_{z_i}]=1$ for all $i$s.  In this case we will use $y_1'$ instead of $y_1$ in all the previous discussions.  By carefully  iterating this argument, we obtain that there exists a nontrivial element $y$ of odd order such that $\langle y\rangle\triangleleft\Gamma$, a contradiction.

We are now reduced to the case when $\G_{x_1}$ acts on $M_1\langle t_1\rangle$ as a subgroup of $N_{M_1\langle t_1\rangle}(T)$. Notice, that because of symmetry between $x_1$ and $x_2$, to finish the analysis it remains to consider the case when $\G_{x_2}$ acts on $M_2\langle t_2\rangle$ as a subgroup of $N_{M_2\langle t_2\rangle}(T)$.  But in this case $\langle \G_{x_1}, \G_{x_2}\rangle \leq N$. Hence, we may move to the next vertex $y$ on our graph. However, again, the only possible case will be $\G_y\leq N$, and so on and so forth. Therefore, in the end of this case, the only possible conclusion will be $\G\leq N$, which is a contradiction as  $N$ is not a cocompact lattice of $G$, nor does it contain one.  This completes the proof of Lemma \ref{l:minimalcovnoncomcase} when $q \equiv 1 \pmod 4$.

\subsubsection{Proof of Lemma \ref{l:minimalcovnoncomcase} when $q \equiv 1$ (mod 4)}\label{s:subcase 2 q = 1}

Again consider $L_i=M_iT$, $i=1,2$. As in the case $q \equiv 3 \pmod 4$, $M_i\triangleleft L_i$, $T_i$ is a homomorphic image of $C_{q-1}$ and $T_i/Z(M_i)$ is a cyclic group of order dividing $\frac{q-1}{2}$. But this time if $x$ is an involution in $L_i\cap T$, then $x\in M_iT_i$.   As usual we may suppose that $Z(G)=1$. As before, it follows  that $T_i$ must act faithfully on $M_j$, $\{i, j\}=\{1, 2\}$.

Recall that $Q_i\in\mathcal{S}yl_2(Z(L_i))$ and $Q_i^2$ is the unique subgroup of $Q_i$ of index $2$.  Since we assume that $G$ admits an edge-transitive lattice as in the conclusion \eqref{i:PGL q = 1} of Theorem~\ref{t:cocompact existence} and that $Z(G)=1$, $T_i=Z_i\times Q_i$ where $|Q_i|\leq 2$, $|Z_i|$ is odd and $M_i\cong PSL_2(q)$. If $|Q_i|=2$,  $\delta=4$, while if $|Q_i|=1$, $\delta=2$.  The covolume of $\Gamma_0$ is then $\frac{2}{(q+1)\delta}$. 
Notice that if $|Q_i|=1$, the lattice in the conclusion \eqref{i:PGL q = 1 other} of Theorem~\ref{t:cocompact existence} also exists, but its covolume is twice $\mu(\G_0\bs G)$.
 Put $\delta_0:=\frac{\delta}{2}$.

As in the previous case, let us assume that $\Gamma$ is a cocompact lattice in $G$ whose covolume is strictly less than that of $\Gamma_0$.  If $\Gamma_{x_i}\cap T_i\leq Q_i$, then $$\G_{x_i}Q_i/Q_i \cong\G_{x_i}T_i/T_i\leq M_iT/T_i\cong PGL_2(q).$$  In particular,  $|\G_{x_i}|\leq (q+1)\delta_0$. Therefore if $\G_{x_i}\cap T_i\leq Q_i$ for both $i=1,2$, it follows that $\mu(\G\bs G)\geq\mu(\G_0 \bs G)$, a contradiction.  Hence  without loss of generality we may assume that there exists $1\neq y_1\in\G_{x_1}\cap T_1$ such that $\langle y_1\rangle=O(\G_{x_1}\cap T_1)$. Note that $1\neq o(y_1)$ is odd, $o(y_1)\mid\frac{q-1}{2}\delta_0$ and $\langle y_1\rangle$ acts faithfully  on $M_2$.  Recall that as $y_1\in T\cap\G_{x_1}$, $y_1\in\G_{x_2}$. And so there are several possibilities  for the action of $\G_{x_2}$ on $M_2$: either as a subgroup of  a normaliser of a split torus, or as a subgroup of $K_2$ where $K_2\in\{ S_4, A_5\}$, in which case $o(y_1)\leq 5$.

Assume first that $\G_{x_2}\cap T_2\leq Q_2$. Then $\G_{x_2}$ is isomorphic to a subgroup of $2PGL_2(q)$.  If $\G_{x_2}$ is isomorphic to a subgroup of $ K_2$, we obtain the usual contradiction with the minimality of covolume, as $\frac{1}{60}>\frac{2}{(q+1)\delta}$ for $q\geq 60$.  Hence  $\G_{x_2}$ acts on $M_2$ a subgroup of a normaliser of a split torus and $O(\G_{x_2})\neq 1$.

Since we have full information about the action of $L_i$ on the set of edges $E_X(x_i)$ coming out of $x_i$, $i=1,2$, let us first consider the action of $\G_{x_1}$ on $E_X(x_1)$.  Since $T_1$ fixes $E_X(x_1)$ pointwise, we are interested in the action of the projection of $\G_{x_1}$ on $M_1\langle t_1\rangle$.  Assume first that $\G_{x_1}$ acts transitively on $E_X(x_1)$.  Then in the quotient graph $S$,  $x_1$ has a unique neighbouring vertex $x_2$.  Consider the set of  neighbours of $x_2$ in $S$.  If $VS=\{x_1, x_2\}$, the quotient graph corresponds to the edge-transitive lattice, and we are done.  Therefore, there exists $x_2 \in VS$ which has at least two neighbouring vertices in $S$.  Let $z_1, \ldots , z_k$ be representatives of the neighbouring vertices of $x_2$ in $S$ other than $x_1$.  If each $\G_{z_i}$ acts transitively on the set of edges $E_X(z_i)$ coming out of $z_i$, then we are looking at the whole  $VS=\{z_1, \ldots , z_k, x_2, x_1\}$.  If for some $i$, $C_{\Gamma_{z_i}}(M_{z_i})\leq Q_i$,  then $|\Gamma_{z_i}|\leq(q+1)\delta$. Looking at $\frac{1}{|\G_{x_2}|}+\frac{1}{|\G_{z_i}|}$, we obtain that $\mu(\Gamma\bs G)\geq\mu(\Gamma_0\bs G)$, a contradiction.  Therefore for all $i$, $C_{\Gamma_{z_i}}(M_{z_i})\not\leq Q_i$.  Hence there exists $y_{z_i}$, an element of $\Gamma_{z_i}$, such that $\langle y_{z_i}\rangle=O(C_{\Gamma_{z_i}}(M_{z_i}))$.  Then just like $y_1$, $y_{z_i}$ acts faithfully on $M_2$ and $y_{z_i}\in\Gamma_{x_2}$.  But $\Gamma_{x_2}$ is a subgroup of the normaliser of the split torus of $M_2$, and so $[y_1, y_{z_i}]=1$.  Again using the fact that we know how $L_2$ acts on the set $E_X(x_2)$, we observe that $y_1$ must fix the edge $[x_2, z_i]$.  It follows that $y_1$ fixes all vertices of $S$, and so is in the kernel of the action of $\Gamma$, a contradiction. Therefore, either $\Gamma_{x_1}$, or some $\Gamma_{z_i}$ does not act transitively on the corresponding set of edges. Thus by taking $x_1=z_i$ if necessary, we reduce to the following case: $\Gamma_{x_1}$ has at least two orbits on $E_X(x_1)$.

Therefore either $x_1$ has a unique neighbouring vertex $x_2$ in $S$ but the number of edges between them in $S$ is greater than one, or there are at least two neighbouring vertices of $x_1$ in the quotient graph.  Let us discuss those two cases.  In the former one, if $|VS|=2$ then $\langle y_1\rangle\triangleleft\Gamma$, which is a contradiction since $\Gamma$ is a faithful lattice. It follows that if $x_1$ has a unique neighbouring vertex $x_2$ is $S$ then $|VS|>2$, and so $x_2$ has more than one neighbouring vertex in $S$. Let us denote these neighbours (other than $x_1$) by $z_1, \ldots , z_k$.  Again not to contradict the minimality of covolume of $\Gamma$, we may assume that $C_{\Gamma_{z_i}}(M_{z_i})\not\leq Q_i$. Using the same argument as above we obtain that $y_1\in\Gamma_{z_i}$.  Assume that $y_1$ acts on $ M_{z_i}$ as an element of $K_{z_i}$ where $K_{z_i}\in\{S_4, A_5\}$.  If $o(y_1)\leq 15$, let us evaluate the covolume of $\Gamma$. Let $a_2\in\Gamma_{x_2}$ be such that $\langle a_2\rangle=O(\Gamma_{x_2})$. Then $y_1\in\langle a_2\rangle$.  If $o(y_1)=o(a_2)$, $|\Gamma_{x_2}|\leq 30\delta$ which immediately contradicts the minimality of covolume of $\Gamma$ ($\frac{1}{30\delta}>\frac{2}{(q+1)\delta}$ for $q>59$).  Hence, $o(a_2)>o(y_1)$.  Now, consider the action of $\Gamma_{x_1}$ on $M_1$.  If it acts as a subgroup of $K_1$ where $K_1\in\{S_4, A_5\}$, then $|\Gamma_{x_1}|\leq 15\cdot 60\delta_0$, and arguing as in the previous case ($q\equiv 3\pmod 4$), we get $|\G_{x_2}|\leq 150\delta$ leading to the usual contradiction with the minimality of covolume of $\Gamma$ as long as $q\geq 257$.  If it acts as a subgroup of a normaliser of a non-split torus, then together with the action of $a_2$, we will get that $\Gamma_{x_1}$ will act as a subgroup containing $p$--elements, a contradiction. Thus it can only act on $M_1$ as a subgroup of a normaliser of a split torus, and so $|\Gamma_{x_1}|\leq o(y_1)2\frac{o(a_2)}{o(y_1)}\delta\leq (q-1)\delta$. Since $|\Gamma_{x_2}|\leq (q-1)\delta$, we again get a contradiction with the minimality of the covolume of $\Gamma$.  Now suppose that  $o(y_1)>15$, then there exists $y_1'$ such that $y_1'\in\langle y_1\rangle$ and $y_1'$ centralises all the $M_{z_i}$s. We may use it instead of $y_1$.  Therefore for all $i$, either $y_1$ centralises $M_{z_i}$ or acts on it as a normal subgroup of a normaliser of a split torus.  In both cases, $\langle y_1\rangle$ is normal in $\Gamma_{z_i}$, $i=1, \ldots , k$.  Continuing with this argument we obtain that $\langle y_1\rangle$ is a normal subgroup of $\Gamma$ which implies that $y_1$ is in the kernel of the action of $\Gamma$ on $X$, a contradiction. 

Therefore we may suppose that $x_1$ has more than one neighbouring vertex in $S$. One is $x_2$, and let $z$ be among the other neighbouring vertices  of $x_1$.  If $C_{\Gamma_{z}}(M_{z})\leq Q_z$, then using $|\Gamma_{x_2}|$ and $|\Gamma_z|$, we obtain a contradiction with the minimality of $\mu(\Gamma\bs G)$: $\mu(\G\bs G)\geq\frac{1}{|\Gamma_{x_2}|}+\frac{1}{|\Gamma_z|}\geq\frac{2}{(q+1)\delta}$.  Therefore $C_{\Gamma_z}(M_z)\not\leq Q_z$ and we may take $x_2=z$. Therefore, we may assume that $\Gamma_{x_2}\cap T_2\neq 1$. It follows that there exists $y_2\in O(\Gamma_2\cap T_2)$ such that  $1\neq \langle y_2\rangle=O(\Gamma_2\cap T_2)$.

As before notice that $y_2\in\Gamma_{x_1}$ and $\langle y_2\rangle $ acts faithfully on $M_1$ via inner automorphisms.  Also, the subgroup structure of $L_1$ implies that either $\Gamma_{x_1}$ acts on $M_1$ as a subgroup of $K_1$ where $K_1\in\{S_4, A_5\}$, or as a subgroup of a normaliser of a split torus of $M_1\langle t_1\rangle$.

Assume the former, i.e., let $\G_{x_1}$ act on $M_1$ as a subgroup of $K_1$.  If $o(y_1)=2$, then $|\G_{x_1}|\leq 2\cdot 60$ which as always implies that $\mu(\G\bs G)$ is not minimal.  Hence, $O(\G_{x_1}\cap T_1)\neq 1$ and we may choose $y_1$ to be such that $\langle y_1\rangle=O(\G_{x_1}\cap T_1)$.  As always notice that $y_1\in\G_{x_2}$ and acts faithfully on $M_2$ via inner automorphisms.  If $\Gamma_{x_2}$ acts on $M_2$ as a subgroup of $K_2$ where $K_2\in\{S_4, A_5\}$, then $\mu(\Gamma\bs G)\geq\frac{2}{5\cdot 60\delta_0}>\frac{2}{(q+1)\delta}=\mu(\Gamma_0\bs G)$ for $q\geq 300$, a contradiction.

Hence $\Gamma_{x_2}$ acts on $M_2$ as a subgroup of a normaliser of a split torus. Recall that $|VS|>2$ and $x_1$ has neighbouring vertices other than $x_2$ in $S$. Let us call them $v_1, \ldots , v_k$.  Arguing for $v_i$ as we did for $x_2$, we obtain that $\langle y_1\rangle$ acts faithfully on $M_{v_i}$ and that $\Gamma_{v_i}$ acts on $M_{v_i}$ as a normaliser of a non-split torus.  It follows that $\langle y_1\rangle$ is normal in each $\Gamma_{v_i}$.  Now let $z_1, \ldots , z_m$ be the neighbouring vertices of $x_2$ in $VS-\{x_1\}$. Let us consider $\Gamma_{z_i}=\Gamma\cap P_{z_i}$.  If $C_{\Gamma_{z_i}}(M_{z_i})\leq Q_{z_i}$ for some $i$, then there is at most one such vertex, for otherwise we would contradict the minimality of covolume of $\Gamma$.  If there exists such a vertex, then without loss of generality assume that  such $i=1$, i.e., $C_{\Gamma_{z_1}}(M_{z_1})\leq Q_{z_1}$, and for $i>1$, $C_{\Gamma_{z_i}}(M_{z_i})\not\leq Q_{z_i}$. Further, we may assume that $T\leq P_{z_1}$.  Then $y_1, y_2\in\Gamma_{z_1}$.  If $\Gamma_{z_1}$ acts on $M_{z_1}$ as a subgroup of $K_{z_1}$ where $K_{z_1}\in\{S_4, A_5 \}$, then $|\Gamma_{z_1}|\leq 60\delta_0$, which leads to the usual contradiction as long as $q>59$. Thus $\Gamma_{z_1}$ acts on $M_{z_1}$ as a subgroup of  a normaliser of a split torus, and so $\langle y_1\rangle\triangleleft \Gamma_{z_1}$.  Now, let us look at $\Gamma_{z_i}$ for $i\geq 2$ (or it is possible that the previous case does not happen, then we are looking at $\Gamma_{z_i}$ for $i\geq 1$).  It follows  that for each $i\geq 1$ (or strictly greater than 1 if $o(y_{z_1})\mid 2$) there exists $y_{z_i}$ such that $\langle y_{z_i}\rangle=O(C_{\Gamma_{z_i}}(M_{z_i}))$.

Since $y_{z_i}$ sits in the kernel of $L_{z_i}$ in its action on $E_X(z_i)$, $y_{z_i}$ fixes every vertex neighbouring $z_i$, and in particular fixes $x_2$. It follows that $y_{z_i}\in\Gamma_{x_2}$. On the other hand, by the usual argument, $y_2$ acts faithfully on $M_{z_i}$ implying $[y_2, y_{z_i}]=1$. Therefore, $\langle y_{z_i}\rangle\triangleleft\Gamma_{x_2}$.  Finally, as  $C_{\Gamma_{x_2}}(M_{2})$ stabilises every edge $[x_2, z_i]$, it follows that  $C_{\Gamma_{x_2}}(y_{z_i})\leq\Gamma_{z_i}$. It follows that $y_1\in\Gamma_{z_i}$ and so $\langle y_1, y_2\rangle\leq\Gamma_{z_i}$.  Suppose $y_1$ acts on $M_{z_i}$ as an element of a subgroup of $K_{z_i}$ with $K_{z_i}\in\{S_4, A_5\}$. By the argument identical to the one  for $q\equiv 3\pmod 4$, we obtain that it is not possible and that there exists a nontrivial element $y$ of odd order such that $\langle y\rangle\triangleleft\Gamma$, a contradiction.

We are now reduced to the case when $\G_{x_1}$ acts on $M_1\langle t_1\rangle$ as a subgroup of $N_{M_1\langle t_1\rangle}(M_1\langle t_1\rangle\cap T)$. Notice that because of symmetry between $x_1$ and $x_2$, to finish the analysis it remains to consider the case when $\G_{x_2}$ acts on $M_2\langle t_2\rangle$ as a subgroup of $N_{M_2\langle t_2\rangle}(M_2\langle t_2\rangle\cap T)$.  But in this case $\langle \G_{x_1}, \G_{x_2}\rangle \leq N$. Hence we may move to the next vertex $y$ on our graph. However, again, the only possible case will be $\G_y\leq N$, and so on and so forth. Therefore, in the end of this case, the only possible conclusion will be $\G\leq N$, which is a contradiction as $N$ is not a cocompact lattice of $G$, nor does it contain one.  This completes the proof of Lemma \ref{l:minimalcovnoncomcase}, and thus the proof of Theorem \ref{t:covolumes} when $G$ admits an edge-transitive lattice.
\end{proof}

\subsection{Non-edge transitive case}\label{s:minimality non-edge-transitive}

Let $G$ be as in Theorem~\ref{t:covolumes} above, and assume that $G$ does not admit any edge-transitive lattices.  In this section we compute the covolume of the lattice $\G'$ constructed in Section~\ref{s:construction} above, which is not edge-transitive, and then prove that for $q \geq 514$, the lattice $\G'$ is the cocompact lattice of minimal covolume in $G$.  The important difference with the edge-transitive case is that  $$q\equiv 1\pmod 4 \quad \mbox{while} \quad Q_i^2\not\leq Z(G).$$  \noindent Hence certain arguments that worked very well in the cases $q=2^a$ and $q\equiv 3\pmod 4$  and even for $q\equiv 1\pmod 4$ above, now fail and/or become considerably longer and more delicate.  Thus we produce them carefully and at length. 

From the construction of $\G'$ in Section \ref{s:construction} and the discussion of covolumes in Section~\ref{s:lattices}, it follows that the covolume of $\G'$ is given by  $$\mu(\Gamma'\bs G)=\frac{1}{|\mathcal{A}_{x_1}|}+\frac{1}{|\mathcal{A}_{x_2}|}=\frac{1}{|SH_1|}+\frac{1}{|SH_2|}.$$  Recall that $S\cap H_i=Z(M_i)\leq Q_0$ and $|H_i:Z(M_i)|=\frac{q+1}{2}$.  Hence \[|SH_i|=\frac{|S||H_i|}{|S\cap H_i|}=|S|\frac{q+1}{2}=|Z(G)||Q_0:(Q_0\cap Z(G))|\frac{q+1}{2}.\]  Since $ |Q_0:(Q_0\cap Z(G))|=2\delta$ where $\delta\in\{ 1,2\}$ and its precise value depends on $G$, we obtain that
\begin{equation}
\label{eq::volumeGamma0}
\mu(\Gamma'\bs G)=\frac{2}{\delta|Z(G)|(q+1)}\ \ \mbox{with} \ \delta\in\{1,2\}\ \
\mbox{depending\ on}\ \  G.
\end{equation}

Let us now discuss the issue of minimality of covolume of $\G'$.  Our goal now is to prove the following statement.

\begin{lemma}
\label{l:minimalcovnoncomcase1}
Let $G$ be as in Theorem~\ref{t:cocompact existence}, $q\geq 514$  and assume that $G$ does not admit an edge-transitive lattice. Assume further that Conjecture~\ref{c:cocompact} holds.  Then $\Gamma'$ is a cocompact lattice in $G$ of minimal covolume.
\end{lemma}

Assume that there is a cocompact lattice $\G$ in $G$ whose covolume $\mu(\G\bs G)$ is strictly smaller than the covolume of $\G'$ given above. 
% Let $Y \subset X$ be a connected fundamental domain for $\G$ and let $A$ be the graph $A= \G \bs X$.
%Since $\Gamma$ has at least two orbits of vertices,  $Y$ contains at least two vertices $x_1$ and $x_2$ (connected by at least one edge), such that without loss of generality
%$G_{x_i} \leq P_i$ for $i=1,2$. By the discussion in Section~\ref{ss:lattices} above, $$\mu(\Gamma\bs G) =\sum_{y\in VY} \frac{1}{|\Gamma_y|}\geq \frac{1}{|\Gamma_{x_1}|}+\frac{1}{|\Gamma_{x_2}|}\ \ $$
%Since $\Gamma$ is discrete, $|\Gamma_{x_i}|$ is finite.  But  $\Gamma$ is cocompact, and so Proposition~\ref{p:nop-elements} above implies that, in fact, we may suppose that $\Gamma_{x_i}$ is a subgroup of $L_i$ of order co-prime to $p$ (where $L_i$ is a Levi complement of the parabolic $P_i$, $i=1,2$).

\begin{Remark}\label{Remark::TintGamma} Notice  that  $T\leq P_1\cap P_2$ together with $\G\cap P_i=\G_{x_i}$ yields $\G_{x_1}\cap T=\G_{x_2}\cap T$. \end{Remark}

As in the edge-transitive case, our plan is to subdivide further into the subcases: $$ \bf{Subcase\ 1}:\ \it{For} \ i=1,2,  L_i/Z(L_i)\cong PSL_2(q),\ \mbox{and}\ \bf{Subcase\ 2}:\ \it{For}  \ i=1,2,\  L_i/Z(L_i)\cong PGL_2(q).$$  We consider Subcase 1 in Section \ref{s:non subcase 1} and Subcase 2 in Section \ref{s:non subcase 2}.  Again, to simplify arguments we assume that $G$ has trivial centre, that is, the finite group $Z(G)$ satisfies $|Z(G)| = 1$.

\subsection{Non-edge-transitive case, Subcase 1}\label{s:non subcase 1}

In this case  $L_i=M_i\circ T_i$, that is, $L_i$ is a central (commuting) product of $M_i$ and  $T_i=C_T(M_i)$.  Moreover, if an element of $T$ centralises a non-split torus of $M_i$, then from the structure of $M_i$ and $L_i$, it follows immediately that it centralises $M_i$.  Now, $Z(G)=1$ implies that $T_i\cap T_j=1$ and $T_i$ acts faithfully on $M_j$ with $\{i,j\}=\{1,2\}$.  Let us make a few more comments about the structure of the $L_i$s.  Recall that for a finite group $F$, $O_2(F)$ denotes the largest normal $2$--subgroup of $F$.  

Suppose first that $L_i=M_i\times T_i$. Assume that $Z(M_i)\neq 1$, i.e., $M_i\cong SL_2(q)$.  Then $1\neq Q_0\leq C_G(M_i)$ for $i=1,2$, and so $Q_0\leq C_G(\la M_1, M_2\ra)\leq Z(G)=1$, a contradiction.  Thus if $L_i=M_i\times T_i$,  $M_i\cong PSL_2(q)$. Moreover, as far as the value of our parameter $\delta$ is concerned, it follows immediately that $|T_i|$ is odd whenever $\delta=1$, and  $|T_i|$ is even whenever $\delta=2$.
(In particular, in the key example $G=PSL_2(\F_q((t)))$, $|T_i|=1$ and $\delta=1$.)  Furthermore, $T_i$ must act faithfully on $M_j$ and so $T_i$ must be isomorphic to a subgroup of $M_j$. It follows that $|T_i|$ divides $\frac{q-1}{2}$.  Suppose now that $Z(M_i)\neq 1$, $T_i\cap M_i\neq 1$, $4\mid |T_i|$ and $L_i=M_i\circ _{\la -I\ra }T_i$. In particular, $M_i\cong M_j\cong SL_2(q)$.  Choose an element $g_i\in T\cap M_i$ of order $(q-1)$. Then $\langle g_i^{\frac{q-1}{2}}\rangle=Z(M_i)$. Since $g_i\in T$, it follows that $g_i\in L_j$ and $ g_i^{\frac{q-1}{2}}$ must act faithfully on $M_j$. Thus  $O_2(\la g_i\ra)$ acts faithfully on $M_j$ via inner automorphisms, which is a contradiction since $O_2(\la g_i\ra)\cong C_{2-part\ of\ (q-1)}$ while $\Inn(SL_2(q))=PSL_2(q)$ does not contain such a subgroup.  Therefore, this case does not happen.  Hence $M_i\cong PSL_2(q)$ and $L_i\cong T_i\times M_i$.

If $|\G_{x_i}\cap T_i|\leq \delta$, then \begin{equation}\label{eq::*} \G_{x_i}/\G_{x_i}\cap T_i\cong\G_{x_i}T_i/T_i\leq M_iT_i/T_i\cong M_i\cong PSL_2(q). \end{equation}  Notice that $|\G_{x_i}\cap T_i|=2$ implies that $\delta=2$.  By Dickson's Theorem, it follows that $|\G_{x_i}|\leq\delta(q+1)$. Therefore if $|\G_{x_i}\cap T_i|\leq\delta$ for both $i=1,2$, it follows that $\mu(\G\bs G)\geq\mu(\G'\bs G)$, a contradiction (this is precisely the case in~\cite[Example (6.2)]{LW} implying the minimality of the lattice constructed there).  Hence, without loss of generality we may assume that there exists $1\neq y_1\in\G_{x_1}\cap T_1$ such that $\langle y_1\rangle=\G_{x_1}\cap T_1$ with $o(y_1)>\delta$.  Then $o(y_1)\mid \frac{q-1}{2}$ and $\langle y_1\rangle$ acts faithfully on $M_2$ via  inner automorphisms. Notice that if $\delta=1$, $o(y_1)\neq 1$ is odd, and so for $\delta\in\{ 1, 2\}$, $o(y_1)\geq  3$.  As noticed in Remark~\ref{Remark::TintGamma}, since  $y_1\in \G_{x_1}\cap T$, $y_1\in\G_{x_2}$. Thus $\G_{x_2}$ acts non-trivially on $M_2$.  Now Dickson's Theorem asserts that $\G_{x_2}$ must act on $M_2$ either as a subgroup of a normaliser of a split torus of $M_2$, or as a subgroup of $K_2$ with $K_2\in\{ S_4, A_5\}$ (notice that in this case $|o(y_1)|\leq 5$).

Assume first that $|\G_{x_2}\cap T_2|\leq\delta$.  Then (\ref{eq::*}) implies that  $\G_{x_2}$ is actually isomorphic to a subgroup of $PSL_2(q)\times C_{\delta}$.   If $\G_{x_2}$ acts on $M_2$ as a subgroup of $K_2$, then using the previous paragraph we obtain that $\mu(\G\bs G)\geq\frac{1}{5\cdot(q+1)}+\frac{1}{60\delta}>\frac{2}{\delta(q+1)}=\mu(\G'\bs G)$ for $q>107$. Since this obviously contradicts the minimality of covolume of $\G$, $\G_{x_2}$ must be acting on $M_2$ as a subgroup of a normaliser of a split torus of $M_2$.  It follows that $\langle y_1\rangle$ is normal in $\G_{x_2}$.

We are now interested in the action of $\G_{x_1}$ on $M_1$.  By abuse of notation, we identify $x_i$ with its image in the quotient graph $A = \G \bs X$ for $i = 1,2$.  Then in the graph $A$, $x_1$ is a neighbour of $x_2$.  If $|A|=2$, it follows immediately that $\langle y_1\rangle\triangleleft\G$, a contradiction. And so $|A|>2$.  Let $z_1,\ldots, z_k$ be representatives of the other neighbouring vertices of $x_2$ in $A$.  If  for some $i$, $|C_{\Gamma_{z_i}}(M_{z_i})|\leq\delta$,  $|\Gamma_{z_i}|\leq (q+1)\delta$.  Hence, $\mu(\Gamma\bs G)\geq\mu(\Gamma'\bs G)$, a contradiction. Therefore for all $i$,  $|C_{\Gamma_{z_i}}(M_{z_i})|>\delta$.  Denote by $y_{z_i}$ an element of $\Gamma_{z_i}$ such that $\langle y_{z_i}\rangle=C_{\Gamma_{z_i}}(M_{z_i})$. We may use exactly the same arguments for $y_{z_i}$ as we did for $y_1$.  Then  $o(y_{z_i})\geq 3$, just like $y_1$, $y_{z_i}$ acts faithfully on $M_2$ and $y_{z_i}\in\Gamma_{x_2}$. But  $\Gamma_{x_2}$ acts on $M_2$ as a subgroup of the normaliser of the split torus.  Hence, $[y_1, y_{z_i}]=1$. Using the fact that we know how $L_2$ acts on $E_X(x_2)$,  we observe that $y_1$ then must fix the edge $[x_2, z_i]$.  It follows that $y_1$ fixes $x_1, x_2, z_1, \ldots , z_k$ and that $y_1\in\Gamma_{z_i}$.

Assume that $y_1$ acts on $M_{z_i}$ as a non-trivial element of $K_{z_i}\in\{S_4, A_5\}$.  Suppose first that $o(y_1)=2^{\beta_1}\cdot 3^{\beta_2}\cdot 5^{\beta_3}$ where either $\beta_1\leq 2$ and $\beta_i\leq 1$ for $i=2,3$ (this corresponds to $K_{z_i}\cong A_5$), or $\beta_1\leq 3$, $\beta_2\leq 1$ while $\beta_3=0$ (this is when $K_{z_i}\cong S_4$).  Let us evaluate the covolume of $\Gamma$. Let $a_2\in\Gamma_{x_2}$ be such that $\langle a_2\rangle$ acts faithfully on $M_2$ and $\langle a_2\rangle\leq\Gamma_{x_2}\leq(\langle a_2\rangle\rtimes\langle s_2\rangle)(\G_{x_2}\cap T_2)$ where $\langle a_2\rangle\rtimes\langle s_2\rangle\cong D_{2\cdot o(a_2)}$, $|\G_{x_2}\cap T_2|\leq\delta$ and $\langle a_2\rangle\times (\G_{x_2}\cap T_2)=\Gamma_2\cap T$.  Then $a_2\in \Gamma_{x_1}$ and without loss of generality we may assume that $y_1\in\langle a_2\rangle$. If $o(y_1)=o(a_2)$, then $|\Gamma_{x_2}|\leq (2^2\cdot 3\cdot 5)\cdot 2\delta$ which immediately contradicts the minimality of covolume of $\Gamma$ ($\frac{1}{120\delta}\geq \frac{2}{(q+1)\delta}$ for $q\geq 240$).  Hence, $o(a_2)>o(y_1)$ and so $\Gamma_{x_1}$ acts on $M_1$ as either a subgroup of $K_1\in\{S_4, A_5\}$, or as a subgroup of $N_{M_1}(M_1\cap T)$.

In the former case $|\Gamma_{x_1}|\leq 2^{\beta_1}\cdot 3^{\beta_2}\cdot 5^{\beta_3}\cdot 60$.  In particular, $|\G_{x_1}\cap T|\leq 2^{\beta_1}\cdot 3^{\beta_2}\cdot 5^{\beta_3}\cdot 5$.  As we noticed earlier,  $\G_{x_1}\cap T=\G_{x_2}\cap T$, and so $|\G_{x_2}\cap T|\leq 2^{\beta_1}\cdot 3^{\beta_2}\cdot 5^{\beta_3+1}$.  Since $\G_{x_2}$ acts on $M_2$ as a subgroup of a normaliser of a split torus of $M_2$, $|\G_{x_2}|\leq (2^{\beta_1}\cdot 3^{\beta_2}\cdot 5^{\beta_3+1})\cdot 2$. Note that if $\beta_1=1$, $\delta=2$.  If follows that $\mu(\G\bs G)>\frac{1}{2^{\beta_1+2}\cdot 3^{\beta_2+1}\cdot 5^{\beta_3+1}}+ \frac{1}{2^{\beta_1+1}\cdot 3^{\beta_2}\cdot 5^{\beta_3+1}}\geq  \frac{2}{(q+1)\delta}$ for $q\geq 514$, again a contradiction.

In the latter case, $|\Gamma_{x_1}|\leq o(y_1)2\frac{o(a_2)\delta}{o(y_1)}\leq (q-1)\delta$. Since $|\Gamma_{x_2}|\leq (q-1)\delta$, we again get a contradiction with the minimality of covolume of $\Gamma$. Therefore either $o(y_1)$ is divisible by $2^{\beta_1}\cdot3^{\beta_2}\cdot 5^{\beta_3}$ with  either $\beta_1\geq 3$ or $\beta_i\geq 2$ for some $i=2,3$, or by $2^{\beta_1}\cdot3^{\beta_2}$ with $\beta_1\geq 4$ or $\beta_2\geq 2$, or $o(y_1)$ is divisible by $\alpha_1\neq 1$ with $(\alpha_1, 30)=1$. In all the cases  there exists $1\neq y_1'\in\langle y_1\rangle$ such that $[y_1',M_{z_i}]=1$ for all $i$'s. Moreover, $o(y_1')\geq 3$.  We may now replace $y_1$ by $y_1'$ if necessary in all the previous subgroups to obtain the following conclusion: either $y_1$ centralises $M_{z_i}$ or acts on it as a subgroup of a normaliser of a split torus of $M_{z_i}$ where $i=1, \ldots, k$. In both situations, $\langle y_1\rangle$ is normal in $\Gamma_{z_i}$ for $i=1,\ldots,k$.  If $x_1$ has no other neighbouring vertices than $x_2$ in $A$, we continue with the argument (i.e., next look at $y_1$ in the $\G$--stabilisers of the neighbouring vertices of the $y_{z_i}$s in $A$ and so on) only to conclude that $\langle y_1\rangle \triangleleft \G$, a contradiction.

Therefore, it is possible that $x_1$ has more than one neighbouring vertex in $A$.  One is $x_2$ and  let $z$ be among the other neighbouring vertices of $x_1$.  If $|\G_z\cap T_z|\leq\delta$, then using $|\Gamma_{x_2}|$ and $|\Gamma_z|$, we obtain a contradiction with the minimality of covolume of $\Gamma$.  Therefore $|\G_z\cap T_z|>\delta$ and we may take $x_2=z$.  Thus whether $x_1$ has one or more neighbouring vertices in $A$, we may assume that  $|\G_{x_2}\cap T_2|>\delta$. Hence, there exists $y_2\in \G_{x_2}\cap T_2$ with $o(y_2)>\delta$ and $\langle y_2\rangle=\G_{x_2}\cap T_2$.  As for $y_1$,  we notice that $o(y_2)\mid\frac{q-1}{2}$,  $o(y_2)\geq 3$,  $y_2\in \G_{x_1}$ and $\langle y_2\rangle$ acts faithfully on $M_1$ via inner automorphisms.  Now  Dickson's Theorem  allows us to conclude that  either $\G_{x_1}$ acts on $M_1$  as a subgroup of $K_1$ where $K_1\in\{S_4, A_5\}$ (in which case  $o(y_2)\leq 5$),  or $\G_{x_1}$ acts on $M_1$ as  a subgroup of $N_{M_1}(M_1\cap T)$.

Let us begin with the case when $\G_{x_1}$ acts on $M_1$ as a subgroup of $K_1$.  If $\G_{x_2}$ acts on $M_2$ as a subgroup of $K_2$, then $\mu(\G \bs G)\geq\frac{2}{60\cdot 5}>\frac{2}{(q+1)\delta}=\mu(\G'\bs G)$ for $q\geq 514$, a contradiction.  Hence, $\G_{x_2}$ acts on $M_2$ as a subgroup $N_{M_2}(M_2\cap T)$ and in particular, $\langle y_1\rangle \triangleleft\G_{x_2}$.  Again, if $|VA|=2$, $\langle y_1\rangle\triangleleft \G$, a clear contradiction.  Thus $|VA|>2$ and let  $v_1, \ldots , v_k$ be the neighbours of $x_1$ in $VA-\{x_2\}$. Since $y_1$ fixes every edge in $E_X(x_1)$, it follows that $y_1$ acts faithfully on $M_{v_i}$ and holding a discussion similar to the above one with $v_i$ in place of $x_2$, we may assume that $\G_{v_i}$ acts on $M_{v_i}$ as a subgroup of a normaliser of a split torus of $M_{v_i}$. It follows that $\langle y_1\rangle $ is normal in each $\G_{v_i}$.  Now let $z_1, \ldots, z_m$ be the neighbours of $x_2$ in $VA-\{x_1\}$.  Let us consider $\G_{z_i}=\G\cap P_{z_i}$.  If $|C_{\G_{z_i}}(M_{z_i})|\leq\delta$,  then there is at most one such vertex, otherwise we would contradict the minimality of covolume of $\Gamma$. Hence, we may assume that if it happens, $i=1$, i.e.,  $|C_{\G_{z_1}}(M_{z_1})|\leq\delta$. Then we may further assume that $T\leq P_{z_1}$.  Thus $y_1, y_2\in\Gamma_{z_1}$. If $\Gamma_{z_1}$ acts on $M_{z_1}$ as a subgroup of $K_{z_1}\in\{ S_4, A_5\}$, then $|\Gamma_{z_1}|\leq 60\delta$, which is a contradiction, as always ($\frac{1}{60\delta}\geq\frac{2}{(q+1)\delta}$ for $q>120$).  Hence, $\Gamma_{z_1}$ acts on $M_{z_1}$ as a subgroup of a normaliser of a split torus of $M_{z_1}$.  It follows that $\langle y_1\rangle$ is a normal subgroup of $\Gamma_{z_1}$.  Now for $i>1$, there exists $y_{z_i}\in  C_G(M_{z_i})$ whose order $o(y_{z_i})>\delta$ (and thus is at least $3$) and does divide $\frac{q-1}{2}$.  But this element sits in the kernel of action of $L_{z_i}$ on $E_X(z_i)$ and therefore, $y_{z_i}\in \Gamma_{x_2}$.  On the other hand by the usual argument, $y_2$ acts faithfully on $M_{z_i}$ and so $[y_2, y_{z_i}]=1$. It follows that $\langle y_{z_i}\rangle$  is normal in $\G_{x_2}$.  Finally, as $C_{\Gamma_{x_2}}(M_2)$ stabilises $(x_2,z_i)$, it follows that $C_{\G_{x_2}}(y_{z_i})\leq\G_{z_i}$.  It follows that $y_1\in\G_{z_i}$ and so $\langle y_1, y_2\rangle\leq \G_{z_i}$.  Assume that $y_1$ acts on $M_{z_i}$ as a subgroup of $K_{z_i}\in\{S_4, A_5\}$.  Using the same argument as before we obtain that  there exists $y_1'\in\langle y_1\rangle$ with $[y_1',M_{z_i}]=1$ for all $i>1$ and with $o(y_1')\geq 3$.  In this case we  will replace $y_1$ by $y_1'$ if necessary in all the previous subgroups to obtain the following conclusion: $\langle y_1\rangle$ is normal in $\G_v$ for all the vertices mentioned so far, i.e., $x_1, x_2, z_1,  \ldots , z_m, v_1, \ldots, v_k$.  By iterating this argument we may show that $\langle y_1\rangle\triangleleft\G$ which is a contradiction.

We are now reduced to the last possible situation: $\G_{x_1}$ acts on $M_1$ as a subgroup of $N_{M_1}(M_1\cap T)$.  Notice that because of the symmetry between $x_1$ and $x_2$ to finish the analysis it remains to consider the case when $\G_{x_2}$ acts on $M_2$ as a subgroup of $N_{M_2}(M_2\cap T)$.  But in this case $\langle \G_{x_1}, \G_{x_2}\rangle \leq N$. Hence we may move to the next vertex $y$ on our graph. Using the previous argument we obtain that again that  the only possible case will be $\G_y\leq N$, and so on and so forth. Therefore, in the end of this case, the only possible conclusion will be $\G\leq N$, which is a contradiction as $N$ is not a cocompact lattice  of $G$, not does it contain any cocompact lattice.

\subsection{Non-edge-transitive case, Subcase 2}\label{s:non subcase 2} We are now in the situation when $T$ induces some non-trivial outer-diagonal automorphisms on $M_i$, that is, $L_i/Z(L_i)\cong PGL_2(q)$.  Consider $L_i=M_iT$. As before $M_i\triangleleft L_i$ and $T_i=C_T(M_i)$. Then there exists an element $t_i\in T-T_iM_i$  such that $t_i^2\in T_iM_i$ and $t_i$ induces an outer diagonal automorphism on $M_i$.  Since $q\equiv 1\pmod 4$, if $x$ is an involution in $L_i\cap T$, then $x\in M_iT_i$.

Recall that  $G$ does not admit any edge-transitive lattice.  Therefore if $Q_i\in\mathcal{S}yl_2(L_i)$ and $Q_i^2$ is its unique subgroup of index $2$, then $Q_i^2\not\leq Z(G)$.  It follows that $|Q_i/Q_i\cap Z(G)|\geq 4$ and so $\delta=2$.

The minimality of covolume of $\G'$ can be now shown by repeating exactly the same sequence of arguments as in Subcase 1 applied to subgroups of $L_i$, $i=1,2$. It turns out that the difference in the structure of $L_i$ (which is now a quotient of $GL_2(q)$) does not significantly affect the argument, and so we omit it here in order to avoid a fairly routine repetition.

This completes the proof of Theorem \ref{t:covolumes}.

\end{document}